\documentclass[11pt]{article}
\usepackage{color,verbatim,graphicx,color,ifthen}
\usepackage{amsmath,amsthm,amssymb,amsfonts,accents}
\DeclareMathAlphabet\mathbfcal{OMS}{cmsy}{b}{n}
\newboolean{UnBlinded}
\setboolean{UnBlinded}{true}
\makeatletter
\def\widebreve{\mathpalette\wide@breve}
\def\wide@breve#1#2{\sbox\z@{$#1#2$}%
     \mathop{\vbox{\m@th\ialign{##\crcr
\kern0.08em\brevefill#1{0.8\wd\z@}\crcr\noalign{\nointerlineskip}%
                    $\hss#1#2\hss$\crcr}}}\limits}
\def\brevefill#1#2{$\m@th\sbox\tw@{$#1($}%
  \hss\resizebox{#2}{\wd\tw@}{\rotatebox[origin=c]{90}{\upshape(}}\hss$}
\makeatletter
\def\QuasiTheorem{Result}
\def\COMPQUANT{\bC}
\def\bUall{\bU_{\mbox{\tiny all}}}
\def\bGamma{\boldsymbol{\Gamma}}
\def\bOmega{\boldsymbol{\Omega}}

\def\bb{\boldsymbol{b}}
\def\bw{\boldsymbol{w}}
\def\diag{\mbox{diag}}
\def\vech{\mbox{\rm vech}}
\def\bG{\boldsymbol{G}}

\def\Var{\mbox{Var}}
\def\betahat{{\widehat\beta}}
\def\BinRViidj{B_{ii'j}}
\def\bQmmd{\bQ_{mm'}}
\def\bigVertF{\big\Vert_{\mbox{\tiny $F$}}}
\def\VertF{\Vert_{\mbox{\tiny $F$}}}
\def\BigVertF{\Big\Vert_{\mbox{\tiny $F$}}}
\def\smbbeta{\mbox{$\tiny{\bbeta}$}}
\def\smbpsi{\mbox{$\tiny{\bpsi}$}}
\def\bbetahat{{\widehat\bbeta}}
\def\bpsihat{{\widehat\bpsi}}
\def\bone{\boldsymbol{1}}
\def\AsyCov{\mbox{Asy.Cov}}
\def\AsyVar{\mbox{Asy.Var}}
\def\trQuant{T}
\def\naturalNumbers{{\mathbb N}}
\def\bXAiid{\bX_{\mbox{\scriptsize A}ii'}}
\def\bXAiidu{\bX_{\mbox{\scriptsize A}i\idunder}}
\def\bXAiuid{\bX_{\mbox{\scriptsize A}\iunder i'}}
\def\tinybullet{\tiny\mbox{$\bullet$}}

\def\ndotdot{n_{\tinybullet\tinybullet}}
\def\iunder{{\undertilde{i}}}
\def\idunder{{\undertilde{i}}'}
\def\bXiid{\bX_{ii'}}
\def\bXiidu{\bX_{i\idunder}}
\def\bXiuid{\bX_{\iunder i'}}
\def\dAwindow{\dA^{\boxplus}}
\def\oonoo{\textstyle\frac{1}{n_{11}}}
\newcommand{\bXuptrione}{\uptri{\boldsymbol{X}}_1}
\def\bQtilde{{\widetilde \bQ}}
\def\bPsi{\boldsymbol{\Psi}}
\def\verytinyblackdiamond{\mbox{\fontsize{1.4mm}{1em}\selectfont{$\bigstar$}}}
\def\smallblackdiamond{\mbox{\fontsize{2.2mm}{1em}\selectfont{$\bigstar$}}}
\def\tinyblackdiamond{\mbox{\fontsize{2.2mm}{1em}\selectfont{$\bigstar$}}}
\def\smalluptri{\mbox{\fontsize{2.3mm}{1em}\selectfont{$\blacktriangle$}}}
\newcommand*{\uptri}[1]{\accentset{\stackrel{\mbox{$\smalluptri$}}{\null}}{#1}}

\newcommand{\bXuptrii}{\uptri{\boldsymbol{X}}_i}
\newcommand{\bXuptriiu}{\uptri{\boldsymbol{X}}_{\iunder}}
\newcommand{\bXuptriAi}{\uptri{\boldsymbol{X}}_{\mbox{\scriptsize A}i}}
\newcommand{\bXuptriAiu}{\uptri{\boldsymbol{X}}_{\mbox{\scriptsize A}\iunder}}
\newcommand*{\diamd}[1]{\accentset{\stackrel{\mbox{$\smallblackdiamond$}}{\null}}{#1}}
\newcommand{\bXdiamd}{\diamd{\boldsymbol{X}}}
\def\bAtilde{{\widetilde \bA}}
\def\bBtilde{{\widetilde \bB}}
\def\quarter{{\textstyle{1\over4}}}
\def\bXtilde{{\widetilde \bX}}
\def\bXdiamdrv{\bXdiamd_{\circ}}
\def\Xrv{X_{\circ}}
\def\smalldot{\mbox{\fontsize{0.1mm}{0.5em}\selectfont{$\bullet$}}}
\def\bXdot{\overset{\ \smalldot}{\bX}}
\newcommand*{\dotontop}[1]{\accentset{\stackrel{\mbox{$\smalldot$}}{\null}}{#1}}
\newcommand{\bXTdot}{\dotontop{\boldsymbol{X}}^T}
\newcommand*{\diamdontop}[1]{\accentset{\stackrel{\mbox{$\smallblackdiamond$}}{\null}}{#1}}
\newcommand*{\tinydiamdontop}[1]{\accentset{\stackrel{\mbox{$\tinyblackdiamond$}}{\null}}{#1}}

\newcommand{\ddiamd}{\tinydiamdontop{d}}
\newcommand{\bXTdiamdrv}{\diamdontop{\boldsymbol{X}_{\circ}}^T}
\newcommand{\ddiamdAsSubscript}{\verytinydiamdontop{d}}
\newcommand*{\verytinydiamdontop}[1]{\accentset{\stackrel{\mbox{$\verytinyblackdiamond$}}{\null}}{#1}}
\def\bXdotdot{\overset{\ \smalldot\smalldot}{\bX}}
\def\blockmatrixdum{\mathop{\mbox{\rm blockmatrix}}}
\def\blockmatrix#1{\blockmatrixdum_{#1}}
\def\bXA{\bX_{\mbox{\scriptsize A}}}
\def\bXB{\bX_{\mbox{\scriptsize B}}}
\def\dA{d_{\mbox{\tiny A}}}
\def\dB{d_{\mbox{\tiny B}}}
\def\bbetaA{\bbeta_{\mbox{\scriptsize A}}}
\def\bbetaB{\bbeta_{\mbox{\scriptsize B}}}
\def\bbetaAzero{\bbetaA^0}
\def\bbetaBzero{\bbetaB^0}
\def\bSigmaZero{\bSigma^0}
\def\bSigmadZero{(\bSigma')^0}
\def\sigsq{\sigma^2}
\def\sigsqZero{(\sigsq)^0}
\def\bbetaAMLE{{\widehat\bbeta}_{\mbox{\scriptsize A}}}
\def\bbetaBMLE{{\widehat\bbeta}_{\mbox{\scriptsize B}}}
\def\bSigmaMLE{{\widehat\bSigma}}
\def\bSigmadMLE{{\widehat\bSigma'}}
\def\bthetaMLE{{\widehat\btheta}}
\def\sigsqMLE{\sigmahat^2}
\def\bbetaAZero{\bbetaA^0}
\def\bbetaBZero{\bbetaB^0}
\def\bbetaZero{\bbeta^0}
\def\bpsiZero{\bpsi^0}
\def\bthetaZero{\btheta^0}

\definecolor{ruppertgreen}{rgb}{0,.7,.25}

\def\bbeta{\boldsymbol{\beta}}
\def\bpsi{\boldsymbol{\psi}}
\def\btheta{\boldsymbol{\theta}}
\def\bSigma{\boldsymbol{\Sigma}}
\def\sigmahat{{\widehat\sigma}}
\def\oom{{\textstyle\frac{1}{m}}}
\def\oomd{{\textstyle\frac{1}{m'}}}
\def\argmaxdum{\mathop{\mbox{\rm argmax}}}
\def\argmax#1{\argmaxdum_{#1}}
\def\sumim{\sum_{i=1}^m}
\def\sumidmd{\sum_{i'=1}^{m'}}
\def\bLambda{\boldsymbol{\Lambda}}
\def\bX{\boldsymbol{X}}
\def\ba{\boldsymbol{a}}
\def\banorm{\ba_{\mbox{\tiny norm}}}
\def\Nmart{N_{\mbox{\tiny mart}}}
\def\bn{\boldsymbol{n}}
\def\bs{\boldsymbol{s}}

\def\bx{\boldsymbol{x}}

\def\bz{\boldsymbol{z}}
\def\bA{\boldsymbol{A}}
\def\bD{\boldsymbol{D}}
\def\bW{\boldsymbol{W}}
\def\ellR{\ell_{\mbox{\tiny R}}}
\def\bwAr{\bw_{\mbox{\scriptsize A}r}}
\def\bwBr{\bw_{\mbox{\scriptsize B}r}}
\def\bwAone{\bw_{\mbox{\scriptsize A}1}}
\def\bwBone{\bw_{\mbox{\scriptsize B}1}}
\def\bwAdA{\bw_{\mbox{\scriptsize A}\dA}}
\def\bwBdB{\bw_{\mbox{\scriptsize B}\dB}}
\def\bVhalfLoose{\bV^{1/2}_{\mbox{\tiny loose}}}
\def\vecof{\mbox{\rm vec}}

\def\real{{\mathbb R}}
\def\convdist{\stackrel{D}{\longrightarrow}}
\def\convprob{\stackrel{P}{\longrightarrow}}
\def\bzero{\boldsymbol{0}}
\def\bI{\boldsymbol{I}}
\def\bV{\boldsymbol{V}}
\def\blockdiagdum{\mathop{\mbox{\rm blockdiag}}}
\def\blockdiag#1{\blockdiagdum_{#1}}
\def\subsubsubsection#1{\vskip2mm\noindent$\underline{\mbox{\em #1}}$\vskip2mm}
\def\smhalf{{\textstyle{\frac{1}{2}}}}
\def\tr{\mbox{tr}}
\def\myand{\&\ }
\def\bthetahat{{\widehat\btheta}}
\def\bXTArv{\bXArv^T}
\def\bXTBrv{\bXBrv^T}
\def\bLambdaNest{\bLambda_{\mbox{\scriptsize nest}}}
\def\bLambdaCross{\bLambda_{\mbox{\scriptsize cross}}}
\def\bJ{\boldsymbol{J}}
\def\bK{\boldsymbol{K}}
\def\bL{\boldsymbol{L}}
\def\bM{\boldsymbol{M}}
\def\bS{\boldsymbol{S}}
\def\bT{\boldsymbol{T}}
\def\bU{\boldsymbol{U}}
\def\bO{\boldsymbol{O}}
\def\bQ{\boldsymbol{Q}}
\def\be{\boldsymbol{e}}
\def\bB{\boldsymbol{B}}
\def\bC{\boldsymbol{C}}
\def\bR{\boldsymbol{R}}
\def\bXAi{\bX_{\mbox{\scriptsize A}i}}
\def\bXBi{\bX_{\mbox{\scriptsize B}i}}
\def\Cov{\mbox{\rm Cov}}
\def\Fsc{{\mathcal F}}
\def\Psc{{\mathcal P}}
\def\Xsc{{\mathcal X}}
\def\bX{\boldsymbol{X}}
\def\bY{\boldsymbol{Y}}
\def\bZ{\boldsymbol{Z}}
\def\bib{\vskip12pt\par\noindent\hangindent=1 true cm\hangafter=1}
\def\bXArv{\bX_{\mbox{\scriptsize A}\circ}}
\def\bXBrv{\bX_{\mbox{\scriptsize B}\circ}}
\def\bXrv{\bX_{\circ}}
\def\bXTrv{\bX^T_{\circ}}
\def\simind{\stackrel{{\tiny \mbox{ind.}}}{\sim}}
\def\stackdum{\mathop{\mbox{\rm stack}}}
\def\stack#1{\stackdum_{#1}}
\def\bXBiid{\bX_{\mbox{\scriptsize B}ii'}}

\def\iunder{{\undertilde{i}}}
\def\iunderTwo{\iunder^*}
\def\bbetaAMLET{({\widehat\bbeta}_{\mbox{\scriptsize A}})^T}
\def\bbetaBMLET{({\widehat\bbeta}_{\mbox{\scriptsize B}})^T}
\def\bXAione{\bX_{\mbox{\scriptsize A}i1}}
\def\bXAiuone{\bX_{\mbox{\scriptsize A}\iunder1}}
\def\bXAiidone{\bX_{\mbox{\scriptsize A}ii'1}}
\def\bXAiidj{\bX_{\mbox{\scriptsize A}ii'j}}
\def\bXAiidniid{\bX_{\mbox{\scriptsize A}ii'n_{ii'}}}
\def\bXBiidone{\bX_{\mbox{\scriptsize B}ii'1}}
\def\bXBiidj{\bX_{\mbox{\scriptsize B}ii'j}}
\def\bXBiidniid{\bX_{\mbox{\scriptsize B}ii'n_{ii'}}}
\def\bXAiuTwoOne{\bX_{\mbox{\scriptsize A}\iunderTwo1}}
\def\bXAoneone{\bX_{\mbox{\scriptsize A}\tiny{11}}}
\def\bXAtwoone{\bX_{\mbox{\scriptsize A}\tiny{21}}}
\def\gothicT{\mathfrak{T}}
\def\bIddiamnd{\bI_{\scriptsize{\ddiamdAsSubscript}}}
\def\bVnest{\bV_{\mbox{\scriptsize nest}}}
\def\bVcross{\bV_{\mbox{\scriptsize cross}}}
\def\lambdaNest{\lambda_{\mbox{\scriptsize nest}}}
\def\lambdaCross{\lambda_{\mbox{\scriptsize cross}}}
\def\Isc{{\mathcal I}}
\def\bLbreve{\widebreve{\bL}}
\def\bLbrevedash{\widebreve{\bL}'}
\def\bR{\boldsymbol{R}}
\def\bS{\boldsymbol{S}}
\def\TOPRIGHT{\bR_m}
\def\BOTRIGHT{\bS_m}
\newtheorem{result}{\textbf{\QuasiTheorem}}
\newtheorem{lemma}{\textbf{Lemma}}

\begin{document}
\vskip5mm
\centerline{\Large\bf Precise Asymptotics for Linear Mixed}
\vskip2mm
\centerline{\Large\bf Models with Crossed Random Effects}
\vskip5mm
\ifthenelse{\boolean{UnBlinded}}{
\centerline{\normalsize\sc Jiming Jiang$\null^1$, 
Matt P. Wand$\null^2$ and Swarnadip Ghosh$\null^3$}
\vskip5mm
\centerline{\textit{
$\null^1$University of California, Davis,
$\null^2$University of Technology Sydney
and $\null^3$Radix Trading}}}
{\null}

\vskip6mm
\centerline{11th September 2025}

\vskip8mm
\centerline{\large\bf Abstract}
\vskip3mm

We obtain an asymptotic normality result
that reveals the precise asymptotic behavior
of the maximum likelihood estimators of parameters for a very general
class of linear mixed models containing cross random effects. 
In achieving the result, we overcome theoretical
difficulties that arise from random effects being crossed as opposed 
to the simpler nested random effects case. Our new theory is for a class 
of Gaussian response linear mixed models which includes crossed
random slopes that partner arbitrary multivariate predictor effects
and does not require the cell counts to be balanced.
Statistical utilities include confidence interval construction, 
Wald hypothesis test and sample size calculations.

\noindent
\textit{Keywords:} Asymptotic normality, maximum likelihood estimation,
sample size calculations. 

\section{Introduction}\label{sec:introduction}

Linear mixed models with crossed random effects are useful for the analysis of 
regression-type data that are cross-classified according to two or more grouping
mechanisms. Baayen \textit{et al.} (2008), for example, use the terms \emph{subjects}
and \emph{items} for groupings that are typical in psychology studies.
Specific examples discussed in Baayen \textit{et al.} (2008) have subjects
corresponding to human participants in a psycholinguistic experiment and 
items corresponding to words in a particular language. 
Gao \myand Owen (2020) and Ghosh \textit{et al.} (2022)
is concerned with electronic commerce and related applications
involving crossed random effects, and is such that
subjects and items correspond to customers and products.

Despite the widespread use of linear mixed models with crossed random effects, theory
concerning the asymptotic behaviors of model parameter estimators is scant. 
This is largely due to the complicated mathematical forms that arise from
random effects being crossed. Unlike the nested random effects case, the marginal
covariance matrix of the response vector does not have a block diagonal form, which
makes theoretical analyses significantly more challenging. For Gaussian response
linear mixed models with nested random effects precise asymptotics are relatively
straightforward as conveyed by, for example, Section 3.5 of 
McCulloch \textit{et al.}\ (2008).
Recently Jiang \textit{et al.}\ (2022) obtained a precise asymptotic
normality result for the joint distribution of all model parameters in a
generalized linear mixed model with nested random effects. In this article
we derive an analogous result for Gaussian response linear mixed models
with crossed random effects.

Some early contributions to asymptotic theory for linear mixed models with crossed
random effects structures are Hartley \myand Rao (1967) and Miller (1977). 
Indeed, the second example in Section 4 of Miller (1977) corresponds to a special
case of the class of linear mixed models considered 
in the present article when his $c_{ij}$ term is omitted. Further details concerning
this example are in Sections 6.1 and 6.2 of Miller (1973), and includes an 
expression for the asymptotic covariance matrix of the maximum likelihood estimator
of the vector of variance parameters. Asymptotic normality of the maximum likelihood
estimators is also established in Miller (1973, 1977). However, the explicit results
in these seminal articles are confined to balanced linear mixed models that 
are devoid of predictor data. Jiang (1996) focused on restricted maximum likelihood
(REML) estimation of variance parameters in a wide class of linear mixed models that
include those containing crossed random effects and obtained conditions under
which asymptotic normality of the REML estimators hold. The results in Jiang (1996)
are expressed in terms of generic Fisher information matrices rather than the
explicit asymptotic forms provided by Jiang \textit{et al.} (2022).
Lyu \textit{et al.} (2024) is a recent article
that is also concerned asymptotic normality of estimators 
in a crossed random effects setting. Connections between
Lyu \textit{et al.} (2024) and this paper are described below.

In this article we obtain precise asymptotics, in a similar vein to those
of Jiang \textit{et al.} (2022), for Gaussian response linear mixed models 
with crossed random effects. Our results apply to a wide class of situations 
that include unbalanced designs, predictor data and multivariate crossed random 
effects. They reveal that asymptotic covariance matrices of the estimators 
parameter vectors are quite similar to those that arise for nested random
effects despite inherent differences due to effects being crossed.
For example, the estimates of fixed effects parameters that are unaccompanied
by random effects have the same asymptotic variances regardless of whether
the model contains nested or crossed random effects. However, as we shall see,
the pathway towards establishing such results for the crossed random
effects case is much longer and involved.

The majority of the research in this article was done 
concurrently with and independently of the
Lyu \textit{et al.} (2024) research and we became
aware of their article after devising \QuasiTheorem\ \ref{res:mainResult}.
The linear mixed model treated by Lyu \textit{et al.} (2024) does not assume that the 
responses are Gaussian. They also include a random interaction term, which our 
model does not have. In the case of Gaussian responses and additive crossed random 
effects, our main result extends the theoretical
findings of Lyu \textit{et al.} (2024) in the following two ways: 
(1) multivariate random slopes are included and (2) unbalanced 
cell counts are accommodated. Each of (1) and (2) are quite
important in practice, but require lengthy matrix algebraic
and convergence in probability arguments since the 
deterministic Kronecker product forms used in 
Miller (1973) and  Lyu \textit{et al.} (2024) no longer apply.

Contemporary data sets for which linear mixed models with crossed random effects
provide a useful vehicle for analysis vastly differ in terms of the density of the
observations. For some applications, the cell counts arising from subject/item
cross-classification are all non-zero. As an example, the illustration
given in Section 6 of Menictas \textit{et al.}\ (2023) for the U.S. National
Education Longitudinal Study has $8,564\times 24=205,488$ cells with a 
few observations per cell. 
The rows correspond to $8,564$ U.S. school students and the columns
correspond to $24$ items such as reading, mathematics
and science ability. The responses correspond to the scores for
each student/item combination. The students were followed 
longitudinally, which resulted in higher cell counts. 
Predictor data such as gender, time spent on homework 
and parental education were also recorded.
Other data sets, such as those that motivate
Ghosh \textit{et al.} (2022), have total number of observations much lower
than the number of cells. Ghosh \textit{et al.} describe an example concerning
customer ratings from the clothing company Stitch Fix with
$762,752\times 6,318$ cells. The rows correspond to $762,752$
customers and the columns correspond to $6,318$ clothing items.
There are five million ratings, which means that the average 
cell count is approximately $0.001$.
In this article we focus on dense data 
situations where the cell counts are non-zero and growing in our asymptotic
analyses. Relaxation to various sparse data situations is certainly of interest
but, with conciseness and closure in mind, this is left aside in this article's
theoretical study.

Generalized linear mixed models with crossed random effects
are particularly challenging theoretically and it was not until Jiang (2013)
that a consistency proof was established. As pointed out at the end of Section 4.5.7 
of Jiang \myand Nguyen (2021), there is no existing asymptotic distribution
theory for maximum likelihood estimators in the non-Gaussian version of such
models. We only treat the Gaussian version here.

The linear mixed model with crossed random effects that we study is described in Section 
\ref{sec:modelDescrip}, as well as maximum likelihood estimation of the model 
parameters. An asymptotic normality result that reveals the precise asymptotic 
behavior of all maximum likelihood estimators is given in Section \ref{sec:asyNorm}.
A key finding in Section \ref{sec:asyNorm} is that the leading terms are very 
similar to those arising in nested random effects models. In Section \ref{sec:heur} 
we provide some heuristic arguments that help explain these similarities.
Section \ref{sec:statsUtility} discusses statistical utility of the new theory.
Some concluding remarks are made in Section \ref{sec:concluding}.
An online supplement provides derivational details of the central result.

\section{Model Description and Maximum Likelihood Estimation}\label{sec:modelDescrip}

Consider the following crossed random effects linear mixed models:
\begin{equation}
\begin{array}{c}
\bY_{ii'}|\bU_i,\bU'_{i'},\bXAiid,\bXBiid\simind
N\big(\bXAiid(\bbetaAZero+\bU_i+\bU'_{i'})+\bXBiid\bbetaBZero,
\sigsqZero\bI\big),\\[2ex]
\bU_{i}  \simind N(\bzero, \bSigmaZero),\quad
1 \le i\le m, \quad \bU'_{i'}\simind N(\bzero, \bSigmadZero), \quad
1 \le i'\le m'
\end{array}
\label{eq:theModel}
\end{equation}
where here, and throughout this article, $\simind$ stands for ``independently distributed as''.

The dimensions of the matrices in (\ref{eq:theModel}) are:
\begin{eqnarray*}
\bY_{ii'}\ \mbox{is}\ n_{ii'}\times1,\ \bXAiid\ \mbox{is}\ n_{ii'}\times \dA,
\ \bbetaAZero\ \mbox{is}\ \dA\times1,\ \bU_i\ \mbox{is}\ \dA\times1,
\ \bU'_{i'}\ \mbox{is}\ \dA\times1\\[1ex]
\bXBiid\ \mbox{is}\ n_{ii'}\times \dB,\ \bbetaBZero\ \mbox{is}\ \dB\times1,
\ \bSigmaZero\ \mbox{is}\ \dA\times\dA\ \mbox{and}\ \bSigmadZero\ \mbox{is}\ \dA\times\dA.
\end{eqnarray*}
Here $n_{ii'}$ is the number of response measurements in the $(i,i')$th cell. If
$n_{ii'}=0$ then each of $\bY_{ii'}$, $\bXAiid$ and $\bXBiid$ are null.
The focus of this article is the precise asymptotic properties of the maximum likelihood
estimators of the model parameters when $m$, $m'$ and the $n_{ii'}$ all diverge to $\infty$.
Therefore, from now onwards, we assume that $n_{ii'}>0$ for all $1\le i\le m$ and $1\le i'\le m'$.

In (\ref{eq:theModel}), let the rows of $\bXAiid$ and $\bXBiid$ be
defined according to the notation
$$\bXAiid=\left[
\begin{array}{c}
\bXAiidone^T\\[1ex]
\vdots\\
\bXAiidniid^T
\end{array}
\right]
\quad\mbox{and}\quad
\bXBiid=\left[
\begin{array}{c}
\bXBiidone^T\\[1ex]
\vdots\\
\bXBiidniid^T
\end{array}
\right].
$$
We assume that the $\bXAiidj$, $1\le i\le m$, $1\le i'\le m'$, $1\le j\le n_{ii'}$
are independent and identically distributed $\dA\times1$
random vectors having the same distribution as $\bXArv$.
Similarly, the  $\bXBiidj$ over the same index set are 
independent and identically distributed $\dB\times1$
random vectors having the same distribution as $\bXArv$.

The following matrix assembly notation is useful for describing the maximum likelihood estimators
and their asymptotic properties. Firstly,
$$\stack{1 \le i \le d}(\bA_i)\equiv
\left[
\begin{array}{c}
\bA_1\\
\vdots\\
\bA_d
\end{array}
\right]
\quad\mbox{and}\quad
\blockdiag{1 \le i \le d}(\bA_i)\equiv
\left[
\begin{array}{cccc}
\bA_1 & \bO   & \cdots & \bO\\
\bO   & \bA_2 & \cdots & \bO\\
\vdots& \vdots& \ddots & \vdots\\
\bO   & \bO& \cdots &\bA_d 
\end{array}
\right]
$$
for matrices $\bA_1,\ldots,\bA_d$. The first of these definitions require 
that $\bA_i$, $1\le i\le d$, each have the same number of columns.
Next, define 
$$\blockmatrix{1\le i,\iunder\le d}(\bB_{i\iunder})\equiv
\left[
\begin{array}{ccc}
\bB_{11} & \cdots & \bB_{1d} \\
\vdots   & \ddots & \vdots   \\ 
\bB_{d1} & \cdots & \bB_{dd}
\end{array}
\right]
$$
for matrices $\bB_{i\iunder}$, $1\le i,\iunder\le d$, each having the same
numbers of rows and columns. If we then define 
$$\ndotdot\equiv{\displaystyle\sumim\sumidmd} n_{ii'},\quad
\bY\equiv\stack{1 \le i \le m}\left\{\stack{1 \le i' \le m'}(\bY_{ii'})\right\},$$
\begin{equation}
\bXA\equiv\stack{1 \le i \le m}\left\{\stack{1 \le i' \le m'}(\bXAiid)\right\}
\quad\mbox{and}\quad
\bXB\equiv\stack{1 \le i \le m}\left\{\stack{1 \le i' \le m'}(\bXBiid)\right\}
\label{eq:GrandPark}
\end{equation}
then standard manipulations show that 
$$\bY|\bXA,\bXB\sim N\Big(\bXA\bbetaAZero+\bXB\bbetaBZero,
\bV\big(\bSigmaZero,\bSigmadZero,\sigsqZero\big)\Big)$$
where
\begin{equation}
{\setlength\arraycolsep{1pt}
\begin{array}{rcl}
\bV(\bSigma,\bSigma',\sigsq)&\equiv&{\displaystyle\blockdiag{1\le i\le m}
\left\{\blockmatrix{1\le i',\idunder\le m'}(\bXAiid\bSigma\bXAiidu^T)\right\}}\\[2ex]
&&\quad+{\displaystyle\blockmatrix{1\le i,\iunder\le m}
\left\{\blockdiag{1\le i'\le m'}(\bXAiid\bSigma'\bXAiuid^T)\right\}}
+\sigsq\bI_{\ndotdot}.
\end{array}
}
\label{eq:Vmatrix}
\end{equation}
Therefore, the conditional log-likelihood is
\begin{equation}
{\setlength\arraycolsep{1pt}
\begin{array}{rcl}
\ell(\bbetaA,\bbetaB,\bSigma,\bSigma',\sigsq)&=&
-\smhalf\ndotdot\log(2\pi)-\smhalf\log\big|\bV(\bSigma,\bSigma',\sigsq)\big|\\[2ex]
&&\qquad-\smhalf(\bY-\bXA\bbetaA-\bXB\bbetaB)^T
\bV(\bSigma,\bSigma',\sigsq)^{-1}(\bY-\bXA\bbetaA-\bXB\bbetaB).
\end{array}
}
\label{eq:FLbillboards}
\end{equation}
The maximum likelihood estimator of 
$\big(\bbetaAzero,\bbetaBzero,\bSigmaZero,\bSigmadZero,\sigsqZero\big)$ is 
$$(\bbetaAMLE,\bbetaBMLE,\bSigmaMLE,\bSigmadMLE,\sigsqMLE)
\equiv\argmax{\bbetaA,\bbetaB,\bSigma,\bSigma',\sigsq}\ell(\bbetaA,\bbetaB,\bSigma,\bSigma',\sigsq).$$

\section{Asymptotic Normality \QuasiTheorem}\label{sec:asyNorm}

We now present the article's main centerpiece: an asymptotic normality result that reveals
the precise asymptotic behavior of the maximum likelihood
estimation of $(\bbetaAMLE,\bbetaBMLE,\bSigmaMLE,\bSigmadMLE,\sigsqMLE)$ for data corresponding
to (\ref{eq:theModel}).

Define 
$$n\equiv\frac{\ndotdot}{mm'}=\mbox{average of the within-cell sample sizes}$$
and
$$\COMPQUANT_{\bbetaB}\equiv\mbox{lower right $\dB\times\dB$ block of}\ 
\left\{E\big(\bXrv\bXrv^T\big)\right\}^{-1}
\ \mbox{where}\ 
\bXrv\equiv
\left[
\begin{array}{c}
\bXArv\\
\bXBrv
\end{array}
\right].
$$
Let $\bD_d$ denote the matrix of zeroes and ones such that $\bD_d\vech(\bA)=\vecof(\bA)$
for all $d\times d$ symmetric matrices $\bA$. The Moore-Penrose inverse of $\bD_d$ is
$\bD_d^+=(\bD_d^T\bD_d)^{-1}\bD_d^T$. 

The result relies on the following assumptions:

\begin{itemize}
\item[]
\begin{itemize}
\item[(A1)] The cell dimensions $m$ and $m'$ diverge to $\infty$ in such a way that 
$m=O(m')$ and $m'=O(m)$.
\item[(A2)] The within-cell sample sizes $n_{ii'}$ diverge to $\infty$ in such a way that
$${\displaystyle\max_{1\le i\leq m, 1\le i'\le m'}}\big|n_{ii'}/n - C_{ii'}\big| 
\rightarrow 0\quad\mbox{as}\ m,m' \rightarrow\infty
$$
for positive constants $C_{ii'}$, $1\le i\le m$, $1\le i'\le m'$,
that are bounded above and away from zero. Also, $n/m\to0$ as $m$ and $n$ diverge.
\item[(A3)] All entries of both $\bXArv$ and $\bXBrv$ are not degenerate at zero and
have finite second moment.
\end{itemize}
\end{itemize}

\begin{result}
Assume that (A1)--(A3) and some additional regularity conditions hold. Then
$$
\left[
\begin{array}{c}
\left\{\displaystyle{\frac{\bSigma^0}{m}+\frac{(\bSigma')^0}{m'}}\right\}^{-1/2}
\Big(\bbetaAMLE-\bbetaAZero\Big)\\[2ex]
\left\{\displaystyle{\frac{(\sigma^2)^0\COMPQUANT_{\bbetaB}}{mm'n}}\right\}^{-1/2}
\Big(\bbetaBMLE-\bbetaBZero\Big)\\[2ex]
\left\{\displaystyle{\frac{2\bD_{\dA}^+(\bSigmaZero\otimes\bSigmaZero)\bD_{\dA}^{+T}}{m}}\right\}^{-1/2}
\vech\big(\bSigmaMLE-\bSigmaZero\big)\\[3ex]
\left\{\displaystyle{\frac{2\bD_{\dA}^+(\bSigmadZero\otimes\bSigmadZero)\bD_{\dA}^{+T}}{m'}}\right\}^{-1/2}
\vech\big(\bSigmadMLE-\bSigmadZero\big)\\[3ex]
\left[\displaystyle{\frac{2\{\sigsqZero\}^2}{mm'n}}\right]^{-1/2}\big\{\sigsqMLE-\sigsqZero\big\}
\end{array}
\right]
\convdist N(\bzero,\bI).
$$
\label{res:mainResult}
\end{result}

\vskip3mm\par
Some remarks concerning \QuasiTheorem\ \ref{res:mainResult} are:
\begin{itemize}
\item[1.] \QuasiTheorem\ \ref{res:mainResult} provides 
following asymptotic covariance matrices
of the maximum likelihood estimators:
$$\AsyCov(\bbetaAMLE)=\frac{\bSigma^0}{m}+\frac{(\bSigma')^0}{m'},\ \ 
\AsyCov(\bbetaBMLE)=\frac{\sigsqZero\COMPQUANT_{\bbetaB}}{mm'n},
\ \ \AsyCov(\bSigmaMLE)=\frac{2\bD_{\dA}^+(\bSigmaZero\otimes\bSigmaZero)\bD_{\dA}^{+T}}{m},$$
$$\AsyCov(\bSigmadMLE)=\frac{2\bD_{\dA}^+(\bSigmadZero\otimes\bSigmadZero)\bD_{\dA}^{+T}}{m'}
\quad\mbox{and}\quad\AsyVar(\sigsqMLE)=\frac{2\{\sigsqZero\}^2}{mm'n}.$$
Notation such as $\AsyVar(\sigsqMLE)$ is based on the fact that, for large $m$, $m'$ and $n$,
$\sigsqMLE$ has an approximate Normal distribution with mean $\sigsqZero$ and variance $\AsyVar(\sigsqMLE)$.
There are marked differences in the rates of convergence. For example, the entries of $\bbetaAMLE$
have order $m^{-1}$ asymptotic variances, whilst those of $\bbetaBMLE$ have order $(mm'n)^{-1}$
asymptotic variances. Note that $\bbetaAZero$ and $\bbetaBZero$ differ in that the former is
partnered by crossed random effects in (\ref{eq:theModel}).
\item[2.] The asymptotic normality results for $\bSigmaMLE$ and $\bSigmadMLE$ can be
converted to forms that are more amenable to interpretation and confidence
interval construction using the Multivariate Delta Method (e.g. Agresti, 2013,
Section 16.1.3). For example, if $\dA=2$ and the entries of $\bSigma$ are
parameterized as 
$$\bSigma=\left[
\begin{array}{cc}
\sigma_1^2           & \rho\,\sigma_1\sigma_2\\[1ex]
\rho\,\sigma_1\sigma_2 & \sigma_2^2
\end{array}
\right]
$$
then \QuasiTheorem\ \ref{res:mainResult} implies the following asymptotic normality
results for standard transformations of the first
standard deviation parameter and correlation parameter:
$$\sqrt{m}\big\{\log(\sigmahat_1)-\log(\sigma_1^0)\big\}\convdist N(0,\smhalf)
\quad\mbox{and}\quad
\sqrt{m}\big\{\tanh^{-1}(\widehat{\rho})-\tanh^{-1}(\rho^0)\big\}\convdist N(0,1).$$
Analogous results hold for $\sigmahat_2$ and $\bSigmadMLE$.
\item[3.] There is asymptotic orthogonality between each pair of random vectors within the set
$$\big\{\bbetaAMLE,\bbetaBMLE,\vech(\bSigmaMLE),\vech(\bSigmadMLE),\sigsqMLE\big\}.$$
\item[4.]
Outside of \QuasiTheorem\ \ref{res:mainResult} and Lyu \textit{et al.} (2024),
we are not aware of results for linear mixed models with crossed 
random effects that provide the precise asymptotic covariances given by 
\QuasiTheorem\ \ref{res:mainResult}
for estimation of fixed effects, even for simplified versions of (\ref{eq:theModel}) such as 
those having $\bXAiid=\bone_{n_{ii'}}$ and $\bXBiid$ null. 
In this special case, in which the only fixed effect is the intercept parameter,
the $\bSigma/m + \bSigma'/m'$ leading term behaviour is also apparent from
Theorem 1 of Lyu \textit{et al.} (2024) when their variable $\eta$ is in the
interior of the positive half-line. The predictor set-ups differ between 
the two articles, which hinders succinct comparison of the fixed effects
results for more general cases.
\item[5.] \QuasiTheorem\ \ref{res:mainResult} extends the results of Miller (1973) 
and Lyu \textit{et al.} (2024), concerning asymptotic distributions of 
variance component estimators, to covariance matrices of arbitrary dimension.
\item[6.] Under (A1) $m$ and $m'$ diverge to $\infty$ at the same rate. In some circumstances
this assumption may not be realistic and other assumptions concerning $m$ and $m'$ divergence
may be more appropriate. The subsequent modification of 
\QuasiTheorem\ \ref{res:mainResult} is straightforward.
For example, if $m'=o(m)$ then the component concerning $\bbetaAMLE$ becomes
$$\left\{\displaystyle{\frac{(\bSigma')^0}{m'}}\right\}^{-1/2}
\Big(\bbetaAMLE-\bbetaAZero\Big)\convdist N(\bzero,\bI)
\quad\mbox{leading to}\quad \AsyCov(\bbetaAMLE)=\frac{(\bSigma')^0}{m'}.
$$
\item[7.] The asymptotic covariances for linear mixed models
with crossed random effects have forms that are very similar to those with two-level 
nested random effects. See, for example, the Gaussian 
special case of Theorem 1 of Jiang \textit{et al.} (2022). 
At first glance, this result is somewhat surprising 
and intriguing since the two types of linear mixed models have fundamental differences.
In Section \ref{sec:heur} we provide some heuristic arguments that help explain
this interesting phenomenon.
\item[8.] For the special case $\bXArv=1$ and $\bXBrv=\Xrv$ we have
$$\AsyVar(\betahat_B)=\frac{\sigsqZero}{\Var(\Xrv)(\mbox{total sample size})}.$$
This matches the well-known expression for the asymptotic variance of the slope parameter
in the simple linear regression model. Analogous results arise when $\bXBrv$ is
multivariate. Despite the presence of crossed random effects, the asymptotic 
behaviors of the estimators of slope parameters that are unaccompanied
by random effects are the same as in the ordinary multiple regression situation.
The heuristics in Section \ref{sec:heur} provide some insight into
this phenomenon.
\item[9.] The presence of multivariate random slopes in  
the crossed random effects model (\ref{eq:theModel}) leads to 
considerable challenges in the establishment of the 
\QuasiTheorem\ \ref{res:mainResult} precise asymptotic normality statement.
Detailed and delicate arguments, not given here, would be required
to obtain sufficient regularity conditions under
which \QuasiTheorem\ \ref{res:mainResult} holds.
\item[10.] \emph{Restricted} maximum likelihood estimation is
a commonly used alternative to maximum likelihood estimation
in linear mixed models-based analyses. For model (\ref{eq:theModel}), it involves
replacement of (\ref{eq:FLbillboards}) by the restricted log-likelihood
$$\ellR(\bbetaA,\bbetaB,\bSigma,\bSigma',\sigsq)\equiv
\ell(\bbetaA,\bbetaB,\bSigma,\bSigma',\sigsq)
-\smhalf\log\big|[\bXA\ \bXB]^T\bV(\bSigma,\bSigma',\sigsq)^{-1}
[\bXA\ \bXB]\big|.
$$
The extra term invokes a finite sample adjustment to
the estimators. Result \ref{res:mainResult}, which is concerned with
large sample behavior, also applies to the 
restricted maximum likelihood estimators
of the parameters in (\ref{eq:theModel}).
\item[11.] The establishment of \QuasiTheorem\ \ref{res:mainResult} requires
complicated and long-winded arguments, and are deferred to 
an online supplement.
\end{itemize}

\section{Heuristics on Nested/Crossed Asymptotics Similarities}\label{sec:heur}

We now address the fact that the asymptotic covariance expressions
in \QuasiTheorem\ \ref{res:mainResult} are quite similar to those arising
in the two-level nested case. This involves heuristic arguments
that show that the fixed effects maximum likelihood estimators admit 
quite similar forms when sample means are replaced by population
means. Throughout this section we write $\bbeta$ rather than
$\bbetaZero$. A similar convention is used for $\bSigma$, $\bSigma'$ and $\sigma^2$. 
This suppression of the ``true value'' notation is to aid exposition.

Gaussian response linear mixed models have the following general form:
\begin{equation}
\bY|\bU\sim N\big(\bX\bbeta+\bZ\bU,\bR\big),\quad
\bU\sim N(\bzero,\bG).
\label{eq:LMMgenForm}
\end{equation}
For the crossed random effects model (\ref{eq:theModel})
$$\bX=[\bX_A\ \bX_B],\quad \bZ=
\left[\blockdiag{1\le i\le m}\left\{\stack{1\le i'\le m'}(\bXAiid)\right\}
\ \ \stack{1\le i\le m}\left\{\blockdiag{1\le i'\le m'}(\bXAiid)\right\}\right],
$$
$$
\bG=\mbox{blockdiag}\big(\bI_m\otimes\bSigma,\bI_{m'}\otimes\bSigma'\big)
\quad\mbox{and}\quad\bR=\sigma^2\bI
$$
where $\bX_A$ and $\bX_B$ are given by (\ref{eq:GrandPark}).

The Gaussian version of the class of \emph{nested} linear mixed
models studied by Jiang \textit{et al.} (2022) is
\begin{equation}
\begin{array}{c}
\bY_i|\bU_i,\bXAi,\bXBi\simind
N\big(\bXAi(\bbetaA+\bU_i)+\bXBi\bbetaB,\sigma^2\bI\big),\\[2ex]
\bU_{i}  \simind N(\bzero, \bSigma),\quad
1 \le i \le m,
\end{array}
\label{eq:theNestModel}
\end{equation}
which is a special case of (\ref{eq:LMMgenForm}) with 
$$\bX=\stack{1\le i\le m}[\bXAi\ \bXBi],\quad 
\bZ=\blockdiag{1\le i\le m}(\bXAi),\quad 
\bG=\bI_m\otimes\bSigma
\quad\mbox{and}\quad\bR=\sigma^2\bI.
$$
Analogous to the set-up for model (\ref{eq:theModel}), we assume
that the transposes of the rows of $\bXAi$, $1\le i\le m$, are independent and 
identically distributed $\dA\times1$
random vectors having the same distribution as $\bXArv$.
A similar assumption applies to the $\bXBi$.

In terms of the notation in (\ref{eq:LMMgenForm}), the fixed effects maximum
likelihood estimator has the following generalized least squares form:
$$\bbetahat=(\bX^T\bV^{-1}\bX)^{-1}\bX^T\bV^{-1}\bY
\quad\mbox{where}\quad
\bV\equiv\bZ\bG\bZ^T+\bR.
$$
If $\Xsc$ denotes the predictor data in the $\bX$ and $\bZ$ matrices then
the conditional covariance matrix of the fixed effects estimator is
$$\Cov(\bbetahat|\Xsc)=\big(\bX^T\bV^{-1}\bX\big)^{-1}.$$

For the remainder of this section we assume that the data are balanced.
In the crossed case this corresponds to $n_{ii'}=n$ for all $1\le i\le m$
and $1\le i'\le m'$. For the nested case $n_i=n$ for all $1\le i\le m$.

\subsection{The $\bX=\bone$ Special Case}\label{sec:CaliforniaFaded}

Consider the following special case of (\ref{eq:LMMgenForm}):
$$\bY|\bU\sim N\big(\bone\beta_0+\bZ\bU,\bR\big),\quad
\bU\sim N(\bzero,\bG).
$$
for which $\bX=\bone$, which is such that the only fixed effect 
effect is the intercept parameter $\beta_0$.

A further simplification is
\begin{equation}
\bZ=\left\{
\begin{array}{cl}
\big[\bI_m\otimes\bone_{m'n}\ \ \bone_m\otimes\bI_{m'}\otimes\bone_n\big]
& \mbox{for the crossed case,}\\[1ex]
\bI_m\otimes \bone_n & \mbox{for the nested case,}
\end{array}
\right.
\label{eq:orangeCars}
\end{equation}
which corresponds to the random intercept-only models.
Let $\bVcross$ and $\bVnest$ respectively denote the
$\bV$ matrix for the crossed and nested cases based 
on the versions of $\bZ$ given in (\ref{eq:orangeCars}).
Bringing in the commonly used notation $\bJ_d\equiv\bone_d\bone_d^T$
we then have
$$\bVcross=\Sigma(\bI_m\otimes\bJ_{m'n})+\Sigma'(\bJ_m\otimes\bI_{m'}\otimes\bJ_n)
+\sigma^2\bI_{mm'n}
\quad\mbox{and}\quad
\bVnest=\Sigma(\bI_m\otimes\bJ_n)+\sigma^2\bI_{mn}$$
where $\Sigma\equiv\bSigma$ and $\Sigma'\equiv\bSigma'$ are scalars
in the current random intercept special cases. 
The following results are key:
\begin{equation}
\bVcross\bone=\lambdaCross\bone\quad\mbox{and}\quad
\bVnest\bone=\lambdaNest\bone,
\label{eq:BeanieF}
\end{equation}
where $\bone$ denotes a vector of ones with appropriate size,
\begin{equation}
\lambdaCross\equiv\Sigma m'n+\Sigma' mn+\sigma^2\quad\mbox{and}
\quad\lambdaNest\equiv n\Sigma+\sigma^2.
\label{eq:PalmDrive}
\end{equation}
The fact that $\bone$ is an eigenvector of both $\bVcross$ and
$\bVnest$ leads the fixed effects estimators having simpler
and similar forms. A key step involves the inverse eigenvalue results
$$\bVcross^{-1}\bone=\big(1/\lambdaCross\big)\bone\quad\mbox{and}\quad
\bVnest^{-1}\bone=\big(1/\lambdaNest\big)\bone.
$$
We then obtain
$$\betahat_0=(\bone^T\bV^{-1}\bone)^{-1}\bone^T\bV^{-1}\bY
=(\bone^T\bone)^{-1}\bone^T\bY=\mbox{average of the response data}
$$
for both $\bV=\bVcross$ and $\bV=\bVnest$. 
We also have
\begin{equation}
\Var(\betahat_0)=\frac{\lambda}{\mbox{total sample size}}
\label{eq:SFCentral}
\end{equation}
where $\lambda=\lambdaCross$ in the crossed case and 
$\lambda=\lambdaNest$ in the nested case. Results (\ref{eq:PalmDrive})
and (\ref{eq:SFCentral}) then lead to the exact expressions
$$\Var(\betahat_0)=\left\{
\begin{array}{ll}
{\displaystyle\frac{\Sigma}{m}+\frac{\Sigma'}{m'}+\frac{\sigma^2}{mm'n}} 
& \mbox{in the crossed case},\\[1.5ex]
{\displaystyle\frac{\Sigma}{m}+\frac{\sigma^2}{mn}} & \mbox{in the nested case}\\
\end{array}
\right.
$$
which are in keeping with the leading term expression in (\ref{res:mainResult})
and the analogous result in Jiang \textit{et al.} (2022).

In this subsection, we have seen that the 
eigenvalue/eigenvector results given by (\ref{eq:BeanieF})
lead to the fixed effects estimator reducing to 
ordinary least squares form in both cases. Therefore, the $\beta_0$
estimators behave quite similarly despite the ostensible
differences between the crossed and nested cases.

\subsection{Heuristics for the General $\bX$ Crossed Case}

We commence by noting the following exact result:
\begin{eqnarray*}
&&\bVcross\bX=\Bigg[\stack{1\le i\le m}\left[\left\{\stack{1\le i'\le m'}(\bXAiid)\right\}
\left(\bSigma\sum_{i'=1}^{m'}\bXAiid^T\bXAiid 
+ \bSigma'\sum_{i=1}^{m}\bXAiid^T\bXAiid\right)\right]\\
&&\qquad\qquad
\stack{1\le i\le m}\left[\left\{\stack{1\le i'\le m'}(\bXAiid)\right\}
\left(\bSigma\sum_{i'=1}^{m'}\bXAiid^T\bXBiid 
+ \bSigma'\sum_{i=1}^{m}\bXAiid^T\bXBiid\right)\right]
\Bigg]+\sigma^2\bX.
\end{eqnarray*}
Then results such as 
$$\quad\frac{1}{mn}\sum_{i=1}^{m}\bXAiid^T\bXAiid\convprob E(\bXArv\bXTArv)
\ \mbox{and}\ \quad
\frac{1}{mn}\sum_{i=1}^{m}\bXAiid^T\bXBiid\convprob E(\bXArv\bXTBrv)
$$
for all $1\le i'\le m'$ lead to the approximation
$$\bVcross\bX\approx \bX\bLambdaCross$$
where
\begin{equation}
\bLambdaCross\equiv
\left[
\begin{array}{cc}
n\big(m'\bSigma+m\bSigma'\big)E(\bXArv\bXTArv)+\sigma^2\bI_{\dA}&
\quad 
n\big(m'\bSigma+m\bSigma'\big)E(\bXArv\bXTBrv)\\[1ex]
\bO & \sigma^2\bI_{\dB}
\end{array}
\right].
\label{eq:PaloHummus}
\end{equation}
We then have
\begin{equation}
\bbetahat\approx\bLambdaCross(\bX^T\bX)^{-1}\bLambdaCross^{-T}\bX^T\bY
\quad\mbox{and}\quad
\Cov(\bbetahat|\Xsc)\approx\bLambdaCross(\bX^T\bX)^{-1}.
\label{eq:CayugaHeights}
\end{equation}
A simple consequence of (\ref{eq:PaloHummus}) and (\ref{eq:CayugaHeights})
is
\begin{equation}
{\setlength\arraycolsep{1pt}
\begin{array}{rcl}
\Cov(\bbetaBMLE|\Xsc)&\approx&
\sigma^2\Big\{\mbox{lower right $\dB\times\dB$ block of\ } (\bX^T\bX)^{-1}\Big\}\\[2ex]
&\approx&\left(\displaystyle{\frac{\sigma^2}{mm'n}}\right)
\Big[\mbox{lower right $\dB\times\dB$ block of}\ 
\left\{E\big(\bXrv\bXrv^T\big)\right\}^{-1}\Big].\\[2ex]
&=&\left(\displaystyle{\frac{\sigma^2}{\mbox{total sample size}}}\right)
\Big[\mbox{lower right $\dB\times\dB$ block of}\ 
\left\{E\big(\bXrv\bXrv^T\big)\right\}^{-1}\Big].
\end{array}
}
\label{eq:raftTrip}
\end{equation}
The asymptotic covariance matrix of $\Cov(\bbetaAMLE|\Xsc)$ has 
a similar derivation based on (\ref{eq:PaloHummus}) and 
(\ref{eq:CayugaHeights}).

\subsection{Heuristics for the General $\bX$ Nested Case}

For the nested model (\ref{eq:theNestModel}) we have the exact expression
$$\bVnest\bX=\bX_A\bSigma
\Big(\stack{1\le i\le m}\big[\bXAi^T\bXAi\ \ \ \bXAi^T\bXBi\big]\Big)
+\sigma^2\bX.
$$
As $n\to\infty$ and for each $1\le i\le m$ we have 
$$\frac{1}{n}\bXAi^T\bXAi\convprob E(\bXArv\bXTArv)\quad\mbox{and}\quad
\frac{1}{n}\bXAi^T\bXBi\convprob E(\bXArv\bXTBrv)
\quad\mbox{as}\ n\to\infty.
$$
Therefore
$$\bVnest\bX\approx \bX\bLambdaNest\quad\mbox{where}\quad 
\bLambdaNest\equiv\left[
\begin{array}{cc}
n\bSigma E(\bXArv\bXTArv)+\sigma^2\bI_{\dA}&
\quad 
n\bSigma E(\bXArv\bXTBrv)\\[1ex]
\bO  & \sigma^2\bI_{\dB}
\end{array}
\right]
$$
which then leads to
$$\bbetahat\approx\bLambdaNest(\bX^T\bX)^{-1}\bLambdaNest^{-T}\bX^T\bY
\quad\mbox{and}\quad
\Cov(\bbetahat|\Xsc)\approx\bLambdaNest(\bX^T\bX)^{-1}.
$$
The bottom $\dB$ rows of $\bLambdaNest$ have the same
simple form as $\bLambdaCross$ and we obtain
$$\Cov(\bbetaBMLE|\Xsc)\approx
\left(\frac{\sigma^2}{\mbox{total sample size}}\right)
\Big[\mbox{lower right $\dB\times\dB$ block of}\ 
\left\{E\big(\bXrv\bXrv^T\big)\right\}^{-1}\Big]
$$
which matches (\ref{eq:raftTrip}) and, indeed, the asymptotic
covariance matrix form that arises in ordinary multiple regression.

\subsection{Closing Discussion on the Asymptotic Similarities}

In this section we have provided heuristic justifications for 
the fixed effects estimators and their asymptotic covariance
matrices having the approximate forms
\begin{equation}
\bbetahat\approx\bLambda(\bX^T\bX)^{-1}\bLambda^{-T}\bX^T\bY
\quad\mbox{and}\quad
\Cov(\bbetahat|\Xsc)\approx\bLambda(\bX^T\bX)^{-1}.
\label{eq:dumplingSandwich}
\end{equation}
for \emph{both} the crossed random effects model (\ref{eq:theModel})
and the nested model (\ref{eq:theNestModel}). The common approximate
forms in (\ref{eq:dumplingSandwich}) provide a reasonable explanation
for the asymptotic covariance matrices in \QuasiTheorem\ \ref{res:mainResult}
having forms similar to the nested case.

The approximate $\bbetahat$ expression in (\ref{eq:dumplingSandwich})
is intriguingly close to the well-known ordinary least 
expression. In the special case of $\bX$ being a column vector,
$\bLambda$ is scalar and cancels to give the ordinary
least squares form. Such reduction occurred in Section 
\ref{sec:CaliforniaFaded} for the $\bX=\bone$ case.
However there is no such cancellation in general.

The heuristics in the general $\bX$ cases involve approximations
having generic form
\begin{equation}
\bV\bX\approx \bX\bLambda.
\label{eq:SylvesterForm}
\end{equation}
In the special case where $\bX=\bx$ is a column vector and 
$\bLambda=\lambda$ is scalar then (\ref{eq:SylvesterForm}) 
becomes $\bV\bx\approx \bx\lambda$ which corresponds,
approximately, to $\lambda$ being an eigenvalue of $\bV$ 
with eigenvector $\bx$. For general $\bX$ and $\bLambda$
(\ref{eq:SylvesterForm}), with ``$=$'' instead of ``$\approx$'',
is an instance of \emph{Sylvester's equation}
(e.g. Stewart \myand Sun, 1990; Chapter V, Section 1.2).

\section{Statistical Utility}\label{sec:statsUtility}

\QuasiTheorem\ \ref{res:mainResult} provides a great deal of statistical utility 
concerning inference and design. Confidence intervals and Wald
hypothesis tests based on studentization are immediate consequences. 
Another is sample size calculations, for which we provide some details 
in this section.

For illustration of sample size calculations arising from 
\QuasiTheorem\ \ref{res:mainResult},
consider the following special case of (\ref{eq:theModel}):
\begin{equation}
\begin{array}{c}
Y_{ii'j}|\BinRViidj,X_{ii'j},U_i,U'_{i'}\simind 
N\Big(\beta^0_0+U_i+U'_{i'}+\beta^0_1\BinRViidj+\beta^0_2 X_{ii'j}+\beta^0_3 \BinRViidj
X_{ii'j},\sigma^2\Big),\\[2ex]
U_i\simind N(0,\Sigma^0),\ \ U'_{i'}\simind N\big(0,(\Sigma')^0\big),
\ \ 1\le i\le m,\ \ 1\le i'\le m',\ \ 1\le j\le n,
\end{array}
\label{eq:RoosterInAlbania}
\end{equation}
where the $\BinRViidj\simind\mbox{Bernoulli}(p)$ and the $X_{ii'j}$ being 
independently and identically distributed the same as a general random variable 
$\Xrv$ having finite second moment. Consider the one-sided hypotheses
\begin{equation}
H_0:\beta^0_3=0\quad\mbox{versus}\quad H_1:\beta^0_3>0
\label{eq:PickandsInPenn}
\end{equation}
corresponding to a possibly positive interaction effect between the two
predictors. Let $\Delta>0$ be a particular alternative value of $\beta_3^0$
and let $P$ be the corresponding power. Then \QuasiTheorem\ \ref{res:mainResult}
and standard arguments lead to the following sample size formula:
\begin{equation}
m={\Bigg\lceil}\frac{\{\Phi^{-1}(\alpha)+\Phi^{-1}(1-P)\}^2}{(\Delta/\sigma^0)^2
p(1-p)\Var(\Xrv)m'n}
{\Bigg\rceil}
\label{eq:WestOfVenus}
\end{equation}
where, for any $x\in\real$, $\lceil x\rceil$ denotes the smallest integer 
greater than or equal to $x$ and $\Phi^{-1}$ is the $N(0,1)$ quantile function.

Now consider a psychological study such that model (\ref{eq:RoosterInAlbania}) and
hypotheses (\ref{eq:PickandsInPenn}) apply with $m'=25$ items and $n=1$ observation
per subject-item combination. How many subjects should be recruited to potentially
detect a smallest meaningful interaction effect of $\Delta=0.25$ with power $0.9$
from a $0.05$ level of significance test? If it is further be assumed that $p=\smhalf$ 
and $\Var(X)=\frac{1}{12}$ then from (\ref{eq:WestOfVenus}) we should recruit
$$m=53\ \mbox{subjects if the error standard deviation is}\ \sigma^0=0.4.$$
Table \ref{tab:powSSexamp} below provides the required $m$ values for some
other values of $\sigma^0$.

In contemporary Gaussian response linear mixed model software, 
such as the function \texttt{lmer()} within the package \textsf{lme4} 
(Bates \textit{et al.}, 2015), standard errors are typically obtained  
using exact observed Fisher information rather than the approximation
to the (expected) Fisher information on which
(\ref{eq:WestOfVenus}) is based. This raises the question as to whether
the number of subjects chosen according to the \QuasiTheorem\ \ref{res:mainResult}
approximation to the standard error of $\betahat_3$ leads to the advertized
power for exact Fisher information-based hypothesis tests. We addressed
this question by running a simulation study that involved replication
of $1,000$ simulated data sets corresponding to (\ref{eq:RoosterInAlbania})
with various noise levels according to $\sigma^0\in\{0.2,0.4,0.8,1.6\}$.
The $\BinRViidj$ and $X_{iid'j}$ data were generated from $\mbox{Bernoulli}(\smhalf)$
and $\mbox{Uniform}(0,1)$ distributions, respectively. As above, we set
$(m',n,\Delta,\alpha,P)=(20,1,0.25,0.05,0.9)$ and determined $m$ using
(\ref{eq:WestOfVenus}). For each simulated data set we carried out a test
of (\ref{eq:PickandsInPenn}) using calls to \texttt{lmer()}, with 
rejection of $H_0$ if the t-statistic corresponding to $\beta_3^0$ exceeded
$\Phi^{-1}(1-\alpha)=\Phi^{-1}(0.95)$. Table \ref{tab:powSSexamp} shows the
empirical estimates of $P=0.9$ and corresponding 95\% confidence intervals.
For this example we see that the sample size formula (\ref{eq:WestOfVenus})
performs well with regards to the actual power delivered.

\begin{table}[!h]
\begin{center}
\begin{tabular}{lcccc}
\hline\\[-1.8ex]
Error standard deviation ($\sigma^0$): & 0.2   & 0.4 & 0.8 & 1.6          \\[0ex]
Minimum number of subjects ($m$):      & 14  & 53  & 211 & 842  \\[0ex]
Empirical estimate of power:           & 0.889 & 0.902 & 0.878 & 0.885  \\[0ex]
95\% confidence interval of power:     &$(0.870,0.908)$  & $(0.884,0.920)$ 
                                       & $(0.858,0.898)$ & $(0.865,0.905)$
\\[1ex]
\hline
\end{tabular}
\end{center}
\caption{\textit{The results from the illustrative sample size calculation and 
corresponding empirical power checks for the simulation study described in the text.
The number of subjects (}$m$\textit{) values correspond to an advertized power of $0.9$.}}
\label{tab:powSSexamp}
\end{table}

The example in this section demonstrates the statistical utility of 
\QuasiTheorem\ \ref{res:mainResult}. We are not aware of previous results in
the literature for linear mixed models with crossed random effects that 
readily provide the sample size formula (\ref{eq:WestOfVenus}).

\section{Concluding Remarks}\label{sec:concluding}

\QuasiTheorem\ \ref{res:mainResult} provides the precise leading term behaviours of the 
maximum likelihood estimators for a general class of linear mixed models containing
crossed random effects and enables statistical utilities such
as Wald tests for all model parameters and sample size calculations.
It complements the recent contributions of Lyu \textit{et al.} (2024) via
extensions to random slopes and unbalanced designs.
In comparison with the nested random effects situation, the establishment of leading
term results in the presence of crossed random effects is lengthy and arduous -- even
when the responses are Gaussian. The leading terms have similar or identical 
forms to those arising in nested models, 
and we have provided some heuristic arguments for this phenomenon.
We conjecture that the two-term asymptotic covariance matrices for $\bbetaAMLE$,
$\bSigmaMLE$ and $\bSigmadMLE$ in the Section \ref{sec:modelDescrip} set-up are similar 
or identical to those appearing in Section 3.3.1 of Maestrini \textit{et al.} (2024) 
for the nested case, but such an investigation would require a great deal of 
additional effort. Lastly, there are questions of what precise asymptotic results,
if any, could be obtained for non-Gaussian and sparse data versions of linear mixed 
models containing crossed random effects. 
The current article may pave the way for such future endeavours.

\section*{Supplementary Material}

The Supplementary Material contains the derivational details of Result 
\ref{res:mainResult}.

\ifthenelse{\boolean{UnBlinded}}{
\section*{Acknowledgements}

We are grateful to two reviewers for their
improvement suggestions. We also thank
John Duchi, Iain Johnstone, Ziyang Lyu, 
Luca Maestrini, Art Owen, Nicola Sartori and Alan Welsh 
for advice related to this research. This research was 
supported by the U.S. National Science Foundation grant 
DMS-2210569 and the Australian Research Council Discovery Project 
DP230101179. The second author is grateful for hospitality from 
the Department of Statistical Sciences, University of Padua, 
Italy, during part of this research.
}
{\null}

\section*{References}

\bib
Agresti, A. (2013). \textit{Categorical Data Analysis.} Hoboken, New Jersey:
John Wiley \& Sons.

\bib
Baayen, R.H., Davidson, D.J. \myand Bates, D.M. (2008).
Mixed-effects modeling with crossed random effects for 
subjects and items. {\it Journal of Memory and Language},
{\bf 59}, 390--412.

\bib
Bates, D., Maechler, M., Bolker, B. and Walker, S. (2015).
Fitting linear mixed-effects models using \textsf{lme4}. 
\textit{Journal of Statistical Software}, {\bf 67(1)}, 1--48.

\bib
Gao, K. \myand Owen, A.B. (2020).
Estimation and inference for very large linear mixed effects
models. \textit{Statistica Sinica}, \textbf{30}, 1741--1771.

\bib
Ghosh, S., Hastie, T. \myand Owen, A.B. (2022).
Backfitting for large scale crossed random effects regressions.
\textit{The Annals of Statistics}, \textbf{50}, 560--583.

\bib
Hartley, H.O. \myand Rao, J.N.K. (1967).
Maximum likelihood estimation for the mixed analysis of variance model.
\textit{Biometrika}, \textbf{54}, 93--108.

\bib
Jiang. J. (1996). REML estimation: asymptotic behavior 
and related topics.\ \textit{The Annals of Statistics}, \textbf{24}, 255--286.

\bib
Jiang, J. (2013). The subset argument and consistency of MLE in GLMM:
answer to an open problem and beyond.
\textit{The Annals of Statistics}, \textbf{41}, 177--195.

\bib
Jiang, J. \myand Nguyen, T. (2021). \textit{Linear and Generalized 
Linear Mixed Models and Their Applications, Second Edition.} New York: Springer.

\bib
Jiang, J., Wand, M.P. \myand Bhaskaran, A. (2022).
Usable and precise asymptotics for generalized linear
mixed model analysis and design.
\textit{Journal of the Royal Statistical Society, Series B},
\textbf{84}, 55--82.

\bib
Lyu, Z., Sisson, S.A. \myand Welsh, A.H. (2024).
Increasing dimension asymptotics for two-way crossed mixed effect models.
\textit{The Annals of Statistics}, \textbf{52}, 2956--2978.

\bib
Maestrini, L., Bhaskaran, A. \myand Wand, M.P. (2024).
Second term improvement to generalised linear mixed model asymptotics.
\textit{Biometrika}, \textbf{111}, 1077--1084.

\bib
Miller, J.J. (1973). Asymptotic properties and computation of maximum likelihood 
estimates in the mixed model of the analysis of variance. 
Technical Report No. 12, Department of Statistics, Stanford University,
Stanford, California, U.S.A.

\bib
Miller, J.J. (1977). Asymptotic properties of maximum likelihood estimates
in the mixed model of the analysis of variance. \textit{The Annals of 
Statistics}, \textbf{4}, 746--762.

\bib
Stewart, G.W. \myand Sun, J. (1990). \textit{Matrix Perturbation Theory}.
San Diego: Academic Press.

\vfill\eject

%
%
\renewcommand{\theequation}{S.\arabic{equation}}
\renewcommand{\thesection}{S.\arabic{section}}
\renewcommand{\thetable}{S.\arabic{table}}
\setcounter{equation}{0}
\setcounter{table}{0}
\setcounter{section}{0}
\setcounter{page}{1}
\setcounter{footnote}{0}

\centerline{\Large Supplement for:}
\vskip5mm
\centerline{\Large\bf Precise Asymptotics for Linear Mixed}
\vskip2mm
\centerline{\Large\bf Models with Crossed Random Effects}
\vskip5mm
%
\ifthenelse{\boolean{UnBlinded}}{
\centerline{\normalsize\sc Jiming Jiang$\null^1$, 
Matt P. Wand$\null^2$ and Swarnadip Ghosh$\null^3$}
\vskip5mm
\centerline{\textit{
$\null^1$University of California, Davis,
$\null^2$University of Technology Sydney
and $\null^3$Radix Trading}}
}
{\null}

\section{Derivation of \QuasiTheorem\ \ref{res:mainResult}}

In this section we provide a derivation of \QuasiTheorem\ \ref{res:mainResult}, 
starting with notation.

\subsection{Notation}

For any matrix $\bM$ let 
$$\bM^{\otimes2}\equiv\bM\bM^T\quad\mbox{and}\quad 
\Vert\bM\VertF\equiv\{\tr(\bM^T\bM)\}^{1/2}.$$ 
The latter definition is often called the \emph{Frobenius norm} of $\bM$.

The matrix $\bV\big(\bSigma,\bSigma',\sigsq\big)$ given by (\ref{eq:Vmatrix}) is central to the
derivations. Throughout this appendix, we omit the dependence on the covariance matrix parameters
by simply writing it as $\bV$. Define the following partitioning of the inverse of $\bV$:
$$\bV^{-1}
=\left[
\begin{array}{cccc}
\bV^{11} & \bV^{12} & \cdots & \bV^{1m} \\
\bV^{21} & \bV^{22} & \cdots & \bV^{2m} \\
\vdots   &  \vdots   & \ddots & \ddots   \\
\bV^{m1} & \bV^{m2}  & \cdots & \bV^{mm} 
\end{array}
\right]
\ \ \mbox{where}\ \ \bV^{i\iunder}\ \  \mbox{is}
\ \ \left(\sum_{i'=1}^{m'} n_{ii'}\right)\times
\left(\sum_{i'=1}^{m'} n_{\iunder i'}\right).
$$

If $\Psc$ is a logical proposition then $I(\Psc)=1$ if $\Psc$ is true.
Otherwise, $I(\Psc)=0$.

\subsection{Lemmas}

The upcoming Fisher information approximations rely on four lemmas, 
which we present here.

\subsubsection{A Lemma that Provides a Simple Kronecker Product Form}

\begin{lemma}
Let $\bA_d$ be a symmetric $d\times d$ matrix with $(r,s)$th entry 
denoted by $A_{rs}$. Also, let $\bB_d$ be the $\smhalf d(d+1)\times\smhalf d(d+1)$
matrix with entries determined according to the following table:
%
\begin{table}[!h]
\begin{center}
\begin{tabular}{cc}
\hline\\[-1.8ex]
\textit{entry of $\vech(\bA_d)\vech(\bA_d)^T$} &\textit{entry of $\bB_d$ in the same position} \\[0ex]
\hline\\[-1.8ex]
$A_{rr}A_{tt}$ & $A_{rt}^2$\\[0ex]
$A_{rr}A_{tu}, t\ne u$  & $2A_{rt}A_{ru}$\\[0ex]
$A_{rs}A_{tu}, r\ne s,t\ne u$  & $2(A_{rt}A_{su}+A_{ru}A_{st})$\\[1ex]
\hline
\end{tabular}
\end{center}
\caption{\textit{Definition of the matrix $\bB_d$, a function
of a $d\times d$ symmetric matrix $\bA_d$.}}
\label{tab:cracker}
\end{table}
%

\noindent
Then 
$$\bB_d=\bD_d^T(\bA_d\otimes\bA_d)\bD_d.$$
\label{lem:Kron}
\end{lemma}

\subsubsection{Three Lemmas Stating Key Matrix Identities}

The following three lemmas state some matrix identities which play key 
roles in the derivation of \QuasiTheorem\ \ref{res:mainResult}.

\begin{lemma}
Let $\lambda>0$, $\bA$ be a invertible $d\times d$ matrix and $\bX$, $\bXdot$ and $\bXdotdot$
each be $n\times d$ matrices, where $n,d\in\naturalNumbers$.
Then, assuming that all required matrix inverses exist,
{\setlength\arraycolsep{2pt}
\begin{eqnarray*}
\bXTdot(\bX\bA\bX^T+\lambda\bI)^{-1}\bXdotdot&=&(1/\lambda)\bXTdot
\{\bI-\bX(\bX^T\bX)^{-1}\bX^T\}\bXdotdot\\[1ex]
&&\qquad+\bXTdot\bX(\bX^T\bX)^{-1}\{\bA+\lambda(\bX^T\bX)^{-1}\}^{-1}
(\bX^T\bX)^{-1}\bX^T\bXdotdot.
\end{eqnarray*}
}
\label{lem:keyIdents}
\end{lemma}

\noindent
Lemma \ref{lem:keyIdents} has the following immediate corollary:

\noindent
\textbf{Corollary 2.1.}
\textit{If $\lambda$, $\bA$, $\bX$, $\bXdot$ and $\bXdotdot$ are as 
defined in Lemma \ref{lem:keyIdents} then, under the  Lemma \ref{lem:keyIdents}
assumptions:}
\begin{itemize}
\item[\mbox{\rm{(a)}}] 
$\bX^T(\bX\bA\bX^T+\lambda\bI)^{-1}\bXdotdot=\{\bA+\lambda(\bX^T\bX)^{-1}\}^{-1}(\bX^T\bX)^{-1}\bX^T\bXdotdot.$
\item[\mbox{\rm{(b)}}] 
$\bXTdot(\bX\bA\bX^T+\lambda\bI)^{-1}\bX=\bXTdot\bX(\bX^T\bX)^{-1}\{\bA+\lambda(\bX^T\bX)^{-1}\}^{-1}.$
\item[\mbox{\rm{(c)}}] 
$\bX^T(\bX\bA\bX^T+\lambda\bI)^{-1}\bX=\{\bA+\lambda(\bX^T\bX)^{-1}\}^{-1}.$
\end{itemize}

The following related matrix identity is also important:
\begin{lemma}
\textit{If $\lambda$, $\bA$ and $\bX$ are as defined in Lemma \ref{lem:keyIdents} then,
assuming all required matrix inverses exist,}
$$\bX^T(\bX\bA\bX^T+\lambda\bI)^{-2}\bX=\{\bA+\lambda(\bX^T\bX)^{-1}\}^{-1}
(\bX^T\bX)^{-1}\{\bA+\lambda(\bX^T\bX)^{-1}\}^{-1}.
$$
\label{lem:keySqIdent}
\end{lemma}

In addition, the derivation of \QuasiTheorem\ \ref{res:mainResult} makes use of: 

\begin{lemma}
\textit{Let $\bA$ and $\bB$ be $d\times d$ matrices such that 
each of}
$$
\left[
\begin{array}{cc}
\bA & \bB \\[2ex]
\bB & \bA
\end{array}
\right]\ \mbox{and}\ \bA+\bB\ \mbox{are invertible.}
$$
\textit{Then}
$$
\left[
\begin{array}{c}
\bI_d\\[1ex]
\bI_d
\end{array}
\right]^T
\left[
\begin{array}{cc}
\bA & \bB \\[2ex]
\bB & \bA
\end{array}
\right]^{-1}
\left[
\begin{array}{c}
\bI_d\\[1ex]
\bI_d
\end{array}
\right]=2(\bA+\bB)^{-1}.
$$
\label{lem:ApsBinvResIden}
\end{lemma}

\subsubsection{Lemmas for Limits of Forms Arising in the Fisher Information Matrix}

Here we provide three convergence in probability lemmas
that are key to dealing with particular forms that arise in the
Fisher information matrix.

First we present Lemma \ref{lem:XApart}  which identifies some key convergence in probability
limits related to predictor summation quantities about the $\bV^{-1}$ matrix.
Let $\bXrv$ be a $d\times1$ random vector and let 
\begin{equation}
\bX_{ii'j},\quad 1\le i\le m,\quad 1\le i'\le m',\quad 1\le j\le n_{ii'},
\label{eq:lemCA}
\end{equation}
be independent and identically distributed random vectors 
having the same distribution as $\bXrv$. Then define
for $1\le i\le m$ and $1\le i'\le m'$:
\begin{equation}
\bX_{ii'}\equiv\left[
\begin{array}{c}
\bX_{ii'1}^T\\[1ex]
\vdots\\
\bX_{ii'n_{ii'}}^T
\end{array}
\right],\ \ 
\bX\equiv\stack{1\le i \le m}\big(\bXuptrii\big)\ \ \mbox{where}\ \ 
\bXuptrii\equiv\stack{1\le i'\le m'}(\bX_{ii'}).
\label{eq:lemCB}
\end{equation}
Next, let 
\begin{equation}
\begin{array}{l}
\bQ_{mm'}\equiv
{\displaystyle
\blockdiag{1\le i\le m}
\left\{\blockmatrix{1\le i',\idunder\le m'}(\bXiid\bM\bXiidu^T)\right\}}\\[1ex]
\qquad\qquad\qquad
+{\displaystyle\blockmatrix{1\le i,\iunder\,\le m}
\left\{\blockdiag{1\le i'\le m'}(\bXiid\bM'\bXiuid^T)\right\}}+\lambda\bI
\end{array}
\label{eq:lemCC}
\end{equation}
where
\begin{equation}
\bM\ \mbox{and}\ \bM'\ \mbox{are}\ d\times d\ 
\mbox{symmetric positive definite matrices and}\ \lambda>0.
\label{eq:lemCD}
\end{equation}
Partition $\bQ_{mm'}^{-1}$ as follows
\begin{equation}
\bQ_{mm'}^{-1}
=\left[
\begin{array}{cccc}
\bQmmd^{11} & \bQmmd^{12} & \cdots & \bQmmd^{1m} \\[2ex]
\bQmmd^{21} & \bQmmd^{22} & \cdots & \bQmmd^{2m} \\[0ex]
\vdots   &  \vdots   & \ddots & \vdots   \\
\bQmmd^{m1} & \bQmmd^{m2}  & \cdots & \bQmmd^{mm} 
\end{array}
\right]
\ \ \mbox{where}\ \ \bQmmd^{i\iunder}\ \  \mbox{is}
\ \ \left(\sum_{i'=1}^{m'} n_{ii'}\right)\times
\left(\sum_{i'=1}^{m'} n_{\iunder i'}\right).
\label{eq:lemCE}
\end{equation}

Introduce the following assumptions:
\begin{itemize}
\item[(A4)] All entries of $\bXrv$ are not degenerate at zero and
have finite second moment.
\item[(A5)] Each of the $n_{ii'}$, $1\le i\le m,\ 1\le i'\le m'$, diverge to $\infty$.
\end{itemize}

\begin{lemma}
Let $\bXrv$ be a $d\times1$ random vector for which (A4) holds.
For $m,m'\in\naturalNumbers$ define $\bX$, $\bXuptrii$, $\bQmmd$ and 
$\bQmmd^{i\iunder}$, $1\le i\le m$, $1\le\iunder\le m'$,  according to 
(\ref{eq:lemCA})--(\ref{eq:lemCE}). Under (A5) we have for
fixed $m,m'\in\naturalNumbers$:
\begin{itemize}
\item[(a)] $\bX^T\bQmmd^{-1}\bX\convprob\left(\oom\bM + \oomd\bM'\right)^{-1}$.
\item[(b)] For all $1\le i\le m$,\ \  $\bXuptrii^T\bQmmd^{ii}\bXuptrii\convprob
\bM^{-1}-\frac{1}{mm'}\bM^{-1}\bM'\left(\oom\bM + \oomd\bM'\right)^{-1}$.
\item[(c)] If $m\ge2$ then for all $1\le i,\iunder\le m$ such that $i\ne\iunder$,
{\setlength\arraycolsep{1pt}
\begin{eqnarray*}
&&\bXuptrii^T\bQmmd^{i\iunder}\bXuptriiu\convprob\,
-{\textstyle\frac{1}{mm'}}\bM^{-1}\bM'\left(\oom\bM + \oomd\bM'\right)^{-1}.
\end{eqnarray*}
}
\item[(d)] $\left(\displaystyle{\sumim\sum_{i'=1}^{m'}}n_{ii'}\right)^{-1}\mbox{\rm tr}
(\bQmmd^{-2})\convprob 1/\lambda^2$.
\end{itemize}
\label{lem:XApart}
\end{lemma}

Let $\bXdiamdrv$ be a $\ddiamd\times1$ random vector and let 
\begin{equation}
\bXdiamd_{ii'j},\quad 1\le i\le m,\quad 1\le i'\le m',\quad 1\le j\le n_{ii'},
\label{eq:DiamondDog}
\end{equation}
be independent and identically distributed random vectors 
having the same distribution as $\bXdiamdrv$. Then define
for $1\le i\le m$ and $1\le i'\le m'$:
\begin{equation}
\bXdiamd_{ii'}\equiv\left[
\begin{array}{c}
\bXdiamd_{ii'1}^T\\[1ex]
\vdots\\
\bXdiamd_{ii'n_i}^T
\end{array}
\right]\quad\mbox{and}\quad
\bXdiamd\equiv\stack{1\le i \le m}\left\{\stack{1\le i'\le m'}(\bXdiamd_{ii'})\right\}.
\label{eq:PhillipOfSpain}
\end{equation}

\begin{lemma}
Let $\bXrv$ be a $d\times1$ random vector and $\bXdiamdrv$ be a $\ddiamd\times1$ 
random vector such that for (A4) holds for both $\bXrv$ and $\bXdiamdrv$.
Define $\bX$ according to (\ref{eq:lemCA})--(\ref{eq:lemCB}),
$\bQmmd$ according to (\ref{eq:lemCB})--(\ref{eq:lemCD})
and
$\bXdiamd$ according to (\ref{eq:DiamondDog})--(\ref{eq:PhillipOfSpain}).
Under (A5) we have for all fixed $m,m'\in\naturalNumbers$:
\begin{itemize}
\item[(a)]
$\left({\displaystyle\sum_{i=1}^m\sum_{i'=1}^{m'}n_{ii'}}\right)^{-1}
\bXdiamd^T\bQmmd^{-1}\bXdiamd\convprob(1/\lambda)
\Big[\mbox{lower right $\ddiamd\times\ddiamd$ block of}\
\big\{E\big([\bXrv\  \bXTdiamdrv]^{\otimes2}\big)\big\}^{-1}\Big]^{-1}
$.
\item[(b)]
$
\bX^T\bQmmd^{-1}\bXdiamd\convprob \left(\oom\bM + \oomd\bM'\right)^{-1}\{E(\bXrv^{\otimes2})\}^{-1}
E(\bXrv\bXdiamdrv^T)
$.
\end{itemize}
\label{lem:XBpart}
\end{lemma}

\subsection{Fisher Information Matrix Approximation}\label{sec:FishInfMatAppr}

The Fisher information matrix of the full vector of unique parameters,
corresponding to the conditional log-likelihood (\ref{eq:FLbillboards}),
is denoted by
\begin{equation}
I\big(\bbetaA,\bbetaB,\vech(\bSigma),\vech(\bSigma'),\sigsq\big).
\label{eq:FishInfoMat}
\end{equation}
We now obtain approximations to each of the sub-blocks of (\ref{eq:FishInfoMat}).

From (A1), $m'$ has the same order of magnitude as $m$. Therefore,
remainder terms such as $o_P(mm'n)$ can be also written as $o_P(m^2n)$.
Throughout this derivation we follow the convention of expressing all remainder terms
that involve $m$ and $m'$ in terms of $m$ only.

\subsubsection{The $(\bbetaA,\bbetaA)$ Diagonal Block}

The $(\bbetaA,\bbetaA)$ diagonal block is $\bXA^T\bV^{-1}\bXA$. From (A3)
and Lemma \ref{lem:XApart}(a), we have for all fixed $m,m'\in\naturalNumbers$ and as $n\to\infty$ 
$$\bXA^T\bV^{-1}\bXA\convprob\left(\frac{\bSigma}{m}+\frac{\bSigma'}{m'}\right)^{-1}.$$
Therefore, under (A1) and (A3), the $(\bbetaA,\bbetaA)$ diagonal block of the
Fisher information matrix is 
$$\left(\frac{\bSigma}{m}+\frac{\bSigma'}{m'}\right)^{-1}+o_P(m)\bone_{\dA}^{\otimes2}.$$

\subsubsection{The $(\bbetaB,\bbetaB)$ Diagonal Block}

The $(\bbetaB,\bbetaB)$ diagonal block is $\bXB^T\bV^{-1}\bXB$. Under (A2)--(A3),
and applying Lemma \ref{lem:XBpart}(a) with $\bX=\bXA$ and $\bXdiamd=\bXB$ we have 
$$\bXB^T\bV^{-1}\bXB=\frac{mm'n\COMPQUANT_{\bbetaB}^{-1}}{\sigma^2}+o_P(mm'n)\bone_{\dB}^{\otimes2}.$$

\subsubsection{The $\big(\vech(\bSigma),\vech(\bSigma)\big)$ Diagonal Block}

From results given in e.g. Section 4.3 of Wand (2002), the $(\bSigma_{rs},\bSigma_{tu})$
entry of the $\big(\vech(\bSigma),\vech(\bSigma)\big)$ diagonal block of the
Fisher information matrix is 
$$\smhalf\tr\left(\bV^{-1}\frac{\partial\bV}{\partial(\bSigma)_{rs}}
\bV^{-1}\frac{\partial\bV}{\partial(\bSigma)_{tu}}\right).
$$
Then note that 
\begin{equation}
\frac{\partial\bV}{\partial(\bSigma)_{rs}}
=\bL_r\bL_s^T+I(r\ne s)\bL_s\bL_r^T\ \ \mbox{where}\ \ 
\bL_r\equiv\blockdiag{1\le i\le m}\big(\bXuptriAi\be_r\big),
\quad \bXuptriAi\equiv\stack{1\le i'\le m'}(\bXAiid)
\label{eq:DeeDeeMyers}
\end{equation}
and $\be_r$ denotes the $\dA\times 1$ matrix with $r$th entry $1$ and all other
entries $0$. Noting the $\tr(\bA\bB)=\tr(\bB\bA)$ identity for all compatible
matrices $\bA$ and $\bB$ and introducing the notation
$$\trQuant_{rstu}\equiv \tr\big\{(\bL_r^T\bV^{-1}\bL_s)^T(\bL_t^T\bV^{-1}\bL_u)\big\}.$$
we then have the following simplifications of the various sub-types of the $(\bSigma_{rs},\bSigma_{tu})$ 
Fisher information blocks:
\begin{equation}
{\setlength\arraycolsep{3pt}
\begin{array}{lcl}
(\bSigma_{rr},\bSigma_{tt}):&$\quad$&\smhalf\trQuant_{rtrt}\\[1ex]
%
%
(\bSigma_{rr},\bSigma_{tu}),\ t\ne u:&$\quad$&\smhalf\Big(\trQuant_{rurt}+\trQuant_{rtru}\Big)\\[1ex]
(\bSigma_{rs},\bSigma_{tt}),\ r\ne s:&$\quad$&\smhalf\Big(\trQuant_{rtst}+\trQuant_{strt}\Big)\\[1ex]
(\bSigma_{rs},\bSigma_{tu}),\ r\ne s,t\ne u:&$\quad$&\smhalf\Big(\trQuant_{rust}+\trQuant_{rtsu}
+\trQuant_{surt}+\trQuant_{stru}\Big).
\end{array}
}
\label{eq:MustardSeed}
\end{equation}
Since
$$\bL_r^T\bV^{-1}\bL_s=\Big[\be_r^T\bXuptriAi^T\bV^{i\iunder}
\bXuptriAiu\be_s\Big]_{1\le i\le m,1\le\,\iunder\,\le m},
$$
we then have
{\setlength\arraycolsep{3pt}
\begin{eqnarray*}
\trQuant_{rstu}&=&\sumim\sum_{\iunder=1}^m\Big(\be_r^T\bXuptriAi^T\bV^{i\iunder}
\bXuptriAiu\be_s\Big)\Big(\be_t^T\bXuptriAi^T\bV^{i\iunder}
\bXuptriAiu\be_u\Big)\\[1ex]
&=&\sumim\Big(\be_r^T\bXuptriAi^T\bV^{ii}
\bXuptriAi\be_s\Big)\Big(\be_t^T\bXuptriAi^T\bV^{ii}\bXuptriAi\be_u\Big)\\[1ex]
&&\quad+\mathop{\sum\sum}_{i\ne\,\iunder}
\Big(\be_r^T\bXuptriAi^T\bV^{i\iunder}
\bXuptriAiu\be_s\Big)\Big(\be_t^T\bXuptriAi^T\bV^{i\iunder}
\bXuptriAiu\be_u\Big).
\end{eqnarray*}
}
Lemma \ref{lem:XApart} (b)--(c) implies that for any fixed $m\in\{2,3,\ldots\}$ 
and $m'\in\naturalNumbers$ we have, as $n\to\infty$,
{\setlength\arraycolsep{3pt}
\begin{eqnarray*}
\trQuant_{rstu}&\convprob&
m\Big(\bSigma^{-1}-{\textstyle\frac{1}{mm'}}\bSigma^{-1}\bSigma'
\left(\oom\bSigma + \oomd\bSigma'\right)^{-1}\Big)_{rs}
\Big(\bSigma^{-1}-{\textstyle\frac{1}{mm'}}\bSigma^{-1}\bSigma'
\left(\oom\bSigma + \oomd\bSigma'\right)^{-1}\Big)_{tu}\\[1ex]
&&\ +{\displaystyle\frac{m(m-1)}{(mm')^2}}
\Big(\bSigma^{-1}\bSigma'\left(\oom\bSigma + \oomd\bSigma'\right)^{-1}\Big)_{rs}
\Big(\bSigma^{-1}\bSigma'\left(\oom\bSigma + \oomd\bSigma'\right)^{-1}\Big)_{tu}.
\end{eqnarray*}
}
Now suppose that $m$ and $m'$ diverge according to (A1). Then straightforward steps show 
that
\begin{equation}
\trQuant_{rstu}=m\big(\bSigma^{-1}\big)_{rs}\big(\bSigma^{-1}\big)_{tu}+O_P(1).
\label{eq:AnnesBeer}
\end{equation}
In view of (\ref{eq:MustardSeed}) and (\ref{eq:AnnesBeer}), under 
(A1) and (A2), the entries of the $\big(\vech(\bSigma),\vech(\bSigma)\big)$ 
diagonal block have the following leading term behavior: 
$$
{\setlength\arraycolsep{3pt}
\begin{array}{lcl}
(\bSigma_{rr},\bSigma_{tt}):&$\quad$&\smhalf m(\bSigma^{-1})^2_{rt}+O_P(1)\\[1ex]
%
%
(\bSigma_{rr},\bSigma_{tu}),\ t\ne u:&$\quad$& m(\bSigma^{-1})_{rt}(\bSigma^{-1})_{ru}+O_P(1)\\[1ex]
(\bSigma_{rs},\bSigma_{tt}),\ r\ne s:&$\quad$&  m(\bSigma^{-1})_{rt}(\bSigma^{-1})_{st}+O_P(1) \\[1ex]
(\bSigma_{rs},\bSigma_{tu}),\ r\ne s,t\ne u:&$\quad$&m\big\{(\bSigma^{-1})_{rt}(\bSigma^{-1})_{su}
+(\bSigma^{-1})_{ru}(\bSigma^{-1})_{st}
\big\}+O_P(1).
\end{array}
}
$$
Application of Lemma \ref{lem:Kron} then leads to the following succinct expression for the
$\big(\vech(\bSigma),\vech(\bSigma)\big)$ Fisher information block:
$$\smhalf m\bD_{\dA}^T(\bSigma^{-1}\otimes\bSigma^{-1})\bD_{\dA}+O_P(1)\bone_{\dA(\dA+1)/2}^{\otimes2}.$$

\subsubsection{The $\big(\vech(\bSigma'),\vech(\bSigma')\big)$ Diagonal Block}

The conditional log-likelihood is unaffected by the interchanging of 
$\bSigma$ and $\bSigma'$. Hence, noting the conclusion of the previous subsection, 
the $\big(\vech(\bSigma'),\vech(\bSigma')\big)$ diagonal block of the 
Fisher information is
$$\smhalf m'\bD_{\dA}^T\big((\bSigma')^{-1}\otimes(\bSigma')^{-1}\big)\bD_{\dA}
+O_P(1)\bone_{\dA(\dA+1)/2}^{\otimes2}.$$ 

\subsubsection{The $\big(\sigsq,\sigsq\big)$ Diagonal Block}

Appealing again to Section 4.3 of Wand (2002), the $(\sigsq,\sigsq)$ diagonal block 
of the Fisher information matrix is 
$$\smhalf\tr\left(\bV^{-1}\frac{\partial\bV}{\partial\sigsq}
\bV^{-1}\frac{\partial\bV}{\partial\sigsq}\right)=\smhalf\tr(\bV^{-2})
=\frac{mm'n}{2\sigma^4}+o_P(m^{-2}n^{-1}),
$$
with the last equality following from Lemma \ref{lem:XApart}(d).

\subsubsection{The $(\bbetaA,\bbetaB)$ Off-Diagonal Block}

The $(\bbetaA,\bbetaB)$ diagonal block is $\bXA^T\bV^{-1}\bXB$. 
From (A3) and Lemma \ref{lem:XBpart}(b), we have for all fixed $m,m'\in\naturalNumbers$ and as $n\to\infty$ 
$$\bXA^T\bV^{-1}\bXB\convprob\left(\frac{\bSigma}{m}+\frac{\bSigma'}{m'}\right)^{-1}
\{E(\bXArv^{\otimes2})\}^{-1}E(\bXArv^T\bXBrv).$$
Therefore, under (A1) and (A3), the $(\bbetaA,\bbetaB)$ diagonal block of the
Fisher information matrix is 
$$\left(\frac{\bSigma}{m}+\frac{\bSigma'}{m'}\right)^{-1}\{E(\bXArv^{\otimes2})\}^{-1}E(\bXArv^T\bXBrv)
+o_P(m).$$

\subsubsection{The $\Big(\big(\bbetaA,\bbetaB\big),
\big(\vech(\bSigma),\vech(\bSigma'),\sigsq\big)\Big)$ 
Off-Diagonal Block}\label{sec:exactOrthog}

From e.g. Section 4.3 of Wand (2002), the 
$$\Big(\big(\bbetaA,\bbetaB\big),\big(\vech(\bSigma),\vech(\bSigma'),\sigsq\big)\Big)$$
off-diagonal block is a matrix having all entries equal to zero. 
In other words, the fixed effects parameters and the covariance
matrix parameters are exactly orthogonal in Gaussian response linear mixed models.

\subsubsection{The $\big(\vech(\bSigma),\vech(\bSigma')\big)$ Off-Diagonal Block}
\label{sec:firstEig}

We commence with the special case of $\dA=1$, $n_{ii'}=n$ and 
$\bXAiid=\bone_n$ for all $1\le i\le m$, $1\le i'\le m'$.
In this case $\bone_{mm'n}$ is an eigenvector of 
$\bV$ with corresponding eigenvalue $m'n\Sigma+mn\Sigma'+\sigma^2$.
This implies that $\bone_{mm'n}$ is also an eigenvector of $\bV^{-1}$
with the just-mentioned eigenvalue reciprocated.  
Relatively straightforward 
manipulations then lead to the following expression
for the $(\Sigma,\Sigma')$ entry of the Fisher information matrix:
\begin{equation}
\smhalf\big[\{\Sigma(m'/m)+\Sigma'+\sigma^2/(mn)\}
\{\Sigma+\Sigma'(m/m')+\sigma^2/(m'n)\}\big]^{-1}
\label{eq:KemblaBowlo}
\end{equation}
which is $O(1)$ under assumption (A1).

Next we treat the general $\dA$, $n_{ii'}$ and $\bXAiid$ situation
with $m\in\naturalNumbers$ and $m'=1$. From e.g. Section 4.3 of Wand (2002), 
the $(\bSigma_{rr},\bSigma'_{tt})$
entry of the $\big(\vech(\bSigma),\vech(\bSigma')\big)$ off-diagonal block of the
Fisher information matrix is 
\begin{equation}
\smhalf\tr\left(\bV^{-1}\frac{\partial\bV}{\partial(\bSigma)_{rr}}
\bV^{-1}\frac{\partial\bV}{\partial(\bSigma')_{tt}}\right)
\label{eq:dotItalo}
\end{equation}
where, noting the current $m'=1$ special case,
\begin{equation}
\frac{\partial\bV}{\partial(\bSigma)_{rr}}=
\blockdiag{1\le i\le m}\big(\bXAione\be_r\be_r^T\bXAione^T\big)
\ \mbox{and}\ 
\frac{\partial\bV}{\partial(\bSigma')_{tt}}
=\blockmatrix{1\le i,\iunder\le m}\big(\bXAione\be_t\be_t^T\bXAiuone^T\big).
\label{eq:thePlug}
\end{equation}
Substitution of (\ref{eq:thePlug}) into (\ref{eq:dotItalo}) and algebraic
manipulations such as those involving the $\tr(\bA\bB)=\tr(\bB\bA)$ identity
lead to 
{\setlength\arraycolsep{1pt}
\begin{eqnarray*}
\tr\left(\bV^{-1}\frac{\partial\bV}{\partial(\bSigma)_{rr}}
\bV^{-1}\frac{\partial\bV}{\partial(\bSigma')_{tt}}
\right)
&=&\sum_{i=1}^m\sum_{\iunder=1}^m\sum_{\iunderTwo=1}^m
\big(\be_t^T\bXAione^T\bV^{i\iunder}\bXAiuone\be_r\big)
\big(\be_r^T\bXAiuone^T\bV^{\iunder\iunder^*}\bXAiuTwoOne\be_t\big)\\[1ex]
&=&\sum_{i=1}^m \big(\be_r^T\bXAione^T\bV^{ii}\bXAione\be_t\big)^2\\[1ex]
&&\ +\mathop{\sum\sum}_{i\ne\iunderTwo}
\big(\be_t^T\bXAione^T\bV^{ii}\bXAione\be_r\big)
\big(\be_r^T\bXAione^T\bV^{i\iunder^*}\bXAiuTwoOne\be_t\big)\\[1ex]
&&\ +\mathop{\sum\sum}_{i\ne\iunder}
\big(\be_t^T\bXAione^T\bV^{i\iunder}\bXAiuone\be_r\big)
\big(\be_r^T\bXAiuone^T\bV^{\iunder\iunder}\bXAiuone\be_t\big)\\[1ex]
&&\ +\mathop{\sum\sum\sum}_{i\ne\iunder\ne\iunderTwo}
\big(\be_t^T\bXAione^T\bV^{i\iunder}\bXAiuone\be_r\big)
\big(\be_r^T\bXAiuone^T\bV^{\iunder\iunder^*}\bXAiuTwoOne\be_t\big).
\end{eqnarray*}
}
Lemma \ref{lem:XApart}(b) and \ref{lem:XApart}(c) then imply that 
{\setlength\arraycolsep{1pt}
\begin{eqnarray*}
\smhalf\tr\left(\bV^{-1}\frac{\partial\bV}{\partial(\bSigma)_{rr}}
\bV^{-1}\frac{\partial\bV}{\partial(\bSigma')_{tt}}
\right)
&\convprob&\smhalf m\Big(\bM^{-1}-\textstyle{\frac{1}{m}}\bM^{-1}\bM'
\left(\textstyle{\frac{1}{m}}\bM+\bM'\right)\Big)_{rt}^2\\[1ex]
&&\quad+\smhalf m(m-1)\Big(\bM^{-1}-\textstyle{\frac{1}{m}}\bM^{-1}\bM'
\left(\textstyle{\frac{1}{m}}\bM+\bM'\right)\Big)_{tr}\\
&&\qquad\quad\times \Big(-\textstyle{\frac{1}{m}}\bM^{-1}\bM'
\left(\frac{1}{m}\bM+\bM'\right)\Big)_{rt}\\[1ex]
&&\quad+\smhalf m(m-1)\Big(-\textstyle{\frac{1}{m}}\bM^{-1}\bM'
\left(\textstyle{\frac{1}{m}}\bM+\bM'\right)\Big)_{tr}\\
&&\qquad\quad\times \Big(\bM^{-1}-\textstyle{\frac{1}{m}}\bM^{-1}\bM'
\left(\frac{1}{m}\bM+\bM'\right)\Big)_{rt}\\[1ex]
&&\quad+\smhalf m(m-1)^2 \Big(-\textstyle{\frac{1}{m}}\bM^{-1}\bM'
\left(\textstyle{\frac{1}{m}}\bM+\bM'\right)\Big)_{tr}\\
&&\qquad\quad\times\Big(-\textstyle{\frac{1}{m}}\bM^{-1}\bM'
\left(\textstyle{\frac{1}{m}}\bM+\bM'\right)\Big)_{rt}\\[1ex]
&=&\textstyle{\frac{1}{2m}}\Big(\left(\frac{1}{m}\bM+\bM'\right)^{-1}\Big)_{rt}^2
\end{eqnarray*}
}
\null\!\!\!
after several algebraic steps and cancellations.
The $r\ne s$ and $t\ne u$ cases are similar.
This confirms that (\ref{eq:KemblaBowlo}) also holds in general, with
the exception of $m'$ being set to $1$. For $m'\ge2$ similar arguments
can be used to show that the summations in (\ref{eq:dotItalo}) 
lead to convergents analogous to those in the $\dA=1$, $n_{ii'}=n$ and 
$\bXAiid=\bone_n$ case and a matrix with order
$O(1)\bone_{\dA(\dA+1)/2}^{\otimes 2}$ under (A1) eventuates.
 
\subsubsection{The $\big(\vech(\bSigma),\sigsq\big)$ Off-Diagonal Block}

We commence with the special case of $\dA=1$, $n_{ii'}=n$ and 
$\bXAiid=\bone_n$ for all $1\le i\le m$, $1\le i'\le m'$.
Using the eigenvalue and eigenvector properties described
near the beginning of Section \ref{sec:firstEig},
relatively straightforward 
manipulations then lead to the following expression
for the $(\Sigma,\sigma^2)$ entry of the Fisher information matrix:
\begin{equation}
{\setlength\arraycolsep{1pt}
\begin{array}{l}
\displaystyle{
\frac{m'(1-1/m)\{\Sigma+\Sigma'(m/m')+\sigma^2/(m'n)\}^2}{2mn
\{\Sigma(m'/m)+\Sigma'+\sigma^2/(mn)\}^2\{\Sigma+\sigma^2/(m'n)\}^2}}\\[2ex]
\qquad\qquad\qquad\qquad\qquad
+\displaystyle{\frac{m'}{2m^2n\{\Sigma(m'/m)+\Sigma'+\sigma^2/(mn)\}^2}}
\end{array}
}
\label{eq:AroundNearPole}
\end{equation}
which is $O(n^{-1})$ under (A1).

Now consider the general $\dA$, $n_{ii'}$ and $\bXAiid$ situation
with $m\in\naturalNumbers$ and $m'=1$. Results in e.g. Section 4.3 of Wand (2002)
imply that the $(\bSigma_{rr},\sigma^2)$ entry of the $\big(\vech(\bSigma),\sigma^2\big)$ 
off-diagonal block of the Fisher information matrix is 
\begin{equation}
\smhalf\tr\left(\bV^{-2}
\blockdiag{1\le i\le m}\Big(\bXAione\be_r\be_r^T\bXAione^T\Big)\right)
=\smhalf
\sum_{i=1}^m\sum_{\iunder=1}^m
\be_r^T\big(\bV^{\iunder\,i}\bXAione\big)^T\big(\bV^{\iunder\,i}\bXAione\big)\be_r.
\label{eq:coffeeWithCannoli}
\end{equation}

For the $m=m'=1$ case (\ref{eq:coffeeWithCannoli}) use of Lemma \ref{lem:keySqIdent}
leads to
{\setlength\arraycolsep{1pt}
\begin{eqnarray*}
n_{11}\bXAoneone^T\bV_{11}^{-2}\bXAoneone
&=&
n_{11}\be_r^T\bXAoneone^T\{\bXAoneone(\bSigma+\bSigma')\bXAoneone^T+\sigma^2\bI\}^{-2}\bXAoneone\be_r\\[1ex]
&=&\be_r^T\big\{\bSigma+\bSigma'+\sigma^2(\bXAoneone^T\bXAoneone)^{-1}\big\}^{-1}
\left(\textstyle{\frac{1}{n_{11}}}\bXAoneone^T\bXAoneone\right)^{-1}\\
&&\qquad\times\big\{\bSigma+\bSigma'+\sigma^2(\bXAoneone^T\bXAoneone)^{-1}\big\}^{-1}\be_r\\[1ex]
&\convprob&\be_r^T(\bSigma+\bSigma')^{-1}E(\bXTrv\bXrv)(\bSigma+\bSigma')^{-1}\be_r^T.
\end{eqnarray*}
}
Hence, the $(\bSigma_{rr},\sigma^2)$ entry of the Fisher information is 
$$\frac{\big((\bSigma+\bSigma')^{-1}E(\bXTrv\bXrv)(\bSigma+\bSigma')^{-1}\big)_{rr}
\{1+o_P(1)\}}{2n_{11}}$$
which extends (\ref{eq:AroundNearPole}) for $\dA\in\naturalNumbers$ and general
predictors for $m=m'=1$. Treatment of the $(\bSigma_{rs},\sigma^2)$ entries for $r\ne s$
is similar and also leads to $O_P(n^{-1})$ leading term behavior
under assumption (A2).

For the $(m,m')=(2,1)$ case, with assistance from Lemmas \ref{lem:keyIdents}
and \ref{lem:keySqIdent}, $2n_{11}$ multiplied by the $(i,\iunder)=(1,1)$ term 
on the right-hand side of (\ref{eq:coffeeWithCannoli}) equals
{\setlength\arraycolsep{1pt}
\begin{eqnarray*}
&&n_{11}\be_r^T\bXAoneone^T(\bV^{11})^2\bXAoneone\be_r\\
&&\qquad=n_{11}\be_r^T\bXAoneone^T
\Bigg(\mbox{upper left $n_{11}\times n_{11}$ block of}\\
&&\qquad\qquad\qquad\qquad\quad \left[
\begin{array}{cc}
\bXAoneone(\bSigma+\bSigma')\bXAoneone^T+\sigma^2\bI& \bXAoneone\bSigma'\bXAtwoone^T \\[1ex]
\bXAtwoone\bSigma'\bXAoneone^T & \bXAtwoone(\bSigma+\bSigma')\bXAtwoone^T+\sigma^2\bI
\end{array}
\right]^{-1}\Bigg)^2\bXAoneone\be_r\\[1ex]
&&\qquad=n_{11}\be_r^T\bXAoneone^T
\Big[\bXAoneone(\bSigma+\bSigma')\bXAoneone^T+\sigma^2\bI\\
&&\qquad\qquad\qquad-\bXAoneone\bSigma'\bXAtwoone^T 
\big\{\bXAtwoone(\bSigma+\bSigma')\bXAtwoone^T+\sigma^2\bI\big\}^{-1}
\bXAtwoone\bSigma'\bXAoneone^T\Big]^{-2}
\bXAoneone\be_r\\[1ex]
&&\qquad=n_{11}\be_r^T\bXAoneone^T
\Big[\bXAoneone(\bSigma+\bSigma')\bXAoneone^T+\sigma^2\bI\\
&&\qquad\qquad\qquad\qquad-\bXAoneone\bSigma'
\big\{\bSigma+\bSigma'+\sigma^2(\bXAtwoone^T\bXAtwoone)^{-1}\bI\big\}^{-1}\bSigma'\bXAoneone^T\Big]^{-2}
\bXAoneone\be_r\\[1ex]
&&\qquad=\be_r^T\Big[\bSigma+\bSigma'-\bSigma'\big\{\bSigma+\bSigma'+\sigma^2(\bXAtwoone^T\bXAtwoone)^{-1}
\big\}^{-1}\bSigma'+\sigma^2(\bXAoneone^T\bXAoneone)^{-1}\Big]^{-1}\\[1ex]
&&\qquad\qquad\times\left(\textstyle{\frac{1}{n_{11}}}\bXAoneone^T\bXAoneone\right)^{-1}\\[1ex]
&&\qquad\qquad\times\Big[\bSigma+\bSigma'-\bSigma'\big\{\bSigma+\bSigma'+\sigma^2(\bXAtwoone^T\bXAtwoone)^{-1}
\big\}^{-1}\bSigma'+\sigma^2(\bXAoneone^T\bXAoneone)^{-1}\Big]^{-1}\be_r\\[1ex]
&&\qquad\convprob
\be_r^T\Big\{\bSigma+\bSigma'-\bSigma'\big(\bSigma+\bSigma'\big)^{-1}\bSigma'\Big\}^{-1}
E(\bXTrv\bXrv)\Big\{\bSigma+\bSigma'-\bSigma'\big(\bSigma+\bSigma'\big)^{-1}\bSigma'\Big\}^{-1}\be_r.
\end{eqnarray*}
}
Similar arguments lead to 
$$n_{21}\be_r^T\bXAtwoone^T(\bV^{22})^2\bXAtwoone\be_r$$
having the same convergence in probability limit.
In addition, and again using Lemmas \ref{lem:keyIdents}
and \ref{lem:keySqIdent},
{\setlength\arraycolsep{1pt}
\begin{eqnarray*}
&&n_{21}\be_r^T\bXAoneone^T(\bV^{21})^T\bV^{21}\bXAoneone\be_r\\
&&\qquad=n_{21}\be_r^T\bXAoneone^T
\Bigg(\mbox{the transposed lower left $n_{12}\times n_{11}$ block of}\\
&&\qquad\qquad\qquad\qquad\quad \left[
\begin{array}{cc}
\bXAoneone(\bSigma+\bSigma')\bXAoneone^T+\sigma^2\bI& \bXAoneone\bSigma'\bXAtwoone^T \\[1ex]
\bXAtwoone\bSigma'\bXAoneone^T & \bXAtwoone(\bSigma+\bSigma')\bXAtwoone^T+\sigma^2\bI
\end{array}
\right]^{-1}\Bigg)^{\otimes2}\bXAtwoone\be_r\\[1ex]
&&\qquad=n_{21}\be_r^T\bXAoneone^T
\big\{\bXAoneone(\bSigma+\bSigma')\bXAoneone^T+\sigma^2\bI\big\}^{-1}
\bXAoneone\bSigma'\bXAtwoone^T
\\
&&\qquad\qquad\qquad\times\Big[\bXAtwoone(\bSigma+\bSigma')\bXAtwoone^T+\sigma^2\bI\\
&&\qquad\qquad\qquad\qquad\qquad-\bXAtwoone\bSigma'\bXAoneone^T 
\big\{\bXAoneone(\bSigma+\bSigma')\bXAoneone^T+\sigma^2\bI\big\}^{-1}
\bXAoneone\bSigma'\bXAtwoone^T\Big]^{-2}\\
&&\qquad\qquad\qquad\times \bXAtwoone\bSigma'\bXAoneone^T
\big\{\bXAoneone(\bSigma+\bSigma')\bXAoneone^T+\sigma^2\bI\big\}^{-1}\bXAoneone\be_r\\[1ex]
&&\qquad=\be_r^T
\big\{\bSigma+\bSigma'+\sigma^2(\bXAoneone^T\bXAoneone)^{-1}\big\}^{-1}\bSigma'
\\
&&\qquad\qquad\qquad\times\Big[\bSigma+\bSigma'+\sigma^2(\bXAtwoone^T\bXAtwoone)^{-1}
-\bSigma'\big\{\bSigma+\bSigma'+\sigma^2(\bXAoneone^T\bXAoneone)^{-1}\big\}^{-1}
\bSigma'\Big]^{-1}\\
&&\qquad\qquad\qquad\times\big(\textstyle{\frac{1}{n_{21}}}\bXAtwoone^T\bXAtwoone\big)^{-1}\\
&&\qquad\qquad\qquad\times\Big[\bSigma+\bSigma'+\sigma^2(\bXAtwoone^T\bXAtwoone)^{-1}
-\bSigma'\big\{\bSigma+\bSigma'+\sigma^2(\bXAoneone^T\bXAoneone)^{-1}\big\}^{-1}
\bSigma'\Big]^{-1}\\
&&\qquad\qquad\qquad\times\bSigma'
\big\{\bSigma+\bSigma'+\sigma^2(\bXAoneone^T\bXAoneone)^{-1}\big\}^{-1}\be_r\\[1ex]
&&\qquad\convprob 
\be_r^T
\big(\bSigma+\bSigma'\big)^{-1}\bSigma'\Big\{\bSigma+\bSigma'
-\bSigma'\big(\bSigma+\bSigma'\big)^{-1}\bSigma'\Big\}^{-1}
E(\bXTrv\bXrv)\\
&&\qquad\qquad\qquad\times\Big\{\bSigma+\bSigma'
-\bSigma'\big(\bSigma+\bSigma'\big)^{-1}\bSigma'\Big\}^{-1}\bSigma'
\big(\bSigma+\bSigma'\big)^{-1}\be_r
\end{eqnarray*}
}
Similar steps lead to $n_{11}\be_r^T\bXAtwoone^T(\bV^{12})^T\bV^{12}\bXAtwoone\be_r$
having the same convergence in probability limit. On combining these results we
obtain the $(\bSigma_{rr},\sigma^2)$ entry of the Fisher information for $(m,m')=(2,1)$ 
having leading term behavior:
{\setlength\arraycolsep{1pt}
\begin{eqnarray*}
&&\smhalf\left(\textstyle{\frac{1}{n_{11}}}+\textstyle{\frac{1}{n_{21}}}\right)
\Bigg(\Big\{\bSigma+\bSigma'-\bSigma'\big(\bSigma+\bSigma'\big)^{-1}\bSigma'\Big\}^{-1}\\[1ex]
&&\qquad\qquad
\times E(\bXTrv\bXrv)\Big\{\bSigma+\bSigma'-\bSigma'\big(\bSigma+\bSigma'\big)^{-1}\bSigma'\Big\}^{-1}
\Bigg)_{rr}\{1+o_P(1)\}\\[1ex]
&&\quad+\smhalf\left(\textstyle{\frac{1}{n_{11}}}+\textstyle{\frac{1}{n_{21}}}\right)
\Bigg(\big(\bSigma+\bSigma'\big)^{-1}\bSigma'\Big\{\bSigma+\bSigma'
-\bSigma'\big(\bSigma+\bSigma'\big)^{-1}\bSigma'\Big\}^{-1}E(\bXTrv\bXrv)\\
&&\qquad\qquad\times\Big\{\bSigma+\bSigma'
-\bSigma'\big(\bSigma+\bSigma'\big)^{-1}\bSigma'\Big\}^{-1}\bSigma'
\big(\bSigma+\bSigma'\big)^{-1}\Bigg)_{rr}\{1+o_P(1)\}
\end{eqnarray*}
}
which, under assumption (A2), has $O_P(n^{-1})$ leading term behavior.
Similar arguments lead to the $O_P(n^{-1})$ property holding for 
the $\big(\bSigma_{rs},\sigma^2\big)$ entries of the
Fisher information matrix for $r\ne s$ when $(m,m')=(2,1)$.

For higher $m$ and $m'$, similar arguments can be used to show that the 
summations in $\big(\vech(\bSigma),\sigma^2\big)$ Fisher information block
lead to convergents that are analogous to those in the $\dA=1$, $n_{ii'}=n$ and 
$\bXAiid=\bone_n$ case and the block satisfies 
$O_P(n^{-1})\bone_{\dA(\dA+1)/2}$ under (A1) and (A2).

This very low order of magnitude of the $\big(\vech(\bSigma),\sigsq\big)$ 
off-diagonal block of the Fisher information matrix is more than 
enough for asymptotic orthogonality between $\bSigma$ and $\sigma^2$.
A larger order of magnitude, such as $O_P(1)\bone_{\dA(\dA+1)/2}$,
would still be sufficient.

\subsubsection{The $\big(\vech(\bSigma'),\sigsq\big)$ Off-Diagonal Block}

In the special case of $\dA=1$, $n_{ii'}=n$ and 
$\bXAiid=\bone_n$ for all $1\le i\le m$, $1\le i'\le m'$
use of the eigenvalue and eigenvector properties described
near the commencement of Section \ref{sec:firstEig}
lead to the $(\Sigma',\sigma^2)$ entry of the Fisher information matrix
having exact expression
$$
{\setlength\arraycolsep{1pt}
\begin{array}{l}
\displaystyle{
\frac{m(1-1/m')\{\Sigma(m'/m)+\Sigma'+\sigma^2/(mn)\}^2}{2m'n
\{\Sigma+\Sigma'(m/m')+\sigma^2/(m'n)\}^2\{\Sigma'+\sigma^2/(mn)\}^2}}\\[3ex]
\qquad\qquad\qquad\qquad\qquad
+\displaystyle{\frac{m}{2(m')^2n\{\Sigma+\Sigma'(m/m')+\sigma^2/(m'n)\}^2}}
\end{array}
}
$$
which has the same form as (\ref{eq:AroundNearPole}) but with the roles of $(M,m)$ and $(M',m')$
reversed. Symmetry considerations dictate that the same happens in the general
setting and the $\big(\vech(\bSigma'),\sigsq\big)$ off-diagonal block is 
$O(n^{-1})\bone_{\dA/(\dA+1)/2}$.

\subsubsection{Assembly of the Fisher Information Sub-Block Approximations}

The Fisher information sub-block approximations obtained in the previous
nine sub-subsections lead to 
{\setlength\arraycolsep{0pt}
\begin{eqnarray*}
&&I\big(\bbetaA,\bbetaB,\vech(\bSigma),\vech(\bSigma'),\sigsq\big)=\\[1ex]
&&\left[
{\setlength\arraycolsep{0.1pt}
\begin{array}{ccccc}
{\displaystyle\left(\frac{\bSigma}{m}+\frac{\bSigma'}{m'}\right)^{-1}}
& O_P(m)\bone_{\dA}\bone_{\dB}^T & \bO & \bO &  \bO  \\[1ex]
\ \ +o_P(m)\bone_{\dA}^{\otimes2} &   &   &   &     \\[2ex]
O_P(m)\bone_{\dB}\bone_{\dA}^T  & {\displaystyle\frac{mm'n\COMPQUANT_{\bbetaB}^{-1}}{\sigma^2}} 
& \bO & \bO &  \bO  \\[1ex]
  & \ \ +o_P(m^2n)\bone_{\dB}^{\otimes2} &   &   &     \\[2ex]
\bO & \bO & {\displaystyle\frac{m\bD_{\dA}^T(\bSigma^{-1}\otimes\bSigma^{-1})\bD_{\dA}}{2}} 
& O_P(1)\bone_{\dAwindow}^{\otimes2} &  O_P(n^{-1})\bone_{\dAwindow}  \\[1ex]
  &   & \ \ +o_P(m)\bone_{\dAwindow}^{\otimes2} &   &    \\[2ex]
\bO & \bO & O_P(1)\bone_{\dAwindow}^{\otimes2}  
& {\displaystyle\frac{m'\bD_{\dA}^T\big((\bSigma')^{-1}\otimes(\bSigma')^{-1}\big)\bD_{\dA}}{2}} 
&  O_P(n^{-1})\bone_{\dAwindow}  \\[1ex]
  &   &  & \ \ +o_P(m)\bone_{\dAwindow}^{\otimes2} &    \\[2ex]
\bO & \bO & O_P(n^{-1})\bone_{\dAwindow}^T 
& O_P(n^{-1})\bone_{\dAwindow}^T &  {\displaystyle\frac{mm'n}{2\sigma^4}}\\[1ex]
  &   &  &  & \ \ +o_P(m^2n) \\[1ex]
\end{array}
}
\right]
\end{eqnarray*}
}
where $\dAwindow\equiv\smhalf\dA(\dA+1)$.

\subsection{Inverse Fisher Information Matrix Approximation}

First note that, since 
$I\big(\bbetaA,\bbetaB,\vech(\bSigma),\vech(\bSigma'),\sigsq\big)$
is block diagonal, its inversion involves the individual 
inversions of the $\big(\bbetaA,\bbetaB\big)$ and 
$\big(\vech(\bSigma),\vech(\bSigma'),\sigsq\big)$ blocks.
These two inversions involve application of well-known
block matrix inversion formulae and keeping track of the
various terms that arise and their orders of magnitude. 
For example, if the sub-blocks of the $\big(\bbetaA,\bbetaB\big)$ block
are denoted as follows: 
$$
\left[
\begin{array}{cc}
\bA_{11}   &  \bA_{12} \\[1ex]
\bA_{12}^T &  \bA_{22}
\end{array}
\right]
\quad\mbox{where}\quad\bA_{11}\ \mbox{is $\dA\times\dA$}
$$
then the upper left $\dA\times\dA$ block of the required
inverse matrix is 
$$\bA_{11}^{-1}+\bA_{11}^{-1}\bA_{12}(\bA_{22}-\bA_{12}^T
\bA_{11}^{-1}\bA_{12})^{-1}\bA_{12}^T\bA_{11}^{-1}.
$$
Appendix A.6 of Jiang \textit{et al.} (2022) contains a detailed
account of this approach for related setting. Analogous steps
for the current setting lead to 
{\setlength\arraycolsep{1pt}
\begin{eqnarray*}
&&I\big(\bbetaA,\bbetaB,\vech(\bSigma),\vech(\bSigma'),\sigsq\big)^{-1}
=I\big(\bbetaA,\bbetaB,\vech(\bSigma),\vech(\bSigma'),\sigsq\big)_{\infty}^{-1}\\[1ex]
&&\ \ +\frac{1}{m}
\left[
{\setlength\arraycolsep{0.5pt}
\begin{array}{ccccc}
o_P(1)\bone_{\dA}^{\otimes2} & O_P(m^{-1}n^{-1})\bone_{\dA}\bone_{\dB}^T  &  \bO & \bO & \bO \\[1ex]
O_P(m^{-1}n^{-1})\bone_{\dB}\bone_{\dA}^T & o_P(m^{-1}n^{-1})\bone_{\dB}^{\otimes2} & \bO & \bO & \bO \\[1ex]
\bO & \bO & o_P(1)\bone_{\dAwindow}^{\otimes2}  
& O_P(m^{-1})\bone_{\dAwindow}\bone_{\dAwindow}^{\otimes2}  & O_P(m^{-2}n^{-1})\bone_{\dAwindow} \\[2ex]
\bO & \bO & O_P(m^{-1})\bone_{\dAwindow}\bone_{\dAwindow}^{\otimes2}  
& o_P(1)\bone_{\dAwindow}^{\otimes2} & O_P(m^{-2}n^{-1})\bone_{\dAwindow}  \\[2ex]
\bO & \bO & O_P(m^{-2}n^{-1})\bone^T_{\dAwindow} & O_P(m^{-2}n^{-1})\bone^T_{\dAwindow} &  o_P(m^{-2}n^{-1})
\end{array}
}
\right]
\end{eqnarray*}
}
where
{\setlength\arraycolsep{1pt}
\begin{eqnarray*}
&&I\big(\bbetaA,\bbetaB,\vech(\bSigma),\vech(\bSigma'),\sigsq\big)_{\infty}^{-1}\\[1ex]
&&\qquad=\left[
{\setlength\arraycolsep{2pt}
\begin{array}{ccccc}
{\displaystyle\frac{\bSigma}{m}+\frac{\bSigma'}{m'}} & \bO & \bO & \bO & \bO \\[2ex]
\bO & {\displaystyle\frac{\sigsq\COMPQUANT_{\bbetaB}}{mm'n}} & \bO & \bO & \bO  \\[2ex]
\bO & \bO & {\displaystyle{\frac{2\bD_{\dA}^+(\bSigma\otimes\bSigma)\bD_{\dA}^{+T}}{m}}} 
& \bO & \bO  \\[2ex]
\bO & \bO & \bO & {\displaystyle{\frac{2\bD_{\dA}^+(\bSigma'\otimes\bSigma')\bD_{\dA}^{+T}}{m'}}}  & \bO  \\[2ex]
\bO & \bO & \bO & \bO & {\displaystyle{\frac{2\sigma^4}{mm'n}}}  \\[2ex]
\end{array}
}
\right].
\end{eqnarray*}
}

\subsection{Asymptotic Normality of the Maximum Likelihood Estimators}

Let 
$$\bbeta\equiv(\bbetaA,\bbetaB)\quad\mbox{and}\quad
\bpsi\equiv\big(\vech(\bSigma),\vech(\bSigma'),\sigsq\big).
$$
As alluded to in Section \ref{sec:exactOrthog}, the Fisher information
has the block diagonal form:
\begin{equation}
I(\bbeta,\bpsi)=\left[
\begin{array}{cc}
I(\bbeta,\bpsi)_{\smbbeta\smbbeta} & \bO \\[2ex]
\bO                            & I(\bbeta,\bpsi)_{\smbpsi\smbpsi} 
\end{array}
\right].
\label{eq:ArchieBunker}
\end{equation}
where $I(\bbeta,\bpsi)_{\smbbeta\smbbeta}$ is the upper left $(\dA+\dB)\times(\dA+\dB)$
block of $I(\bbeta,\bpsi)$ and $I(\bbeta,\bpsi)_{\smbpsi\smbpsi}$ is defined similarly. 
Then, under (A1)--(A3) and some additional regularity conditions
\begin{equation}
\{I(\bbeta^0,\bpsi^0)^{-1}\}^{-1/2}
\left[
\begin{array}{c}
\bbetahat-\bbetaZero\\[1ex]
\bpsihat-\bpsiZero
\end{array}
\right]\convdist N(\bzero,\bI).
\label{eq:GoldCoastAngela}
\end{equation}
Justification for (\ref{eq:GoldCoastAngela}) is given in Section \ref{sec:martCLT}.

\subsection{Convergence Results for Matrix Square Root Discrepancies}

We now deal with the problem of proving that matrix square roots of 
the exact inverse Fisher information matrix and its convergent
$$\{I\big(\bbetaA,\bbetaB,\vech(\bSigma),\vech(\bSigma'),\sigsq\big)^{-1}\}^{1/2}
\quad\mbox{and}\quad
\{I\big(\bbetaA,\bbetaB,\vech(\bSigma),\vech(\bSigma'),\sigsq\big)_{\infty}^{-1}\}^{1/2}
$$
are also sufficiently close to each other as $m$, $m'$ and $n$ diverge.
Using the notation from (\ref{eq:ArchieBunker}), we treat the fixed effects
and covariance parameter diagonal blocks separately. To this end, define
$$I(\bbeta,\bpsi)^{-1}_{\smbbeta\smbbeta,\infty}\equiv
\left[
{\setlength\arraycolsep{0.5pt}
\begin{array}{cc}
{\displaystyle\frac{\bSigma}{m}+\frac{\bSigma'}{m'}} & \bO \\[1ex]
\bO & {\displaystyle\frac{\sigsq\COMPQUANT_{\bbetaB}}{mm'n}} 
\end{array}
}
\right]
$$
and
$$
I(\bbeta,\bpsi)^{-1}_{\smbpsi\smbpsi,\infty}\equiv
\left[
{\setlength\arraycolsep{2pt}
\begin{array}{ccc}
{\displaystyle{\frac{2\bD_{\dA}^+(\bSigma\otimes\bSigma)\bD_{\dA}^{+T}}{m}}} 
& \bO & \bO  \\[2ex]
\bO & {\displaystyle{\frac{2\bD_{\dA}^+(\bSigma'\otimes\bSigma')\bD_{\dA}^{+T}}{m'}}}  & \bO  \\[2ex]
\bO & \bO & {\displaystyle{\frac{2\sigma^4}{mm'n}}}  \\[2ex]
\end{array}
}
\right].
$$
Next note that 
$$m'I(\bbeta,\bpsi)^{-1}_{\smbbeta\smbbeta}
=\left[
{\setlength\arraycolsep{0.5pt}
\begin{array}{cc}
\bK + o_P(\bone_{\dA}^{\otimes2}) & O_P\big((mn)^{-1}\big)\bone_{\dA}\bone_{\dB}^T \\[1ex]
O_P\big((mn)^{-1}\big)\bone_{\dB}\bone_{\dA}^T & \frac{1}{m}\bL + o_P\big((mn)^{-1}\big)\bone_{\dB}^{\otimes2}
\end{array}
}
\right]
\ \mbox{and}\ 
m'I(\bbeta,\bpsi)^{-1}_{\smbbeta\smbbeta,\infty}
=\left[
{\setlength\arraycolsep{0.5pt}
\begin{array}{cc}
\bK & \bO \\[1ex]
\bO & \frac{1}{m}\bL
\end{array}
}
\right]
$$
where
$$\bK\equiv (m'/m)\bSigma + \bSigma'\quad\mbox{and}\quad\bL\equiv
\frac{\sigsq\COMPQUANT_{\bbetaB}}{n}.
$$
Then application of Lemma 2 of Jiang \textit{et al.} (2022) as $m\to\infty$ implies that
\begin{equation}
\Big\Vert\{I\big(\bbeta,\bpsi\big)^{-1}_{\smbbeta\smbbeta,\infty}\}^{-1/2}
\{I\big(\bbeta,\bpsi\big)^{-1}_{\smbbeta\smbbeta}\}^{1/2}-\bI\BigVertF\convprob 0.
\label{eq:SauceBrewery}
\end{equation}
The establishment
\begin{equation}
\Big\Vert\{I\big(\bbeta,\bpsi\big)^{-1}_{\smbpsi\smbpsi,\infty}\}^{-1/2}
\{I\big(\bbeta,\bpsi\big)^{-1}_{\smbpsi\smbpsi}\}^{1/2}-\bI\BigVertF\convprob 0.
\label{eq:MarrickvilleManBuns}
\end{equation}
is very similar.

\subsection{Final Steps for the Derivation of \QuasiTheorem\ \ref{res:mainResult}}\label{sec:ThePenguin}

Let 
$$\btheta\equiv(\bbeta,\bpsi)=\big(\bbetaA,\bbetaB,\vech(\bSigma),\vech(\bSigma'),\sigsq\big)$$
be the full parameter vector. In terms of this new notation, 
result (\ref{eq:GoldCoastAngela}) is
\begin{equation}
\big\{I\big(\bthetaZero\big)^{-1}\big\}^{-1/2}(\bthetaMLE-\bthetaZero)\convdist N(\bzero,\bI)
\label{eq:CLTraw}
\end{equation}
where 
$$\bthetaMLE=\big[\bbetaAMLET\ \bbetaBMLET\ \vech(\bSigmaMLE)^T\ \vech(\bSigmadMLE)^T\ \sigsqMLE\big]^T$$ 
and
$$\btheta^0=\big[(\bbetaAzero)^T\ (\bbetaBzero)^T\ \vech(\bSigmaZero)^T\ \vech(\bSigmaZero)^T
\ \vech(\bSigmadZero)^T\ \sigsqZero\big]^T.$$
It follows from (\ref{eq:CLTraw}) that, for all $(\dA+\dB+2\dAwindow+1)\times 1$ 
vectors $\ba\ne\bzero$, we have
$$
\ba^T\,\{I\big(\bthetaZero\big)^{-1}\}^{-1/2}(\bthetaMLE-\bthetaZero)
\convdist N(0,\ba^T\ba).
$$
As a consequence
\begin{equation}
\ba^T\{I\big(\bthetaZero\big)_{\infty}^{-1}\}^{-1/2}
(\bthetaMLE-\bthetaZero)+r_{mm'n}(\ba)\convdist N(0,\ba^T\ba)
\label{eq:preSlut}
\end{equation}
where
{\setlength\arraycolsep{1pt}
\begin{eqnarray*}
r_{mm'n}(\ba)&\equiv&\ba^T[\{I\big(\bthetaZero\big)^{-1}\}^{-1/2}
-\{I\big(\bthetaZero\big)_{\infty}^{-1}\}^{-1/2}](\bthetaMLE-\bthetaZero)\\[1ex]
&=&\ba^T[\bI-\{I\big(\bthetaZero\big)_{\infty}^{-1}\}^{-1/2}
\{I\big(\bthetaZero\big)^{-1}\}^{1/2}]
\{I\big(\bthetaZero\big)^{-1}\}^{1/2}(\bthetaMLE-\bthetaZero)\\[1ex]
&=&
\left(\big[\{I\big(\bthetaZero\big)_{\infty}^{-1}\}^{-1/2}
\{I\big(\bthetaZero\big)^{-1}\}^{1/2}-\bI\big]^T\ba\right)^T\bZ_{mm'n}
\end{eqnarray*}
}
and $\bZ_{mm'n}\convdist N(\bzero,\bI_{\dA+\dB+2\dAwindow+1})$. Then note that 
$$\Big\Vert\big[\{I\big(\bthetaZero\big)_{\infty}^{-1}\}^{-1/2}
\{I\big(\bthetaZero\big)^{-1}\}^{1/2}-\bI\big]^T\ba\BigVertF
\le 
\big\Vert\{I\big(\bthetaZero\big)_{\infty}^{-1}\}^{-1/2}
\{I\big(\bthetaZero\big)^{-1}\}^{1/2}-\bI\bigVertF\,\Vert\ba\VertF.
$$
As a consequence of (\ref{eq:SauceBrewery}) and (\ref{eq:MarrickvilleManBuns})
we have
\begin{equation}
\big\Vert\{I\big(\bthetaZero\big)_{\infty}^{-1}\}^{-1/2}
\{I\big(\bthetaZero\big)^{-1}\}^{1/2}-\bI\VertF\convprob 0
\label{eq:matrixNormRes}
\end{equation}
and so
$$\big[\{I\big(\bthetaZero\big)_{\infty}^{-1}\}^{-1/2}
\{I\big(\bthetaZero\big)^{-1}\}^{1/2}-\bI\big]\ba\convprob 0.$$
Application of Slutsky's Theorem then gives $r_{mm'n}(\ba)\convprob 0$. From (\ref{eq:preSlut}) 
and another application of Slutsky's Theorem we have 
$$
\ba^T\,\{I\big(\bthetaZero\big)_{\infty}^{-1}\}^{-1/2}(\bthetaMLE-\bthetaZero)
\convdist N(0,\ba^T\ba).
$$
\QuasiTheorem\ \ref{res:mainResult} then follows from the Cram\'er-Wold Device.

\subsection{Justification of (\ref{eq:GoldCoastAngela})}\label{sec:martCLT}

We now provide justification for the asymptotic
normality statement (\ref{eq:GoldCoastAngela})
concerning the maximum likelihood estimators and
the Fisher information matrix.

As in Section \ref{sec:ThePenguin} we let 
$$\btheta\equiv(\bbeta,\bpsi)=\big(\bbetaA,\bbetaB,\vech(\bSigma),\vech(\bSigma'),\sigsq\big)$$
be the full parameter vector. The score vector is 
$$\nabla_{\btheta}\ell(\btheta)
=\left[
\begin{array}{c}
\bXA^T\bV^{-1}(\bY-\bXA\bbetaA-\bXB\bbetaB)\\[1ex]
\bXB^T\bV^{-1}(\bY-\bXA\bbetaA-\bXB\bbetaB)\\[1ex]
\smhalf{\displaystyle\stack{(r,s)\in\Isc_{\dA}}}\left\{
\tr\Big(\bV^{-1}
\bLbreve_{(r,s)}\bV^{-1}(\bY-\bXA\bbetaA-\bXB\bbetaB)^{\otimes2}-\bV^{-1}\bLbreve_{(r,s)}\Big)
\right\}\\[2ex]
\smhalf{\displaystyle\stack{(r,s)\in\Isc_{\dA}}}\left\{
\tr\Big(\bV^{-1}
\bLbrevedash_{(r,s)}\bV^{-1}(\bY-\bXA\bbetaA-\bXB\bbetaB)^{\otimes2}
-\bV^{-1} \bLbrevedash_{(r,s)}\Big)\right\}\\[2ex]
\smhalf\tr\Big(\bV^{-2}(\bY-\bXA\bbetaA-\bXB\bbetaB)^{\otimes2}
-\bV^{-1}\Big)
\end{array}
\right]
$$
where
$$\Isc_{\dA}\equiv\{(1,1),(2,1),\ldots,(\dA,1),(2,2),(3,2),\ldots,(\dA,2),\ldots,(\dA,\dA)\}$$
corresponds to positions on and below the diagonal
of a $\dA\times\dA$ matrix with the $\vech$ operator ordering,
$$\bLbreve_{(r,s)}\equiv\bL_r\bL_s^T+I(r\ne s)\bL_s\bL_r^T$$
with $\bL_r$ as defined by (\ref{eq:DeeDeeMyers}), and
$$\bLbreve'_{(r,s)}\equiv 
{\displaystyle\blockmatrix{1\le i,\iunder\le m}
\left\{\blockdiag{1\le i'\le m'}\Big(\bXAiid\Big(\be_r\be_s^T+I(r\ne s)\be_s\be_r^T\Big)\bXAiuid^T\Big)
\right\}}.$$
Let 
$$\bZ\equiv\left[\blockdiag{1\le i\le m}\left\{\stack{1\le i'\le m'}(\bXAiid)\right\}
\ \ \stack{1\le i\le m}\left\{\blockdiag{1\le i'\le m'}(\bXAiid)\right\}\right],$$
$$\bUall\equiv\left[
\begin{array}{c}
\displaystyle{\stack{1\le i\le m}}(\bU_i)\\[1ex]
\displaystyle{\stack{1\le i'\le m'}}(\bU'_{i'})
\end{array}
\right]
\quad\mbox{and}\quad
\bG\equiv\left[
\begin{array}{cc}
\bI_m\otimes\bSigma & \bO \\[2ex]
\bO                 & \bI_{m'}\otimes\bSigma'
\end{array}
\right].
$$
Next, define
$$\bz\equiv \left[
\begin{array}{cc}
\bG & \bO     \\[1ex]
\bO       & \sigma^2\bI
\end{array}
\right]^{-1/2}
\left[
\begin{array}{c}
\bUall\\[1ex]
\bY-\bXA\bbetaA-\bXB\bbetaB-\bZ\bUall
\end{array}
\right]\ \mbox{and}\ 
\bVhalfLoose\equiv[\bZ\ \bI]\left[
\begin{array}{cc}
\bG & \bO     \\[1ex]
\bO       & \sigma^2\bI
\end{array}
\right]^{1/2}.
$$
The relationship
$$\bVhalfLoose\big(\bVhalfLoose\big)^T=\bV$$
is the reason for the $\bVhalfLoose$ notation since,
loosely (i.e.\ ignoring transposes), it is a matrix square root of $\bV$.
Noting that 
$$\bY-\bXA\bbetaA-\bXB\bbetaB=\bVhalfLoose\bz$$
we can re-write the score vector as
$$\nabla_{\btheta}\ell(\btheta)
=\left[
\begin{array}{c}
\displaystyle{\stack{1\le r\le\dA}}(\bwAr^T\bz)\\[2ex]
\displaystyle{\stack{1\le r\le\dB}}(\bwBr^T\bz)\\[2ex]
\smhalf{\displaystyle\stack{(r,s)\in\Isc_{\dA}}}\left\{
\tr\big(\bW_{(r,s)}(\bz^{\otimes2}-\bI)\big)\right\}\\[2ex]
\smhalf{\displaystyle\stack{(r,s)\in\Isc_{\dA}}}\left\{
\tr\big(\bW'_{(r,s)}(\bz^{\otimes2}-\bI)\big)\right\}\\[2ex]
\smhalf\tr\big(\bW_{\sigma^2}(\bz^{\otimes2}-\bI)\big)
\end{array}
\right]
$$
where
{\setlength\arraycolsep{1pt}
\begin{eqnarray*}
\bwAr&\equiv&\mbox{$r$th column of}\ (\bVhalfLoose)^T
\bV^{-1}\bXA,\quad 1\le r\le \dA,\\[1ex]
\bwBr&\equiv&\mbox{$r$th column of}\ (\bVhalfLoose)^T
\bV^{-1}\bXB,\quad 1\le r\le \dB,\\[1ex]
\bW_{(r,s)}&=&(\bVhalfLoose)^T\bV^{-1}\bLbreve_{(r,s)}\bV^{-1}\bVhalfLoose,\quad (r,s)\in\Isc_{\dA},\\[1ex]
\bW'_{(r,s)}&=&(\bVhalfLoose)^T\bV^{-1}\bLbrevedash_{(r,s)}\bV^{-1}\bVhalfLoose,\quad (r,s)\in\Isc_{\dA}\\[1ex]
\mbox{and}\quad 
\bW_{\sigma^2}&=&(\bVhalfLoose)^T\bV^{-2}\bVhalfLoose.
\end{eqnarray*}
}
Let
$$\bs(m,n)\equiv\left[m\bone_{\dA}\ \ m^2n\bone_{\dB}\ \ m\bone_{\smhalf\dA(\dA+1)}\ \  
m\bone_{\smhalf\dB(\dB+1)}\ \ m^2n\right]^T
$$
be a vector of sample size quantities and accounts for the $m=O(m')$ and $m'=O(m)$ assumptions. 
Then define
$$\banorm\equiv \diag\{\bs(m,n)\}^{1/2}I\big(\bthetaZero\big)^{-1/2}\ba.$$
Letting $\bn$ denote the matrix of $n_{ii'}$ values, note that
$$\ba^TI\big(\bthetaZero\big)^{-1/2}
\nabla_{\btheta}\ell(\bthetaZero)
=\banorm^T\diag\{\bs(m,n)\}^{-1/2}\nabla_{\btheta}\ell(\bthetaZero)
=\sum_{t=1}^{\Nmart}\xi_t(m,m',\bn)
$$
where, for $1\le t\le\Nmart$, 
{\setlength\arraycolsep{1pt}
\begin{eqnarray*}
\xi_t(m,m',\bn)&\equiv&(\banorm)_1 m^{-1/2}(\bwAone^0)_t(\bz)_t
+\ldots+(\banorm)_{\dA}m^{-1/2}(\bwAdA^0)_t(\bz)_t\\[1ex]
&&\quad+(\banorm)_{\dA+1} (m^2n)^{-1/2}(\bwBone^0)_t(\bz)_t
+\ldots+(\banorm)_{\dA+\dB}(m^2n)^{-1/2}(\bwBdB^0)_t(\bz)_t\\[1ex]
&&\quad+\smhalf(\banorm)_{\dA+\dB+1}m^{-1/2}\big(\bW^0_{(1,1)}(\bz^{\otimes2}-\bI)\big)_{tt}\\[1ex]
&&\quad
+\ldots+\smhalf(\banorm)_{\dA+\dB+\smhalf\dA(\dA+1)}m^{-1/2}\big(\bW^0_{(\dA,\dA)}(\bz^{\otimes2}-\bI)\big)_{tt}
\\[1ex]
&&\quad+\smhalf(\banorm)_{\dA+\dB+\smhalf\dA(\dA+1)+1}m^{-1/2}\big((\bW')^0_{(1,1)}(\bz^{\otimes2}-\bI)\big)_{tt}\\[1ex]
&&\quad
+\ldots+\smhalf(\banorm)_{\dA+\dB+\smhalf\dA(\dA+1)+\smhalf\dB(\dB+1)}
m^{-1/2}\big((\bW')^0_{(\dB,\dB)}(\bz^{\otimes2}-\bI)\big)_{tt}
\\[1ex]
&&\quad+\smhalf(\banorm)_{\dA+\dB+\smhalf\dA(\dA+1)+\smhalf\dB(\dB+1)+1} 
(m^2n)^{-1/2}\big(\bW^0_{\sigma^2}(\bz^{\otimes2}-\bI)\big)_{tt}
\end{eqnarray*}
}
and $\Nmart\equiv m+m'+\ndotdot$.
In the definition of the $\xi_t(m,m',\bn)$, the notation $\bwAr^0$ signifies that
each of the model parameters that appear in the definition of $\bwAr$ 
are set to their true values. A similar convention applies to the $\bwBr^0$,
$\bW^0_{(r,s)}$, $(\bW')^0_{(r,s)}$ and $\bW^0_{\sigma^2}$.
Let $\Xsc$ denote the full set of predictor random variables in $\bXA$
and $\bXB$. For $1\le t\le m$, let 
$$\Fsc_t(m,m'\bn)\ \mbox{denote the $\sigma$-field generated by $\Xsc,\bU_1,\ldots,\bU_t$}.$$
For $m\le t\le m+m'$, let 
$$\Fsc_t(m,m',\bn)\ \mbox{denote the $\sigma$-field generated by 
$\Xsc,\bU_1,\ldots,\bU_m,\bU'_1,\ldots,\bU'_t$}.$$
For $m+m'+1\le t\le\Nmart$, let 
\begin{eqnarray*}
&&\Fsc_t(m,m',\bn)\ \mbox{denote the $\sigma$-field generated by} 
\ \Xsc,\bU_1,\ldots,\bU_m,\bU'_1,\ldots,\bU'_{m'},\\
&&\qquad\qquad\qquad\qquad\qquad\qquad\qquad\qquad(\bY-\bXA\bbetaA^0-\bXB\bbetaB^0-\bZ\bUall)_{t-m-m'}.
\end{eqnarray*}
Then
$$\big(\xi_t(m,m',\bn), \Fsc_t(m,m',\bn)\big),\quad1\le t\le\Nmart,$$
is an array of martingale differences.

According to Theorem 3.2 of Hall \myand Heyde (1980), 
\begin{equation}
\ba^TI\big(\bthetaZero\big)^{-1/2}\nabla_{\btheta}\ell(\bthetaZero)
=\sum_{t=1}^{\Nmart}\xi_t(m,m',\bn)\convdist N(0,\ba^T\ba)
\label{eq:scoreCLT}
\end{equation}
if the $\xi_t(m,m',\bn)$ satisfy 
\begin{equation}
\begin{array}{l}
{\displaystyle\max_{1\le t\le\Nmart}}\Big|\xi_t(m,m',\bn)\Big|\convprob 0,\quad
{\displaystyle\sum_{t=1}^{\Nmart}}\xi_t(m,m',\bn)^2\convprob\ba^T\ba\\[3ex]
\quad\mbox{and}\quad
E\left({\displaystyle\max_{1\le t\le\Nmart}}\xi_t(m,m',\bn)^2\right)\ \mbox{is bounded in $(m,m',\bn)$}.
\end{array}
\label{eq:HandHconds}
\end{equation}
Arguments similar to those given in Jiang (1996) 
and Jiang \textit{et al.} (2023) can be used to establish
(\ref{eq:HandHconds}) under conditions such as (A1)--(A3). 
The pathway used in these references
involves studying the asymptotic behaviors of the norms
$$\Vert\bwAr\Vert^2=\big(\bXA^T\bV^{-1}\bXA\big)_{rr},\quad 1\le r\le\dA,
\qquad
\Vert\bwBr\Vert^2=\big(\bXB^T\bV^{-1}\bXB\big)_{rr},\quad 1\le r\le\dA,
$$
$$\Vert\bW_{(r,s)}\VertF^2=\tr\big((\bV^{-1}\bLbreve_{(r,s)})^2\big),
\ \ \Vert\bW'_{(r,s)}\VertF^2=\tr\big((\bV^{-1}\bLbrevedash_{(r,s)})^2\big),
\ \ \Vert\bW_{\sigma^2}\VertF^2=\tr\big(\bV^{-2}\big)$$
for $(r,s)\in\Isc_{\dA}$, as well as the maximum eigenvalues 
of the $\bW_{(r,s)}$, $\bW'_{(r,s)}$ and $\bW_{\sigma^2}$
matrices. From Section \ref{sec:FishInfMatAppr}, the
$\Vert\bwAr\Vert^2$, $\Vert\bW_{(r,s)}\VertF^2$
and $\Vert\bW'_{(r,s)}\VertF^2$ quantities are each $O_P(m)$
under assumption (A1). The $\Vert\bwBr\Vert^2$ and 
$\Vert\bW_{\sigma^2}\VertF^2$ quantities are 
$O_P(m^2n)$ under (A1). The maximum eigenvalue quantities
have similar asymptotic behaviors.

The conditions in (\ref{eq:HandHconds}) follow from results such as 
\begin{equation}
m^{-1}E\big(\Vert\bwAr\Vert^2\big)=O(1)\quad\mbox{and}\quad
(m^2n)^{-1}E\big(\Vert\bW_{\sigma^2}\VertF^2\big)=O(1).
\label{eq:FruitLoops}
\end{equation}
In the case of crossed random intercepts, these matrix norm expectations
follow quickly from the Section \ref{sec:FishInfMatAppr} results.
For the general crossed random effects model (\ref{eq:theModel})
the $\bV$ matrix is random and some additional regularity conditions
are required to ensure that expectations, such as those 
appearing in (\ref{eq:FruitLoops}), have the correct orders
of magnitude and, in turn, provide (\ref{eq:scoreCLT}).
Assuming these regularity conditions, 
the Cram\'er-Wold Device leads to 
$$I\big(\bthetaZero\big)^{-1/2}\nabla_{\btheta}\ell(\bthetaZero)\convdist N(\bzero,\bI).$$
Standard likelihood theory arguments then lead to 
$$I\big(\bthetaZero\big)^{1/2}(\bthetahat-\btheta^0)\convdist N(\bzero,\bI).$$

\subsection{Proofs of Lemmas}

The derivation of \QuasiTheorem\ \ref{res:mainResult} heavily depended on 
Lemmas \ref{lem:Kron}--\ref{lem:XBpart}. We now get to proving them.

\subsubsection{Proof of Lemma \ref{lem:Kron}}

Let  $\be_r$ denote the $d\times1$ matrix with $r$th entry $1$
and all other entries $0$. Then note that
$$A_{rs}=\vech(\be_r\be_s^T)^T\vech(\bA_d)\quad\mbox{for all $r\ge s$.}$$
Therefore
$$A_{rs}A_{tu}=
\vech(\be_r\be_s^T)^T\vech(\bA_d)\vech(\bA_d)^T \vech(\be_t\be_u^T)
\quad\mbox{for all $r\ge s$,\ $t\ge u$.}
$$
Next, note that
$$\bD_d\vech(\be_r\be_s^T)=\bD_d\vech\big(\be_r\be_s^T
+I(r\ne s)\be_s\be_r^T\big)
=\vecof\big(\be_r\be_s^T+I(r\ne s)\be_s\be_r^T\big)
\quad\mbox{for all $r\ge s$}.
$$
Use of the $\vecof(\ba\bb^T)=\bb\otimes\ba$ identity then gives
$$\bD_d\vech(\be_r\be_s^T)=\be_r\otimes\be_s+I(r\ne s)(\be_s\otimes\be_r).$$
We then have for all $r\ge s$ and $t\ge u$
{\setlength\arraycolsep{1pt}
\begin{eqnarray*}
\vech(\be_r\be_s^T)^T\bD_d^T(\bA\otimes\bA)\bD_d\vech(\be_t\be_u^T)
&=&\{\be_r\otimes\be_s+I(r\ne s)(\be_s\otimes\be_r)\}^T(\bA\otimes\bA)\\
&&\quad\times\{\be_t\otimes\be_u+I(t\ne u)(\be_u\otimes\be_t)\}\\[1ex]
&=&(\be_r^T\bA\be_t)(\be_s^T\bA\be_u)
+I(t\ne u)(\be_r^T\bA\be_u)(\be_s^T\bA\be_t)\\
&&\ +I(r\ne s)(\be_s^T\bA\be_t)(\be_r^T\bA\be_u)\\
&&\ +I(r\ne s)I(t\ne u)(\be_s^T\bA\be_u)(\be_r^T\bA\be_t)\\[1ex]
%
%
%
&=&\left\{
\begin{array}{ll}
A_{rt}^2,\quad &r=s,t=u,\\
2A_{rt}A_{ru},\quad &r=s,t>u,\\
2(A_{rt}A_{su}+A_{ru}A_{st}),\quad &r>s,t>u,
\end{array}
\right.\\[1ex]
&=&\vech(\be_r\be_s^T)^T\bB_d\vech(\be_t\be_u^T)
\end{eqnarray*}
}
Therefore the $r\ge s$ and $t\ge u$ entries of  
$\bB_d$ match those of $\bD_d^T(\bA\otimes\bA)\bD_d$.
However, if the roles of $r$ and $s$ are reversed
then each of the expressions involving $A_{vw}$
forms are unaffected and the $r\ge s$ ordering
restriction can be removed. The $t\ge u$ ordering 
restriction can be removed for the same
reason and Lemma \ref{lem:Kron} is established.

\subsubsection{Proof of Lemma \ref{lem:keyIdents}}

We start with a statement of \emph{Woodbury's matrix identity} (Woodbury, 1950).
For invertible matrices $\bS\ (n\times n)$ and  $\bT\ (d\times d)$
and additional matrices $\bU (n\times d)$ and $\bV\ (d\times n)$, 
is
\begin{equation}
(\bS+\bU\bT\bV)^{-1}=\bS^{-1}-\bS^{-1}\bU(\bT^{-1}+\bV\bS^{-1}\bU)^{-1}\bV\bS^{-1}.
\label{eq:WoodburyMatIdent}
\end{equation}
Application of (\ref{eq:WoodburyMatIdent}) with 
$$\bS=\lambda\bI_n,\quad\bT=\bA\quad\mbox{and}\quad\bU=\bV=\bX$$
leads to 
\begin{equation}
(\bX\bA\bX^T+\lambda\bI)^{-1}=(1/\lambda)\bI_n
-(1/\lambda^2)\bX(\bA^{-1}+\bX^T\bX/\lambda)^{-1}\bX^T.
\label{eq:KimberlyAustin}
\end{equation}
Therefore,
{\setlength\arraycolsep{1pt}
\begin{eqnarray*}
\bXTdot(\bX\bA\bXTdot+\lambda\bI)^{-1}\bXdotdot&=&(1/\lambda)\bXTdot\bXdotdot
-(1/\lambda^2)\bXTdot\bX(\bA^{-1}+\bX^T\bX/\lambda)^{-1}\bX^T\bXdotdot\\[1ex]
&=&(1/\lambda)\bXTdot\bXdotdot-(1/\lambda^2)\bXTdot\bX
[(1/\lambda)\bX^T\bX\{\bI_d+\lambda(\bX^T\bX)^{-1}\bA^{-1}\}]^{-1}
\bX^T\bXdotdot\\[1ex]
&=&(1/\lambda)\bXTdot\bXdotdot-(1/\lambda^2)\bXTdot\bX
\{\bI_d+\lambda(\bX^T\bX)^{-1}\bA^{-1}\}^{-1}\lambda(\bX^T\bX)^{-1}\bX^T\bXdotdot\\[1ex]
&=&(1/\lambda)\bXTdot\bXdotdot-(1/\lambda)\bXTdot\bX\{\bI_d+\lambda(\bX^T\bX)^{-1}\bA^{-1}\}^{-1}
(\bX^T\bX)^{-1}\bX^T\bXdotdot.
\end{eqnarray*}
}
Next we apply Woodbury's matrix identity (\ref{eq:WoodburyMatIdent})
to $\{\bI_d+\lambda(\bX^T\bX)^{-1}\bA^{-1}\}^{-1}$ with
$$\bS=\bI_d,\quad\bT=\bA^{-1},\quad\bU=(\bX^T\bX)^{-1}\quad\mbox{and}\quad\bV=\lambda\bI_d$$
to obtain
{\setlength\arraycolsep{1pt}
\begin{eqnarray*}
\{\bI_d+\lambda(\bX^T\bX)^{-1}\bA^{-1}\}^{-1}
&=&\bI_d-(\bX^T\bX)^{-1}\{\bA+\lambda(\bX^T\bX)^{-1}\}^{-1}\lambda.
\end{eqnarray*}
}
Plugging this into the above set of equations we have
{\setlength\arraycolsep{1pt}
\begin{eqnarray*}
\bXTdot(\bX\bA\bX^T+\lambda\bI)^{-1}\bXdotdot
&=&(1/\lambda)\bXTdot\bXdot-(1/\lambda)\bXTdot\bXdot(\bX^T\bXdot)^T
(\bX^T\bX)^{-1}\bX^T\bXdotdot\\[1ex]
&&\quad +\bXTdot\bX(\bX^T\bX)^{-1}\{\bA+\lambda(\bX^T\bX)^{-1}\}^{-1}
(\bX^T\bX)^{-1}\bX^T\bXdotdot\\[1ex]
&=&(1/\lambda)\bXTdot\{\bI_n-\bX(\bX^T\bX)^{-1}\bX^T\}\bXdotdot\\[1ex]
&&\quad +\bXTdot\bX(\bX^T\bX)^{-1}\{\bA+\lambda(\bX^T\bX)^{-1}\}^{-1}
(\bX^T\bX)^{-1}\bX^T\bXdotdot
\end{eqnarray*}
}
and the lemma is proven.

\subsubsection{Proof of Lemma \ref{lem:keySqIdent}}

It follows from (\ref{eq:KimberlyAustin}) that
\begin{equation}
{\setlength\arraycolsep{1pt}
\begin{array}{rcl}
&&\bX^T(\bX\bA\bX^T+\lambda\bI)^{-2}\bX=(1/\lambda^2)\bX^T\bX
-(2/\lambda^3)\bX^T\bX(\bA^{-1}+\bX^T\bX/\lambda)^{-1}\bX^T\bX\\[1ex]
&&\qquad\qquad\qquad\qquad
+(1/\lambda^4)\bX^T\bX(\bA^{-1}+\bX^T\bX/\lambda)^{-1}\bX^T
\bX(\bA^{-1}+\bX^T\bX/\lambda)^{-1}\bX^T\bX
\end{array}
}
\label{eq:BobTheBayesian}
\end{equation}
Steps similar to those given in the proof of Lemma \ref{lem:keyIdents} lead to 
\begin{equation}
(\bA^{-1}+\bX^T\bX/\lambda)^{-1}\bX^T\bX=
\lambda\bI-\lambda^2(\bX^T\bX)^{-1}\{\bA+\lambda(\bX^T\bX)^{-1}\}^{-1}
\label{eq:KuteBelgian}
\end{equation}
Triple substitution of (\ref{eq:KuteBelgian}) into (\ref{eq:BobTheBayesian})
and simplification yields the stated result.

\subsubsection{Proof of Lemma \ref{lem:ApsBinvResIden}}

Lemma \ref{lem:ApsBinvResIden} follows quickly from the following identity:
$$\left[
\begin{array}{cc}
\bA & \bB \\[2ex]
\bB & \bA
\end{array}
\right]
\left[
\begin{array}{c}
\bI_d\\[1ex]
\bI_d
\end{array}
\right]=
\left[
\begin{array}{c}
\bI_d\\[1ex]
\bI_d
\end{array}
\right](\bA+\bB).
$$

\subsubsection{Proof of Lemma \ref{lem:XApart}}

\subsubsubsection{Proof of Lemma \ref{lem:XApart}(a)--(c)}

In the special case of $d=1$, $n_{ii'}=n$ and 
$\bXiid=\bone_n$ for all $1\le i\le m$, $1\le i'\le m'$.
Long-winded, but straightforward, algebraic arguments 
based on the eigenvalue and eigenvector properties described
near the commencement of Section \ref{sec:firstEig}
lead to the exact expression
{\setlength\arraycolsep{0pt}
\begin{eqnarray*}
&&\bone_{m'n}^T\bQ_{mm'}^{ii}\bone_{m'n}=\\[2ex]
&&\qquad
\frac{mm'M'\{(m-1)M'+m' M\}+m'\{(m-1)M' +mM' +m'M\}(\lambda/n)
+m'(\lambda/n)^2}
{(mM'+\lambda/n)\{m'M(m'M+mM')+(mM' +2m'M)(\lambda/n) +(\lambda/n)^2\}}
\end{eqnarray*}
}
for all $1\le i\le m$. This result leads to 
\begin{equation}
\lim_{n\to\infty}\left(\bone_{m'n}^T\bQ_{mm'}^{ii}\bone_{m'n}\right)
=\frac{1}{M}-\frac{M'}{mm' M}\left(\frac{M}{m}+\frac{M'}{m'}\right)^{-1}
\quad\mbox{for all}\quad m,m'\in\naturalNumbers
\label{eq:NashvilleShow}
\end{equation}
which proves Lemma \ref{lem:XApart}(b) in this scalar case. 
Similar calculation lead to
\begin{equation}
\lim_{n\to\infty}\left(\bone_{mm'n}^T\bQ_{mm'}^{-1}\bone_{mm'n}\right)
=\left(\frac{M}{m}+\frac{M'}{m'}\right)^{-1}
\quad\mbox{for all}\quad m,m'\in\naturalNumbers.
\label{eq:OozePaint}
\end{equation}
The result
\begin{equation}
\lim_{n\to\infty}\left(\bone_{m'n}^T\bQ_{mm'}^{i\iunder}\bone_{m'n}\right)
=-\frac{M'}{mm' M}\left(\frac{M}{m}+\frac{M'}{m'}\right)^{-1}
\quad\mbox{for}\quad i\ne\iunder
\label{eq:Dollhouse}
\end{equation}
follows by subtraction and symmetric considerations.

Next, consider the general $d\in\naturalNumbers$ and unrestricted
$n_{ii'}$ setting, but with $m=2$ and $m'=1$. Then
{\setlength\arraycolsep{1pt}
\begin{eqnarray*}
\bQ_{21}=\left[
\begin{array}{cc}
\bX_{11}(\bM+\bM')\bX_{11}^T+\lambda\bI& \bX_{11}\bM'\bX_{21}^T \\[1ex]
\bX_{21}\bM'\bX_{11}^T     & \bX_{21}(\bM+\bM')\bX_{21}^T+\lambda\bI
\end{array}
\right]
\end{eqnarray*}
}
and so, using Corollary 2.1.(c), 
{\setlength\arraycolsep{0pt}
\begin{eqnarray*}
\bX_{11}^T\bQ^{11}_{21}\bX_{11}&=&\bX_{11}^T\big[\bX_{11}(\bM+\bM')\bX_{11}^T\\
&&\qquad\quad-\bX_{11}\bM'\bX_{21}^T\{\bX_{21}(\bM+\bM')\bX_{21}^T+\lambda\bI\}^{-1}
\bX_{21}\bM'\bX_{11}^T+\lambda\bI\big]^{-1}\bX_{11}\\[1ex]
&=&\big[\bM+\bM'-\bM'\{\bM+\bM'+\lambda(\bX_{21}^T\bX_{21})^{-1}\}^{-1}\bM'
+\lambda(\bX_{11}^T\bX_{11})^{-1}\big]^{-1}\\[1ex]
&\convprob&\big\{\bM+\bM'-\bM'(\bM+\bM')^{-1}\bM'\}^{-1}
=\Big[\bM+\bM'\{\bI-(\bM+\bM')^{-1}\bM'\}\Big]^{-1}.
\end{eqnarray*}
}
Noting that $\bI-(\bM+\bM')^{-1}\bM'=(\bM+\bM')^{-1}\bM$ we then have the convergence
in probability limit equalling $\big\{\bM+\bM'(\bM+\bM')^{-1}\bM\big\}^{-1}$.
Application of Woodbury's matrix identity (\ref{eq:WoodburyMatIdent}) with 
$\bS=\bM$, $\bU=\bM'$, $\bT=(\bM+\bM')^{-1}$ and $\bV=\bM$ leads to the limit equalling
$$\bM^{-1}-\bM^{-1}\bM'(\bM+2\bM')^{-1}=\bM^{-1}-\smhalf\bM^{-1}\bM'(\smhalf\bM+\bM')^{-1}$$
which verifies Lemma \ref{lem:XApart}(b) for $m=2$, $m'=1$ and $i=1$. The proof
for $i=2$ is very similar. Then note that, using Corollary 2.1(c),
{\setlength\arraycolsep{0pt}
\begin{eqnarray*}
&&\bX_{11}^T\bQ^{12}_{21}\bX_{21}=\\[1ex]
&&\quad-\bX_{11}^T
\big[\bX_{11}(\bM+\bM')\bX_{11}^T-\bX_{11}\bM'\bX_{21}^T\{\bX_{21}(\bM+\bM')\bX_{21}^T+\lambda\bI\}^{-1}
\bX_{21}\bM'\bX_{11}^T\big]^{-1}\\
&&\qquad\qquad\times\bX_{11}\bM'\bX_{21}^T\{\bX_{21}(\bM+\bM')\bX_{21}^T+\lambda\bI\}^{-1}\bX_{21}\\[1ex]
&&\quad\convprob\,-\big\{\bM+\bM'-\bM'(\bM+\bM')^{-1}\bM'\}^{-1}\bM'(\bM+\bM')^{-1}\\[1ex]
&&\quad=-\smhalf\bM^{-1}\bM'(\smhalf\bM+\bM')^{-1}.
\end{eqnarray*}
}
Hence Lemma \ref{lem:XApart}(c) holds for $m=2$ and $m'=1$. Lemma \ref{lem:XApart}(a) 
for $m=2$ and $m'=1$ follows from summation of the 
Lemma \ref{lem:XApart}(b)--(c) results. This completes verification of 
Lemma \ref{lem:XApart}(a)--(c) for $m=2$ and $m'=1$.

Next we prove Lemma \ref{lem:XApart} for all $m\ge 2$ and $m'=1$
via induction on $m$. 
Let
\begin{equation}
\begin{array}{c}
\bQ_{m+1,1}=\left[
\begin{array}{cc}
\bQ_{m1}     & \TOPRIGHT   \\[1ex]
\TOPRIGHT^T  &  \BOTRIGHT
\end{array}
\right]
\ \mbox{where}\ \BOTRIGHT\equiv\bX_{m+1,1}(\bM+\bM')\bX_{m+1,1}^T+\lambda\bI,\\[4ex] 
\TOPRIGHT\equiv\bX_{1:m,1}\bM'\bX_{m+1,1}^T\quad\mbox{and}\quad  
\bX_{1:m,1}\equiv {\displaystyle\stack{1\le i\le m}}(\bX_{i1}).
\end{array}
\label{eq:QmpOneDefn}
\end{equation}
We then have, with use of Corollary 2.1.(c),
{\setlength\arraycolsep{1pt}
\begin{eqnarray*}
\bX_{m+1,1}^T\bQ_{m+1,1}^{m+1,m+1}\bX_{m+1,1}
&=&\bX_{m+1,1}^T\Big(\BOTRIGHT-\TOPRIGHT^T\bQ_{m1}^{-1}\TOPRIGHT\Big)^{-1}\bX_{m+1,1}\\[1ex]
&=&\bX_{m+1,1}^T\Big\{\bX_{m+1,1}(\bM+\bM')\bX_{m+1,1}^T\\
&&\qquad\qquad-\bX_{m+1,1}\bM' 
\bX_{1:m,1}^T\bQ_{m1}^{-1}\bX_{1:m,1}\bM'\bX_{m+1,1}^T+\lambda\bI      
\Big\}^{-1}\bX_{m+1,1}\\[1ex]
&=&\Big\{\bM+\bM'
-\bM'\bX_{1:m,1}^T\bQ_{m1}^{-1}\bX_{1:m,1}\bM'+\lambda(\bX_{m+1,1}^T\bX_{m+1,1})^{-1}\Big\}^{-1}\\[1ex]
&\convprob&\left\{\bM+\bM'
-\bM'\Big(\frac{1}{m}\bM+\bM'\Big)^{-1}\bM'\right\}^{-1}\\[2ex]
&=&\bM^{-1}-\frac{1}{m+1}\bM^{-1}\bM'\Big(\frac{1}{m+1}\bM+\bM'\Big)^{-1}.
\end{eqnarray*}
}
Analogous arguments for other partitions of $\bQ_{m+1,1}$ lead to the same
convergence in probability limit for $\bX_{i,1}^T\bQ_{m+1,1}^{i,i}\bX_{i,1}$
for each $1\le i\le m+1$. Therefore, by induction, Lemma \ref{lem:XApart}(b) 
holds for all $m\ge 2$ and $m'=1$.

Next, let
$\bQ_{m+1,1}^{1:m,m+1}\equiv{\displaystyle\stack{1\le i\le m}}\big(\bQ_{m+1,1}^{i,m+1}\big)$
and note that
{\setlength\arraycolsep{1pt}
\begin{eqnarray*}
\bX_{1:m,1}^T\bQ_{m+1,1}^{1:m,m+1}\bX_{m+1,1}&=&-\bX_{1:m,1}^T\bQ_{m1}^{-1}\TOPRIGHT
\Big(\BOTRIGHT-\TOPRIGHT^T\bQ_{m1}^{-1}\TOPRIGHT\Big)^{-1}\bX_{m+1,1}\\[1ex]
&=&-\bX_{1:m,1}^T\bQ_{m1}^{-1}\bX_{1:m,1}\bM'\bX_{m+1,1}^T\big\{\bX_{m+1,1}(\bM+\bM')\bX_{m+1,1}^T\\
&&\qquad
-\bX_{m+1,1}\bM'\bX_{1:m,1}^T\bQ_{m1}^{-1}\bX_{1:m,1}\bM'\bX_{m+1,1}^T+\lambda\bI\big\}^{-1}\bX_{m+1,1}\\[1ex]
&=&-\bX_{1:m,1}^T\bQ_{m1}^{-1}\bX_{1:m,1}\bM'\big\{\bM+\bM'
-\bM'\bX_{1:m,1}^T\bQ_{m1}^{-1}\bX_{1:m,1}\bM'\\
&&\qquad\qquad\qquad\qquad\qquad\quad+\lambda(\bX_{m+1,1}^T\bX_{m+1,1})^{-1}
\big\}^{-1}\\[1ex]
&\convprob&-\Big(\frac{1}{m}\bM+\bM'\Big)^{-1}\bM'\Big\{\bM+\bM'
-\bM'\Big(\frac{1}{m}\bM+\bM'\Big)^{-1}\bM'\Big\}^{-1}\\[1ex]
&=&-\Big(\frac{1}{m}\bM+\bM'\Big)^{-1}\bM'
\Big\{\bM^{-1}-\frac{1}{m+1}\bM^{-1}\bM'\Big(\frac{1}{m+1}\bM+\bM'\Big)^{-1}\Big\}\\[1ex]
&=&-\frac{m}{m+1}\bM^{-1}\bM'\Big(\frac{1}{m+1}\bM+\bM'\Big)^{-1}.
\end{eqnarray*}
}
However,
\begin{equation}
\bX_{1:m,1}^T\bQ_{m+1,1}^{1:m,m+1}\bX_{m+1,1}=\sum_{i=1}^m\bX_{i1}^T\bQ_{m+1,1}^{i,m+1}\bX_{m+1,1}
\label{eq:SaintJohn}
\end{equation}
and each term in the summation on the right-hand side of (\ref{eq:SaintJohn}) has the same distribution and,
therefore, the same convergence in probability limit. Hence,
$$\bX_{i1}^T\bQ_{m+1,1}^{i,m+1}\bX_{m+1,1}\convprob -\frac{1}{m+1}\bM^{-1}\bM'\Big(\frac{1}{m+1}\bM+\bM'\Big)^{-1},
\quad 1\le i\le m.
$$
Analogous arguments for other partitions of $\bQ_{m+1}$ lead to 
$$\bX_{i1}^T\bQ_{m+1,1}^{i,\iunder}\bX_{\iunder 1}\convprob 
-\frac{1}{m+1}\bM^{-1}\bM'\Big(\frac{1}{m+1}\bM+\bM'\Big)^{-1},\quad 1\le i\ne\iunder\le m+1,$$
and by induction, Lemma \ref{lem:XApart}(b) and (c) hold for $m\ge 2$ and $m'=1$. 

To establish Lemma \ref{lem:XApart}(a) for $m\ge 2$ and $m'=1$ we sum the results 
just derived for Lemma \ref{lem:XApart}(b) and (c):
{\setlength\arraycolsep{1pt}
\begin{eqnarray*}
\Big\{\stack{1 \le i \le m+1}(\bX_{i1})\Big\}^T
\bQ_{m+1,1}^{-1}\stack{1 \le i \le m+1}(\bX_{i1})
&=&\sum_{i=1}^{m+1}\sum_{\iunder=1}^{m+1}\bX_{i,1}^T\bQ_{m+1,1}^{i,\iunder}\bX_{\iunder,1}\\[1ex]
&\convprob&(m+1)\left\{\bM^{-1}-\frac{1}{m+1}\bM^{-1}\bM'\Big(\frac{1}{m+1}\bM+\bM'\Big)^{-1}\right\}\\[1ex]
&&\quad+m(m+1)\left\{-\frac{1}{m+1}\bM^{-1}\bM'\Big(\frac{1}{m+1}\bM+\bM'\Big)^{-1}\right\}\\[1ex]
&=&\left(\frac{\bM}{m+1}+\bM'\right)^{-1}.
\end{eqnarray*}
}
Induction then leads to Lemma \ref{lem:XApart}(a) holding for $m\ge 2$ and $m'=1$.
This completes verification of Lemma \ref{lem:XApart} for $m\ge2$ and $m'=1$.

For the $m=1$ and $m'=2$ case the matrix of interest is
$$
\bQ_{12}^{11}=\bQ_{12}^{-1}$$
where
{\setlength\arraycolsep{1pt}
\begin{eqnarray*}
\bQ_{12}&=&\left[
\begin{array}{cc}
\bX_{11}(\bM+\bM')\bX_{11}^T& \bX_{11}\bM\bX_{12}^T \\[1ex]
\bX_{12}\bM\bX_{11}^T     & \bX_{12}(\bM+\bM')\bX_{12}^T
\end{array}
\right]+\lambda\bI\\[1ex]
&=&
\left[
\begin{array}{cc}
\bX_{11} & \bO \\[1ex]
\bO      & \bX_{12}
\end{array}
\right]\left[
\begin{array}{cc}
\bM+\bM' &  \bM \\[2ex]
\bM      & \bM+\bM'
\end{array}
\right]
\left[
\begin{array}{cc}
\bX_{11} & \bO \\[1ex]
\bO      & \bX_{12}
\end{array}
\right]^T+\lambda\bI.
\end{eqnarray*}
}
Noting that
$$\bXuptrione\equiv\left[
\begin{array}{c}
\bX_{11}\\[1ex]
\bX_{12}
\end{array}
\right]
=
\left[
\begin{array}{cc}
\bX_{11} & \bO \\[1ex]
\bO      & \bX_{12}
\end{array}
\right]
\left[
\begin{array}{c}
\bI_d\\[1ex]
\bI_d
\end{array}
\right]
$$
we have
{\setlength\arraycolsep{1pt}
\begin{eqnarray*}
\bXuptrione^T\bQ_{12}^{11}\bXuptrione&=&
\left[
\begin{array}{c}
\bI_d\\[1ex]
\bI_d
\end{array}
\right]^T
\left[
\begin{array}{cc}
\bX_{11} & \bO \\[1ex]
\bO      & \bX_{12}
\end{array}
\right]^T
\left\{
\left[
\begin{array}{cc}
\bX_{11} & \bO \\[1ex]
\bO      & \bX_{12}
\end{array}
\right]
\left[
\begin{array}{cc}
\bM+\bM' &  \bM \\[2ex]
\bM      & \bM+\bM'
\end{array}
\right]
\left[
\begin{array}{cc}
\bX_{11} & \bO \\[1ex]
\bO      & \bX_{12}
\end{array}
\right]^T
+\lambda\bI
\right\}^{-1}\\
&&\qquad\qquad\qquad\times\left[
\begin{array}{cc}
\bX_{11} & \bO \\[1ex]
\bO      & \bX_{12}
\end{array}
\right]
\left[
\begin{array}{c}
\bI_d\\[1ex]
\bI_d
\end{array}
\right]
\\[2ex]
&=&\left[
\begin{array}{c}
\bI_d\\[1ex]
\bI_d
\end{array}
\right]^T
\left\{
\left[
\begin{array}{cc}
\bM+\bM' &  \bM \\[2ex]
\bM      & \bM+\bM'
\end{array}
\right]
+\lambda\left[\begin{array}{cc}   
(\bX_{11}^T\bX_{11})^{-1} & \bO\\[1ex]
    \bO & (\bX_{12}^T\bX_{12})^{-1}
\end{array}
\right]
\right\}^{-1}
\left[
\begin{array}{c}
\bI_d\\[1ex]
\bI_d
\end{array}
\right]
\\[1ex]
&\convprob&
\left[
\begin{array}{c}
\bI_d\\[2ex]
\bI_d
\end{array}
\right]^T
\left[
\begin{array}{cc}
\bM+\bM' &  \bM \\[2ex]
\bM      & \bM+\bM'
\end{array}
\right]^{-1}
\left[
\begin{array}{c}
\bI_d\\[2ex]
\bI_d
\end{array}
\right]=(\bM+\smhalf\bM')^{-1}
\end{eqnarray*}
}
where the last equality is due to Lemma \ref{lem:ApsBinvResIden}.

For the $m=2$ and $m'=2$ case the matrix of interest is
{\setlength\arraycolsep{1pt}
\begin{eqnarray*}
\bQ_{22}^{11}&=&\mbox{the top left $(n_{11}+n_{12})\times(n_{11}+n_{12})$ block of}\ 
\left[
\begin{array}{ccc}
\bQ_{12}    & \bR_{12} \\[1ex]
\bR_{12}^T  & \bQtilde_{12}
\end{array}
\right]^{-1}
\end{eqnarray*}
}
where 
$$
\bR_{12}=\left[
\begin{array}{cc}
\bX_{11}\bM'\bX_{12}^T & \bO \\
\bO                   & \bX_{21}\bM'\bX_{22}^T 
\end{array}
\right]\\[1ex]
=\left[
\begin{array}{cc}
\bX_{11} & \bO \\[1ex]
\bO      & \bX_{12}
\end{array}
\right]
\left[
\begin{array}{cc}
\bM' & \bO \\[1ex]
\bO      & \bM'
\end{array}
\right]
\left[
\begin{array}{cc}
\bX_{21} & \bO \\[1ex]
\bO      & \bX_{22}
\end{array}
\right]^T
$$
and
{\setlength\arraycolsep{1pt}
\begin{eqnarray*}
\bQtilde_{12}&=&\left[
\begin{array}{cc}
\bX_{21}(\bM+\bM')\bX_{21}^T+\lambda\bI& \bX_{21}\bM\bX_{22}^T \\[1ex]
\bX_{22}\bM\bX_{21}^T     & \bX_{22}(\bM+\bM')\bX_{22}^T+\lambda\bI
\end{array}
\right]\\[2ex]
&=&\left[
\begin{array}{cc}
\bX_{21} & \bO \\[1ex]
\bO      & \bX_{22}
\end{array}
\right]\left[
\begin{array}{cc}
\bM+\bM' &  \bM \\[2ex]
\bM      & \bM+\bM'
\end{array}
\right]
\left[
\begin{array}{cc}
\bX_{21} & \bO \\[1ex]
\bO      & \bX_{22}
\end{array}
\right]^T+\lambda\bI.
\end{eqnarray*}
}
Therefore, 
{\setlength\arraycolsep{1pt}
\begin{eqnarray*}
\bXuptrione^T\bQ_{22}^{11}\bXuptrione&=&
\bXuptrione^T
\left(\bQ_{12}-\bR_{12}\bQtilde_{12}^{-1}\bR_{12}^T\right)^{-1}
\bXuptrione
\\[1ex]
&=&
\left[
\begin{array}{c}
\bI_d\\[1ex]
\bI_d
\end{array}
\right]^T
\left[
\begin{array}{cc}
\bX_{11} & \bO \\[1ex]
\bO      & \bX_{12}
\end{array}
\right]^T
\Bigg\{\left[
\begin{array}{cc}
\bX_{11} & \bO \\[1ex]
\bO      & \bX_{12}
\end{array}
\right]\left[
\begin{array}{cc}
\bM+\bM' &  \bM \\[2ex]
\bM      & \bM+\bM'
\end{array}
\right]
\left[
\begin{array}{cc}
\bX_{11} & \bO \\[1ex]
\bO      & \bX_{12}
\end{array}
\right]^T  \\
&&\ -\left[
\begin{array}{cc}
\bX_{11} & \bO \\[1ex]
\bO      & \bX_{12}
\end{array}
\right]
\left[
\begin{array}{cc}
\bM' & \bO \\[1ex]
\bO      & \bM'
\end{array}
\right]
\left[
\begin{array}{cc}
\bX_{21} & \bO \\[1ex]
\bO      & \bX_{22}
\end{array}
\right]^T\bQtilde^{-1}_{12}
\left[
\begin{array}{cc}
\bX_{21} & \bO \\[1ex]
\bO      & \bX_{22}
\end{array}
\right]
\left[
\begin{array}{cc}
\bM' & \bO \\[1ex]
\bO      & \bM'
\end{array}
\right]
\left[
\begin{array}{cc}
\bX_{11} & \bO \\[1ex]
\bO      & \bX_{12}
\end{array}
\right]^T
\Bigg\}^{-1}\\
&&\quad\times
\left[
\begin{array}{cc}
\bX_{11} & \bO \\[1ex]
\bO      & \bX_{12}
\end{array}
\right]
\left[
\begin{array}{c}
\bI_d\\[1ex]
\bI_d
\end{array}
\right]
\\[2ex]
&=&
\left[
\begin{array}{c}
\bI_d\\[1ex]
\bI_d
\end{array}
\right]^T
\Bigg\{\left[
\begin{array}{cc}
\bM+\bM' &  \bM \\[2ex]
\bM      & \bM+\bM'
\end{array}
\right]-
\left[
\begin{array}{cc}
\bM' & \bO \\[1ex]
\bO      & \bM'
\end{array}
\right]
\bPsi
\left[
\begin{array}{cc}
\bM' & \bO \\[1ex]
\bO  & \bM'
\end{array}
\right]\\[1ex]
&&\qquad\qquad\qquad\qquad\qquad\qquad
+\lambda\left[\begin{array}{cc}   
(\bX_{11}^T\bX_{11})^{-1} & \bO\\[1ex]
    \bO & (\bX_{12}^T\bX_{12})^{-1}
\end{array}
\right]
\Bigg\}^{-1}
\left[
\begin{array}{c}
\bI_d\\[1ex]
\bI_d
\end{array}
\right]
\\
\end{eqnarray*}
}
where
$$\bPsi\equiv
\left[
\begin{array}{cc}
\bX_{21} & \bO \\[1ex]
\bO      & \bX_{22}
\end{array}
\right]^T\bQtilde^{-1}_{12}
\left[
\begin{array}{cc}
\bX_{21} & \bO \\[1ex]
\bO      & \bX_{22}
\end{array}
\right].
$$
Now note that 
$$
\left[
\begin{array}{cc}
\bM' & \bO \\[1ex]
\bO      & \bM'
\end{array}
\right]
\bPsi
\left[
\begin{array}{cc}
\bM' & \bO \\[1ex]
\bO  & \bM'
\end{array}
\right]
=
\left[
\begin{array}{cc}
\bM'\bX_{21}^T\bQtilde_{12}^{[1,1]}\bX_{21}\bM' &  
\bM'\bX_{21}^T\bQtilde_{12}^{[1,2]}\bX_{22}\bM' \\[1ex]
\bM'\bX_{22}^T\bQtilde_{12}^{[2,1]}\bX_{21}\bM' &  
\bM'\bX_{22}^T\bQtilde_{22}^{[2,2]}\bX_{22}\bM' 
\end{array}
\right]
$$
where
$$\left[
\begin{array}{cc}
\bQtilde_{12}^{[1,1]} & \bQtilde_{12}^{[1,2]} \\[1ex]
\bQtilde_{12}^{[2,1]} & \bQtilde_{12}^{[2,2]}
\end{array}
\right]
$$
is the partition of $\bQtilde_{12}$ such that the sub-blocks have
dimensions:
$$\bQtilde_{12}^{[1,1]}\ \mbox{is}\ n_{21}\times n_{21},
\quad \bQtilde_{12}^{[1,2]}\ \mbox{is}\ n_{21}\times n_{22},
\quad \bQtilde_{12}^{[2,1]}\ \mbox{is}\ n_{22}\times n_{21}
\quad\mbox{and}\quad
\bQtilde_{12}^{[2,2]}\ \mbox{is}\ n_{22}\times n_{22}.
$$
We then have
{\setlength\arraycolsep{1pt}
\begin{eqnarray*}
\bXuptrione^T\bQ_{22}^{11}\bXuptrione&=&
\left[
\begin{array}{c}
\bI_d\\
\bI_d
\end{array}
\right]^T
\left[
\begin{array}{cc}
\bM+\bM'-\bM'\bX_{21}^T\bQtilde_{12}^{[1,1]}\bX_{21}\bM' &  
\bM-\bM'\bX_{21}^T\bQtilde_{12}^{[1,2]}\bX_{22}\bM' \\[1ex]
+\lambda(\bX_{11}^T\bX_{11})^{-1} &                 \\[2ex]
\bM-\bM'\bX_{22}^T\bQtilde_{12}^{[2,1]}\bX_{21}\bM' &  
\bM+\bM'-\bM'\bX_{22}^T\bQtilde_{12}^{[2,2]}\bX_{22}\bM'\\[1ex]
 & +\lambda(\bX_{12}^T\bX_{12})^{-1}
\end{array}
\right]^{-1}
\left[
\begin{array}{c}
\bI_d\\
\bI_d
\end{array}
\right]\\[1ex]
&=&
\left[
\begin{array}{c}
\bI_d\\
\bI_d
\end{array}
\right]^T
\left[
\begin{array}{cc}
\bAtilde &\quad\ \bBtilde \\[1ex]
\bBtilde &\quad\ \bAtilde
\end{array}
\right]^{-1}
\left[
\begin{array}{c}
\bI_d\\
\bI_d
\end{array}
\right]
\end{eqnarray*}
}
where
$$\bAtilde\equiv \bM+\bM'-\smhalf\bM'\left\{\bX_{21}^T\bQtilde_{12}^{[1,1]}\bX_{21}
+\bX_{22}^T\bQtilde_{12}^{[2,2]}\bX_{22}\right\}\bM'\{1+o_P(1)\}
+\lambda(\bX_{11}^T\bX_{11})^{-1}$$
and 
$$\bBtilde\equiv \bM -\smhalf\bM'\left\{\bX_{21}^T\bQtilde_{12}^{[1,2]}\bX_{22}
+\bX_{22}^T(\bQtilde_{12}^{[1,2]})^T\bX_{21}
\right\}\bM'\{1+o_P(1)\}$$
with the $\{1+o_P(1)\}$ factors being justified due to each of
$\bX_{11}$, $\bX_{11}$, $\bX_{21}$ and $\bX_{22}$ containing
random samples from the same distribution.
Application of Lemma \ref{lem:ApsBinvResIden} 
leads to, with $\bXtilde\equiv\big[\bX_{21}^T\ \bX_{22}^T\big]^T$ 
being the $\bQtilde_{12}$ 
version of the $\bX$ matrix from Lemma \ref{lem:XApart}(b) but
for $\bQtilde_{12}$ rather than $\bQ_{12}$, the result
{\setlength\arraycolsep{1pt}
\begin{eqnarray*}
\bXuptrione^T\bQ_{22}^{11}\bXuptrione
&=&2\Big[\bM+\bM'+\bM-\smhalf
\bM'\big(\bXtilde^T\bQtilde_{12}^{-1}\bXtilde\big)\bM'\{1+o_P(1)\}
+\lambda(\bX_{11}^T\bX_{11})^{-1}\{1+o_P(1)\}\Big]^{-1}
\\[1ex]
&\convprob&
2\{2\bM+\bM'-\smhalf\bM'(\bM+\smhalf\bM')^{-1}\bM'\}^{-1}
=\bM^{-1}-\quarter\bM^{-1}\bM'\big(\smhalf\bM+\smhalf\bM'\big)^{-1}.
\end{eqnarray*}
}
which verifies Lemma \ref{lem:XApart}(b) for the $(m,m')=(2,2)$ case.
Induction on $m$ can be used to show that Lemma \ref{lem:XApart}(b) 
holds for general $m\in\naturalNumbers$ and $m'=2$. 

It is apparent from these derivations in the $m'\in\{1,2\}$ cases 
that the behaviors of the summations that lead to
the limits given by (\ref{eq:NashvilleShow})--(\ref{eq:Dollhouse}) in the
$d=1$ and balanced cell counts situation also lead to the
analogous matrix forms for general $m'\in\naturalNumbers$.

\subsubsubsection{Proof of Lemma \ref{lem:XApart}(d)}

In the special case of $d=1$, $n_{ii'}=n$ and 
$\bXiid=\bone_n$ for all $1\le i\le m$, $1\le i'\le m'$.
The eigenvalue and eigenvector properties described
near the start of Section \ref{sec:firstEig} are such that
relatively straightforward manipulations produce the exact expression
{\setlength\arraycolsep{1pt}
\begin{eqnarray*}
\frac{1}{mm'n}\,\tr(\bQmmd^{-2})
&=&\frac{1}{\lambda^2}-\frac{M/(m'n)}{\lambda^2\{M+\lambda/(m'n)\}}
-\frac{M'/(mn)}{\lambda^2\{M'+\lambda/(mn)\}}\\[1ex]
&&\quad+\frac{MM'/(mm'n)}{\lambda^2\{M+\lambda/(m'n)\}\{M(m'/m)+M'+\lambda/(mn)\}}\\[1ex]
&&\quad 
+\frac{MM'/(mm'n)}{\lambda^2 \{M'+\lambda/(mn)\}\{M+M'(m/m')+\lambda/(m'n)\}}\\[1ex]
&&\quad-\frac{M'/\{(mn)^2\}}{\lambda\{M'+\lambda/(mn)\}^2}
-\frac{M/\{(m'n)^2\}}{\lambda\{M+\lambda/(m'n)\}^2}
\\[1ex]
&&\quad +\frac{MM'/\{m'(mn)^2\}}{\lambda\{M+M'(m/m')+\lambda/(m'n)\}\{M'+\lambda/(mn)\}^2}
\\[1ex]
&&\quad+\frac{MM'/\{m'(mn)^2\}}{\lambda\{M(m'/m)+M'+\lambda/(mn)\}^2\{M+\lambda/(m'n)\}}\\[2ex]
&&\quad +\frac{MM'/\{m(m'n)^2\}}{\lambda\{M(m'/m)+M'+\lambda/(mn)\}
\{M+\lambda/(m'n)\}^2}
\\[1ex]
&&\quad+\frac{MM'/\{m(m'n)^2\}}{\lambda\{M+M'(m/m')+\lambda/(m'n)\}^2
\{M'+\lambda/(mn)\}}.
\end{eqnarray*}
}
Hence, under (A5),
\begin{equation}
\frac{1}{mm'n}\,\tr(\bQmmd^{-2})
=\left(\sumim\sum_{i'=1}^{m'}n_{ii'}\right)^{-1}\tr(\bQmmd^{-2})
\to\frac{1}{\lambda^2}
\label{eq:FishINFsigsigScal}
\end{equation}
for all $1\le i\le m$, $1\le i'\le m'$.

Next consider the case of $\dA\in\naturalNumbers$ and $m=m'=1$. Then
$$\bQ_{11}^2=\lambda^2\bI_{n_{11}}+\bX_{11}\bOmega_1\bX_{11}^T
\quad\mbox{where}\quad 
\bOmega_1\equiv(\bM+\bM')\bX_{11}^T\bX_{11}(\bM+\bM')
+2\lambda(\bM+\bM').
$$
Application of Woodbury's matrix identity (\ref{eq:WoodburyMatIdent}) with 
$$\bS=\lambda^2\bI_{n_{11}},\quad \bU\equiv\bX_{11}\bOmega_1,\quad\bT=\bI_{\dA}
\quad\mbox{and}\quad\bV\equiv \bX_{11}^T
$$
then gives
$$\bQ_{11}^{-2}=\lambda^{-2}\bI_{n_{11}}
-\lambda^{-4}\bX_{11}\bOmega_1(\bI+\lambda^{-2}\bX_{11}^T\bX_{11}\bOmega_1)^{-1}\bX^T_{11}$$
and so
{\setlength\arraycolsep{1pt}
\begin{eqnarray*}
\frac{1}{n_{11}}\tr(\bQ_{11}^{-2})&=&\frac{1}{\lambda^2}
-\frac{1}{n_{11}\lambda^4}\tr\Big(
(\bI+\lambda^{-2}\bX_{11}^T\bX_{11}\bOmega_1)^{-1}\bX^T_{11}\bX_{11}\bOmega_1\Big)\\[1ex]
&=&\frac{1}{\lambda^2}
-\frac{1}{n_{11}\lambda^2}\tr\Big(
\big\{\bOmega_1+\lambda^2(\bX_{11}^T\bX_{11})^{-1}\big\}^{-1}\bOmega_1\Big)
\convprob\frac{1}{\lambda^2}.
\end{eqnarray*}
}

For the $\dA\in\naturalNumbers$, $m\in\naturalNumbers$ and $m'=1$ extension we 
note, as given earlier in (\ref{eq:QmpOneDefn}), that
$$
\begin{array}{c}
\bQ_{m+1,1}=\left[
\begin{array}{cc}
\bQ_{m1}     & \TOPRIGHT   \\[1ex]
\TOPRIGHT^T  &  \BOTRIGHT
\end{array}
\right]
\ \mbox{where}\ \BOTRIGHT\equiv\bX_{m+1,1}(\bM+\bM')\bX_{m+1,1}^T+\lambda\bI,\\[4ex] 
\TOPRIGHT\equiv\bX_{1:m,1}\bM'\bX_{m+1,1}^T\quad\mbox{and}\quad  
\bX_{1:m,1}\equiv {\displaystyle\stack{1\le i\le m}}(\bX_{i1})
\end{array}
$$
which gives
$$
\bQ_{m+1,1}^2
=\left[
\begin{array}{cc}
\bQ^2_{m1}+\TOPRIGHT\TOPRIGHT^T &\quad \bQ_{m1}\TOPRIGHT+\TOPRIGHT\BOTRIGHT \\[1ex]
(\bQ_{m1}\TOPRIGHT+\TOPRIGHT\BOTRIGHT)^T & \BOTRIGHT^2+\TOPRIGHT^T\TOPRIGHT
\end{array}
\right].
$$
Then the lower right $n_{m+1,1}\times n_{m+1,1}$ block of 
$\bQ^{-2}_{m1}$ equals
{\setlength\arraycolsep{1pt}
\begin{eqnarray*}
&&\big\{\BOTRIGHT^2+\TOPRIGHT^T\TOPRIGHT-
(\bQ_{m1}\TOPRIGHT+\TOPRIGHT\BOTRIGHT)^T 
(\bQ^2_{m1}+\TOPRIGHT\TOPRIGHT^T)^{-1}
(\bQ_{m1}\TOPRIGHT+\TOPRIGHT\BOTRIGHT)\big\}^{-1}\\[1ex]
&&\qquad\qquad
=\big(\lambda^2\bI_{n_{m+1,1}}+\bX_{m+1,1}^T\bOmega_2\bX_{m+1,1}\Big)^{-1}
\end{eqnarray*}
}
where
{\setlength\arraycolsep{1pt}
\begin{eqnarray*}
\bOmega_2&\equiv&2\lambda(\bM+\bM')+(\bM+\bM')\bX_{m+1,1}^T\bX_{m+1,1}(\bM+\bM')
+\bM'\bX_{1:m,1}^T\bX_{1:m,1}\bM'\\[1ex]
&&\qquad-\bOmega_3^T(\bQ^2_{m1}+\TOPRIGHT\TOPRIGHT^T)^{-1}\bOmega_3
\end{eqnarray*}
}
with
$$\bOmega_3\equiv(\bQ_{m1}+\lambda\bI)\bX_{1:m,1}\bM'
+\bX_{1:m,1}\bM' \bX_{m+1,1}^T\bX_{m+1,1}(\bM+\bM').
$$
Another application of (\ref{eq:WoodburyMatIdent}) with 
$$\bS=\lambda^2\bI_{n_{m+1,1}},\quad \bU\equiv\bX_{m+1,1}\bOmega_2,\quad\bT=\bI_{\dA}
\quad\mbox{and}\quad\bV\equiv\bX_{m+1,1}^T
$$
then gives the lower right $n_{m+1,1}\times n_{m+1,1}$ block of 
$\bQ^{-2}_{m1}$ equalling
$$\lambda^{-2}\bI_{n_{m+1,1}}
-\lambda^{-4}\bX_{m+1,1}\bOmega_2(\bI+\lambda^{-2}\bX_{m+1,1}^T\bX_{m+1,1}\bOmega_2)^{-1}\bX^T_{m+1,1}$$
and so
{\setlength\arraycolsep{1pt}
\begin{eqnarray*}
&&\frac{1}{n_{m+1,1}}\tr\Big(\mbox{lower right $n_{m+1,1}\times n_{m+1,1}$ block of 
$\bQ^{-2}_{m1}$}\Big)\\[1ex]
&&\qquad=\frac{1}{\lambda^2}
-\frac{1}{n_{m+1,1}\lambda^4}\tr\Big(
(\bI+\lambda^{-2}\bX_{m+1,1}^T\bX_{m+1,1}\bOmega_2)^{-1}\bX^T_{m+1,1}
\bX_{m+1,1}\bOmega_2\Big)\\[1ex]
&&\qquad=\frac{1}{\lambda^2}
-\frac{1}{n_{m+1,1}\lambda^2}\tr\Big(
\big\{\bOmega_2+\lambda^2(\bX_{m+1,1}^T\bX_{m+1,1})^{-1}\big\}^{-1}\bOmega_2\Big)
\convprob\frac{1}{\lambda^2}.
\end{eqnarray*}
}
By induction on $m$ we then have, under (A5),
$$\left(\sumim n_{i1}\right)^{-1}\tr(\bQ_{m1}^{-2})
\convprob\frac{1}{\lambda^2}
\quad\mbox{for all}\ m\in\naturalNumbers.
$$

For higher $m'$, similar arguments can be used to show that
the summations in $\tr(\bQ_{mm'}^{-2})$ lead to convergents
that are analogous to those in the $\dA=1$, $n_{ii'}=n$
and $\bXAiid=\bone_n$ case and Lemma \ref{lem:XApart}(d) holds.

\subsubsection{Proof of Lemma \ref{lem:XBpart}}

First we prove Lemma \ref{lem:XBpart} 
for $m=m'=1$, for which the $\bQ$ 
matrix reduces to 
$$\bQ_{11}=\bX_{11}(\bM+\bM')\bX_{11}^T+\lambda\bI.$$
Then, from Lemma \ref{lem:keyIdents},
{\setlength\arraycolsep{1pt}
\begin{eqnarray*}
\bXdiamd_{11}^T\bQ_{11}^{-1}\bXdiamd_{11}
&=&\bXdiamd_{11}^T\big\{\bX_{11}(\bM+\bM')\bX_{11}^T+\lambda\bI\big\}^{-1}\bXdiamd_{11}\\[1ex]
&=&(1/\lambda)\bXdiamd_{11}^T
\{\bI-\bX_{11}(\bX_{11}^T\bX_{11})^{-1}\bX_{11}^T\}\bXdiamd_{11}\\[1ex]
&&\qquad+\bXdiamd^T_{11}\bX_{11}(\bX_{11}^T\bX_{11})^{-1}\{\bM+\bM'
+\lambda(\bX_{11}^T\bX_{11})^{-1}\}^{-1}
(\bX_{11}^T\bX_{11})^{-1}\bX_{11}^T\bXdiamd_{11}.
\end{eqnarray*}
}
Hence 
{\setlength\arraycolsep{1pt}
\begin{eqnarray*}
\oonoo\bXdiamd_{11}^T\bQ_{11}^{-1}\bXdiamd_{11}
&=&(1/\lambda)\left\{\left(\oonoo\bXdiamd_{11}^T\bXdiamd_{11}\right)    
-\left(\oonoo\bXdiamd_{11}^T\bX_{11}\right)
\left(\oonoo\bX_{11}^T\bX_{11}\right)^{-1}
\left(\oonoo\bX_{11}^T\bXdiamd_{11}\right)\right\}\\[1ex]
&&\qquad+\oonoo\left(\oonoo\bXdiamd^T_{11}\bX_{11}\right)\left(\oonoo\bX_{11}^T\bX_{11}\right)^{-1}
\left\{\bM+\bM'
+\lambda(\bX_{11}^T\bX_{11})^{-1}\right\}^{-1}\\
&&\qquad\qquad\qquad\times\left(\oonoo\bX_{11}^T\bX_{11}\right)^{-1}\left(\oonoo\bX_{11}^T\bXdiamd_{11}\right)\\[1ex]
&\convprob&(1/\lambda)\Big[E\Big(\bXdiamdrv^{\otimes2}\Big)-E\Big(\bXdiamdrv\bXrv^T\Big)
\Big\{E\Big(\bXrv^{\otimes2}\Big)\Big\}^{-1}E\Big(\bXdiamdrv\bXdiamdrv^T\Big)\Big]\\[1ex]
&=&(1/\lambda)\left[\mbox{lower right $\ddiamd\times\ddiamd$ block of}\
\big\{E\big([\bXrv\  \bXTdiamdrv]^{\otimes2}\big)\big\}^{-1}\right]^{-1}.
\end{eqnarray*}
}
Thus, Lemma \ref{lem:XBpart} (a) holds for $m=m'=1$. 

To establish Lemma \ref{lem:XBpart}(b) 
for $m=m'=1$ we apply Corollary 2.1(a) to obtain
{\setlength\arraycolsep{1pt}
\begin{eqnarray*}
\bX_{11}^T\bQ_{11}^{-1}\bXdiamd_{11}&=&
\bX_{11}^T\big\{\bX_{11}(\bM+\bM')\bX_{11}^T+\lambda\bI\big\}^{-1}\bXdiamd_{11}\\[1ex]
&=&\{\bM+\bM'+\lambda(\bX_{11}^T\bX_{11})^{-1}\}^{-1}\left(\oonoo\bX_{11}^T\bX_{11}\right)^{-1}
\oonoo\bX_{11}^T\bXdiamd_{11}\\[1ex]
&\convprob&(\bM+\bM')^{-1}\big\{E(\bXrv^{\otimes2})\big\}^{-1}E(\bXrv\bXdiamdrv^T).
\end{eqnarray*}
}
Therefore, Lemma \ref{lem:XBpart} is proven for $m=m'=1$. 

Next we prove that the lemma holds for all $m\ge1$ and $m'=1$ via induction on $m$.
Let $\bQ_{m1}$ denote the $m'=1$ version of (\ref{eq:lemCC}) and consider the partition
of $\bQ_{m+1,1}$ given by (\ref{eq:QmpOneDefn}). Also let 
$$\bXdiamd_{1:m,1}\equiv\stack{1\le i\le m}(\bXdiamd_{i1})\quad\mbox{and}\quad 
\bXdiamd_{1:m+1,1}\equiv\stack{1\le i\le m+1}(\bXdiamd_{i1})
=\left[
\begin{array}{c}
\bXdiamd_{1:m,1}\\[1ex]
\bXdiamd_{m+1,1}
\end{array}
\right].
$$
Then
{\setlength\arraycolsep{1pt}
\begin{eqnarray*}
&&\bXdiamd_{1:m+1,1}^T\bQ_{m+1,1}^{-1}\bXdiamd_{1:m+1,1}
=\bXdiamd_{1:m,1}^T(\bQ_{m1}-\TOPRIGHT\BOTRIGHT^{-1}\TOPRIGHT)^{-1}\bXdiamd_{1:m,1}\\
&&\qquad\qquad -\bXdiamd_{1:m,1}^T(\bQ_{m1}-\TOPRIGHT\BOTRIGHT^{-1}\TOPRIGHT)^{-1}
\TOPRIGHT\BOTRIGHT^{-1}\bXdiamd_{m+1,1}\\
&&\qquad\qquad -\bXdiamd_{m+1,1}^T\BOTRIGHT^{-1}\TOPRIGHT^T
(\bQ_{m1}-\TOPRIGHT\BOTRIGHT^{-1}\TOPRIGHT)^{-1}\bXdiamd_{1:m,1}\\[1ex]
&&\qquad\qquad +\bXdiamd_{m+1,1}^T(\BOTRIGHT
-\TOPRIGHT^T\bQ_{m1}^{-1}\TOPRIGHT)^{-1}
\bXdiamd_{m+1,1}\\[1ex]
&&\quad=\bXdiamd_{1:m,1}^T\Big[\bQ_{m1}-\bX_{1:m,1}\bM'\bX_{m+1,1}^T
\big\{\bX_{m+1,1}(\bM+\bM')\bX^T_{m+1,1}
+\lambda\bI\big\}^{-1}\bX_{m+1,1}\bM'\bX_{1:m,1}^T\Big]^{-1}\bXdiamd_{1:m,1}\\[1ex]
&&\quad\quad -\bXdiamd_{1:m,1}^T \Big[\bQ_{m1}-\bX_{1:m,1}\bM'\bX_{m+1,1}^T
\big\{\bX_{m+1,1}(\bM+\bM')\bX^T_{m+1,1}
+\lambda\bI\big\}^{-1}\bX_{m+1,1}\bM'\bX_{1:m,1}^T\Big]^{-1}\\
&&\qquad\quad\times\bX_{1:m,1}\bM'\bX_{m+1,1}^T\big\{\bX_{m+1,1}(\bM+\bM')\bX^T_{m+1,1}
+\lambda\bI\big\}^{-1}\bXdiamd_{m+1,1}\\[1ex]
&&\quad\quad -\bXdiamd_{m+1,1}^T 
\big\{\bX_{m+1,1}(\bM+\bM')\bX^T_{m+1,1}
+\lambda\bI\big\}^{-1}\bX_{m+1,1}\bM'\bX_{1:m,1}^T\\
&&\qquad\qquad\times\Big[\bQ_{m1}-\bX_{1:m,1}\bM'\bX_{m+1,1}^T\big\{\bX_{m+1,1}(\bM+\bM')\bX^T_{m+1,1}
+\lambda\bI\big\}^{-1}\bX_{m+1,1}\bM'\bX_{1:m,1}^T\Big]^{-1}\bXdiamd_{1:m,1}\\[1ex]
&&\quad\quad+\bXdiamd_{m+1,1}^T
\big\{\bX_{m+1,1}(\bM+\bM')\bX_{m+1,1}^T-\bX_{m+1,1}\bM'\bX_{1:m,1}^T\bQ_{m1}^{-1}
\bX_{1:m,1}\bM'\bX_{m+1,1}^T+\lambda\bI\big\}^{-1}
\bXdiamd_{m+1,1}\\[2ex]
&&\quad=\gothicT_1-\gothicT_2-\gothicT_2^T+\gothicT_3+\gothicT_4
\end{eqnarray*}
}
where
{\setlength\arraycolsep{1pt}
\begin{eqnarray*}
\gothicT_1&\equiv&\bXdiamd_{1:m,1}^T\big(\bQ_{m1}+\bGamma_1\big)^{-1}
\bXdiamd_{1:m,1},\\[1ex]
\gothicT_2&=&\bXdiamd_{1:m,1}^T \big(\bQ_{m1}+\bGamma_1\big)^{-1}
\bX_{1:m,1}\bM'\bGamma_2
(\bX^T_{m+1,1}\bX_{m+1,1})^{-1}\bX^T_{m+1,1}\bXdiamd_{m+1,1}\\[1ex]
\gothicT_3&=&(1/\lambda)\bXdiamd_{m+1,1}^T
\big\{\bI-\bX_{m+1,1}(\bX_{m+1,1}^T\bX_{m+1,1})^{-1}\bX_{m+1,1}^T\big\}\bXdiamd_{m+1,1}
\\[1ex]
\gothicT_4&=&\bXdiamd_{m+1,1}^T\bX_{m+1,1}(\bX_{m+1,1}^T\bX_{m+1,1})^{-1}\\[1ex]
&&\qquad\times\big\{\bM+\bM'
-\bM'\bX_{1:m,1}^T\bQ_{m1}^{-1}\bX_{1:m,1}\bM'+\lambda(\bX_{1:m,1}^T\bX_{1:m,1})^{-1}\big\}^{-1}\\
&&\qquad\times(\bX_{m+1,1}^T\bX_{m+1,1})^{-1}\bX_{m+1,1}^T\bXdiamd_{m+1,1},\\[1ex]
\bGamma_1&\equiv&\bX_{1:m,1}\bM'\bGamma_2\bM'\bX_{1:m,1}^T
\quad\mbox{and}\quad\bGamma_2\equiv-\big\{\bM+\bM'
+\lambda(\bX_{m+1,1}^T\bX_{m+1,1})^{-1}\big\}^{-1}.
\end{eqnarray*}
}
Application of Woodbury's matrix identity (\ref{eq:WoodburyMatIdent})
to $\big(\bQ_{m1}+\bGamma_1\big)^{-1}$ with $\bS=\bQ_{m1}$, 
$\bU=\bX_{1:m,1}\bM'$, $\bV=\bM'\bX_{1:m,1}^T$ and $\bT=\bGamma_2$ leads to 
{\setlength\arraycolsep{1pt}
\begin{eqnarray*}
\gothicT_1&=&\bXdiamd_{1:m,1}^T\bQ_{m1}^{-1}\bXdiamd_{1:m,1}\\[1ex]
&&\qquad-\bXdiamd_{1:m,1}^T\bQ_{m1}^{-1}\bX_{1:m,1}\bM'\big\{\bGamma_2
+\bM' \bX_{1:m,1}^T\bQ_{m1}^{-1}\bX_{1:m,1}\bM'\big\}^{-1}\bM'\bX_{1:m,1}^T\bQ_{m1}^{-1}\bXdiamd_{1:m,1}\\[1ex]
&=&\frac{n_{11}+\ldots+n_{m1}}{\lambda}
\Big[\mbox{lower right $\ddiamd\times\ddiamd$ block of}\
\big\{E\big([\bXrv\  \bXTdiamdrv]^{\otimes2}\big)\big\}^{-1}\Big]^{-1}\{1+o_P(1)\}
\end{eqnarray*}
}
by Lemma \ref{lem:XBpart} and the inductive hypothesis. Similarly, the first three factors
of $\gothicT_2$ are
{\setlength\arraycolsep{1pt}
\begin{eqnarray*}
\bXdiamd_{1:m,1}^T \big(\bQ_{m1}+\bGamma_1\big)^{-1}\bX_{1:m,1}
&=&\bXdiamd_{1:m,1}^T\bQ_{m1}^{-1}\bX_{1:m,1}
-\bXdiamd_{1:m,1}^T\bQ_{m1}^{-1}\bX_{1:m,1}\\[1ex]
&&\ \times\bM'\big\{\bGamma_2
+\bM'\bX_{1:m,1}^T\bQ_{m1}^{-1}\bX_{1:m,1}\bM'\big\}^{-1}
\bM'\bX_{1:m,1}^T\bQ_{m1}^{-1}\bX_{1:m,1}\\[1ex]
\end{eqnarray*}
}
which soon leads to $\gothicT_2$ having all entries being $O_P(m)$.
Next, we have 
$$\gothicT_3=\big(n_{m+1,1}/\lambda\big)
\Big[\mbox{lower right $\ddiamd\times\ddiamd$ block of}\
\big\{E\big([\bXrv\  \bXTdiamdrv]^{\otimes2}\big)\big\}^{-1}\Big]^{-1}\{1+o_P(1)\}$$
and $\gothicT_4$ having all entries being $O_P(1)$. Combining these results
for $\gothicT_1$, $\gothicT_2$, $\gothicT_3$ and $\gothicT_4$ leads to
$$\left(\sum_{i=1}^mn_{i1}\right)^{-1}
\bXdiamd_{1:m+1,1}^T\bQ_{m+1,1}^{-1}\bXdiamd_{1:m+1,1}
\convprob(1/\lambda)\Big[\mbox{lower right $\ddiamd\times\ddiamd$ block of}\
\big\{E\big([\bXrv\  \bXTdiamdrv]^{\otimes2}\big)\big\}^{-1}\Big]^{-1}
$$
which proves Lemma \ref{lem:XBpart} (a) for all $m\in\naturalNumbers$ and $m'=1$.
The proof of Lemma \ref{lem:XBpart} (b) for all $m\in\naturalNumbers$ and $m'=1$
involves a similar set of arguments.

Now we turn our attention to establishing Lemma \ref{lem:XBpart} (a) for $m=1$ and $m'=2$.
Noting that
$$
\bQ_{12}
=
\left[
\begin{array}{cc}
\bX_{11} & \bO \\[1ex]
\bO      & \bX_{12}
\end{array}
\right]\left[
\begin{array}{cc}
\bM+\bM' &  \bM \\[2ex]
\bM      & \bM+\bM'
\end{array}
\right]
\left[
\begin{array}{cc}
\bX_{11} & \bO \\[1ex]
\bO      & \bX_{12}
\end{array}
\right]^T+\lambda\bI.
$$
and
$$\bXdiamd=\left[
\begin{array}{c}
\bXdiamd_{11}\\[1ex]
\bXdiamd_{12}
\end{array}
\right]
=
\left[
\begin{array}{cc}
\bXdiamd_{11} & \bO \\[1ex]
\bO      & \bXdiamd_{12}
\end{array}
\right]
\left[
\begin{array}{c}
\bI_d\\[1ex]
\bI_d
\end{array}
\right]
$$
we have
{\setlength\arraycolsep{1pt}
\begin{eqnarray*}
\bXdiamd^T\bQ_{12}^{-1}\bXdiamd&=&
\left[
\begin{array}{c}
\bIddiamnd\\[1ex]
\bIddiamnd
\end{array}
\right]^T
\left[
\begin{array}{cc}
\bXdiamd_{11} & \bO \\[1ex]
\bO      & \bXdiamd_{12}
\end{array}
\right]^T
\left\{
\left[
\begin{array}{cc}
\bX_{11} & \bO \\[1ex]
\bO      & \bX_{12}
\end{array}
\right]
\left[
\begin{array}{cc}
\bM+\bM' &  \bM \\[2ex]
\bM      & \bM+\bM'
\end{array}
\right]
\left[
\begin{array}{cc}
\bX_{11} & \bO \\[1ex]
\bO      & \bX_{12}
\end{array}
\right]^T
+\lambda\bI
\right\}^{-1}\\
&&\qquad\qquad\qquad\times\left[
\begin{array}{cc}
\bXdiamd_{11} & \bO \\[1ex]
\bO      & \bXdiamd_{12}
\end{array}
\right]
\left[
\begin{array}{c}
\bIddiamnd\\[1ex]
\bIddiamnd
\end{array}
\right]
\\[2ex]
&=&\gothicT_5+\gothicT_6
\end{eqnarray*}
}
where
{\setlength\arraycolsep{1pt}
\begin{eqnarray*}
\gothicT_5
&=&(1/\lambda)\left[
\begin{array}{c}
\bIddiamnd\\[1ex]
\bIddiamnd
\end{array}
\right]^T
\left[
\begin{array}{cc}
\bXdiamd_{11} & \bO \\[1ex]
\bO      & \bXdiamd_{12}
\end{array}
\right]^T\\[1ex]
&&\ \times\left\{\bI-\left[
\begin{array}{cc}
\bX_{11} & \bO \\[1ex]
\bO      & \bX_{12}
\end{array}
\right]
\left(\left[
\begin{array}{cc}
\bX_{11} & \bO \\[1ex]
\bO      & \bX_{12}
\end{array}
\right]^T
\left[
\begin{array}{cc}
\bX_{11} & \bO \\[1ex]
\bO      & \bX_{12}
\end{array}
\right]
\right)^{-1}
\left[
\begin{array}{cc}
\bX_{11} & \bO \\[1ex]
\bO      & \bX_{12}
\end{array}
\right]^T
\right\}
\left[
\begin{array}{cc}
\bXdiamd_{11} & \bO \\[1ex]
\bO      & \bXdiamd_{12}
\end{array}
\right]
\left[
\begin{array}{c}
\bIddiamnd\\[1ex]
\bIddiamnd
\end{array}
\right]\\[1ex]
&=&\frac{n_{11}}{\lambda}
\left\{\left(\frac{1}{n_{11}}\bXdiamd_{11}^T \bXdiamd_{11}\right)   
-\left(\frac{1}{n_{11}}\bXdiamd_{11}^T \bX_{11}\right)
\left(\frac{1}{n_{11}}\bX_{11}^T \bX_{11}\right)^{-1}
\left(\frac{1}{n_{11}}\bX_{11}^T\bXdiamd_{11}\right)\right\}\\[1ex]
&&\ +\frac{n_{12}}{\lambda}
\left\{\left(\frac{1}{n_{12}}\bXdiamd_{12}^T \bXdiamd_{12}\right)   
-\left(\frac{1}{n_{12}}\bXdiamd_{11}^T \bX_{12}\right)
\left(\frac{1}{n_{12}}\bX_{12}^T \bX_{12}\right)^{-1}
\left(\frac{1}{n_{12}}\bX_{12}^T\bXdiamd_{12}\right)\right\}\\[1ex]
&=&\frac{n_{11}+n_{12}}{\lambda}
\left[\mbox{lower right $\ddiamd\times\ddiamd$ block of}\
\big\{E\big([\bXrv\  \bXTdiamdrv]^{\otimes2}\big)\big\}^{-1}\right]^{-1}
\{1+o_P(1)\}.
\end{eqnarray*}
}
and
{\setlength\arraycolsep{1pt}
\begin{eqnarray*}
\gothicT_6
&=&
\left[
\begin{array}{c}
\big(\bXdiamd_{11}^T \bX_{11}\big)\big(\bX_{11}\bX_{11}\big)^{-1}\\[1ex]
\big(\bXdiamd_{12}^T \bX_{12}\big)(\bX_{12}\bX_{12}\big)^{-1}
\end{array}
\right]
\left[
\begin{array}{cc}
\bM+\bM'+\lambda(\bX_{11}^T\bX_{11})^{-1} & \bM\\[1ex]
\bM   & \bM+\bM'+\lambda(\bX_{12}^T\bX_{12})^{-1} 
\end{array}
\right]^{-1}\\[1ex]
&&\quad\times\left[
\begin{array}{c}
\big(\bXdiamd_{11}^T \bX_{11}\big)\big(\bX_{11}\bX_{11}\big)^{-1}\\[1ex]
\big(\bXdiamd_{12}^T \bX_{12}\big)(\bX_{12}\bX_{12}\big)^{-1}
\end{array}
\right].
\end{eqnarray*}
}
Since each of the entries of $\gothicT_6$ are $O_P(1)$ we have 
$$\frac{1}{n_{11}+n_{12}}\bXdiamd^T\bQ_{12}^{-1}\bXdiamd
\convprob(1/\lambda)
\left[\mbox{lower right $\ddiamd\times\ddiamd$ block of}\
\big\{E\big([\bXrv\  \bXTdiamdrv]^{\otimes2}\big)\big\}^{-1}\right]^{-1}
$$
which verifies Lemma \ref{lem:XBpart}(a) for $m=1$ and $m'=2$. An analogous pattern
continues for higher $m$ and $m'$ which leads to the Lemma \ref{lem:XBpart}(a)
result holding generally.

For Lemma \ref{lem:XBpart}(b) in the $m=1$ and $m'=2$ case we instead have,
using Corollary 2.1(a) and Lemma \ref{lem:ApsBinvResIden},
{\setlength\arraycolsep{1pt}
\begin{eqnarray*}
\bX^T\bQ_{12}^{-1}\bXdiamd&=&
\left[
\begin{array}{c}
\bI_d\\[1ex]
\bI_d
\end{array}
\right]^T
\left[
\begin{array}{cc}
\bX_{11} & \bO \\[1ex]
\bO      & \bX_{12}
\end{array}
\right]^T
\left\{
\left[
\begin{array}{cc}
\bX_{11} & \bO \\[1ex]
\bO      & \bX_{12}
\end{array}
\right]
\left[
\begin{array}{cc}
\bM+\bM' &  \bM \\[2ex]
\bM      & \bM+\bM'
\end{array}
\right]
\left[
\begin{array}{cc}
\bX_{11} & \bO \\[1ex]
\bO      & \bX_{12}
\end{array}
\right]^T
+\lambda\bI
\right\}^{-1}\\
&&\qquad\qquad\qquad\times\left[
\begin{array}{cc}
\bXdiamd_{11} & \bO \\[1ex]
\bO      & \bXdiamd_{12}
\end{array}
\right]
\left[
\begin{array}{c}
\bIddiamnd\\[1ex]
\bIddiamnd
\end{array}
\right]
\\[2ex]
&=&\left[
\begin{array}{c}
\bI_d\\[1ex]
\bI_d
\end{array}
\right]^T
\left[
\begin{array}{cc}
\bM+\bM'+\lambda(\bX_{11}^T\bX_{11})^{-1} &  \bM \\[2ex]
\bM      & \bM+\bM'+\lambda(\bX_{12}^T\bX_{12})^{-1} 
\end{array}
\right]^{-1}\\[1ex]
&&\qquad\times
\left[
\begin{array}{c}
\left(\frac{1}{n_{11}}\bX_{11}^T \bX_{11}\right)^{-1}
\left(\frac{1}{n_{11}}\bX_{11}^T\bXdiamd_{11}\right)\\[1ex]
\left(\frac{1}{n_{12}}\bX_{11}^T \bX_{12}\right)^{-1}
\left(\frac{1}{n_{12}}\bX_{12}^T\bXdiamd_{12}\right)
\end{array}
\right]\\[1ex]
&\convprob&
\left[
\begin{array}{c}
\bI_d\\[1ex]
\bI_d
\end{array}
\right]^T
\left[
\begin{array}{cc}
\bM+\bM' &  \bM \\[2ex]
\bM      & \bM+\bM'
\end{array}
\right]^{-1}
\left[
\begin{array}{c}
\bI_d\\[1ex]
\bI_d
\end{array}
\right]\big\{E(\bXrv^{\otimes2})\big\}^{-1}E(\bXrv\bXdiamdrv^T)\\[1ex]
&=&(\bM+\smhalf\bM')^{-1}\big\{E(\bXrv^{\otimes2})\big\}^{-1}E(\bXrv\bXdiamdrv^T)
\end{eqnarray*}
}
which verifies Lemma \ref{lem:XBpart}(b) for $m=1$ and $m'=2$.

For general $m$ and $m'$, note that the behavior of $\bX^T\bQmmd^{-1}\bXdiamd$
mimics that of the $\bX^T\bQmmd^{-1}\bX$ special case, with the 
$\big\{E(\bXrv^{\otimes2})\big\}^{-1}E(\bXrv\bXdiamdrv^T)$ factor being the only 
difference in the convergence in probability limit. The summations that
provide the Lemma \ref{lem:XApart}(a) result have analogous behaviors in
this extended case and lead to Lemma \ref{lem:XBpart}(b) holding generally.

\section*{Additional References}

\bib
Hall, P. \myand Heyde, C.C. (1980). \textit{Martingale Central Limit
Theory and Its Application}. New York: Academic Press.

\bib
Jiang, J., Jiang, W., Paul, D., Zhang, Y. \myand Zhao, H. (2023).
High-dimensional asymptotic behaviour of inference based on GWAS
summary statistic. \textit{Statistica Sinica}, \textbf{33}, 1555--1576.

\bib
Magnus, J.R. \myand Neudecker, H. (1999). \textit{Matrix Differential
Calculus with Applications in Statistics and Econometrics, Revised 
Edition}. Chichester U.K.: Wiley.

\bib
Wand, M.P. (2002). Vector differential calculus in statistics.
\textit{The American Statistician}, \textbf{56}, 55--62.

\bib
Woodbury, M.A. (1950). Inverting modified matrices.
Statistical Research Group, Memorandum Report Number 42,
Princeton University, Princeton, New Jersey, U.S.A.

\end{document}